%% file: main_new.tex
\documentclass[12pt]{article}
\usepackage{amsmath}
\usepackage{amssymb}
\usepackage{graphicx}
\usepackage{enumerate}
\usepackage{url} 
\usepackage{todonotes}
\usepackage{bbold} 

\newtheorem{proposition}{Proposition}
\newtheorem{lemma}{Lemma}
\newtheorem{corollary}{Corollary}
\newtheorem{remark}{Remark}
\newtheorem{algorithm}{Algorithm}
\newtheorem{definition}{Definition}

\usepackage{caption}
\begin{document}

\def\spacingset#1{\renewcommand{\baselinestretch}%
{#1}\small\normalsize} \spacingset{1}

\title{A class of count time series models uniting compound Poisson INAR and INGARCH models}
\author{Johannes Bracher$^{1, 2}$ and Barbora Sobolov\'a$^1$\\
$^1$ Institute of of Statistics, Karlsruhe Institute of Technology, \\ 
$^2$ Computational Statistics Group, Heidelberg Institute for Theoretical Studies}

\newcommand{\juv}{E}
\newcommand{\clustersize}{\theta} 
\newcommand{\disp}{\psi} 
\newcommand{\import}{\varepsilon}
\newcommand{\Reff}{R_\text{e}}
\newcommand{\red}[1]{\textcolor{red}{#1}}

\maketitle

\begin{abstract}
INAR (integer-valued autoregressive) and INGARCH (integer-valued GARCH) models are among the most commonly employed approaches for count time series modelling, but have been studied in largely distinct strands of literature. In this paper, a new class of generalized integer-valued ARMA (GINARMA) models is introduced which unifies a large number of compound Poisson INAR and INGARCH processes. Its stochastic properties, including stationarity and geometric ergodicity, are studied. Particular attention is given to a generalization of the INAR($p$) model which parallels the extension of the INARCH($p$) to the INGARCH($p$, $q$) model. For inference, we consider moment-based estimation and a maximum likelihood inference scheme inspired by the forward algorithm. Models from the proposed class have a natural interpretation as stochastic epidemic processes, which throughout the article is used to illustrate our arguments. In a case study, different instances of the class, including both established and newly introduced models, are applied to weekly case numbers of measles and mumps in Bavaria, Germany. 
\end{abstract}

\bigskip

\begin{center}
\end{center}


\maketitle

\newpage


\section{Introduction}
\label{sec:intro}

Count time series arise in many contexts from hydrology \cite{McKenzie1985} to criminology and traffic studies \cite{Scotto2015}. Numerous modelling approaches for such data exist, including e.g., hidden Markov \cite{Zucchini2009}, generalized linear ARMA \cite{Benjamin2003} and latent Gaussian models \cite{Jia2023}. This diversity led a recent review \cite{Davis2021} to conclude that ``the field developed without a unifying theory''. In the present paper we aim to provide an overarching framework for two particularly influential model classes, namely the INAR (integer-valued autoreggressive) and INGARCH (integer-valued GARCH) classes. These have been highlighted as ``probably the most widely used approaches for stationary count time series'' \cite{Weiss2021}. While INAR models employ thinning operations and resemble branching processes \cite{Dion1995}, INGARCH models take their starting point in generalized linear regression. Despite some known links between the two \cite{Lu2021, Weiss2015}, they have been treated in largely distinct strands of literature. Our contribution to bridging this gap is threefold:

\begin{itemize}
\item Building on generalized INAR models \cite{Latour1998}, we define a broad model class comprising many well-known INAR and INGARCH processes, as well as new models. Its properties are studied with a particular focus on compound Poisson (CP) formulations.
\item As an important special case, we study a generalization of the INAR($p$) model which parallels the extension of the INARCH($p$) to the INGARCH($p, q$).
\item Borrowing ideas from epidemic modelling, we contrast the ``mechanistic'' assumptions of different instances of the class, thus providing a useful language to distinguish them.
\end{itemize}
Indeed, both INAR \cite{Cardinal1999, Pedeli2015} and INGARCH \cite{Bracher2021, Ferland2006} models are commonly applied to infectious disease counts, though often without discussion of the implied assumptions on disease spread (see \cite{Bauer2018} for an exception). Throughout the article we will use the epidemiological interpretation to strengthen intuition, and we will conclude with a case study on measles and mumps in the German state of Bavaria. Similarly to \cite{Kucharski2014}, we will estimate local effective reproductive numbers and the relative importance of imported cases.

The article is structured as follows. In Section \ref{sec:preliminaries}, we provide some background on CP-INAR and INGARCH models. In Section \ref{sec:ginarma} we introduce our general model class, before turning to its INAR-like instance and inference aspects in Section \ref{sec:extension_inar}. In Section \ref{sec:real_data}, the real-data application is presented before Section \ref{sec:discussion} concludes with a discussion.

\section{Preliminaries}
\label{sec:preliminaries}

We start by reviewing relevant fundamentals of (generalized) INAR and INGARCH models.

\subsection{Poisson (G)INAR(1) and INGARCH(1,1) models}
\label{subsec:poisson_ginar_inarch}

The generalized INAR(1) model \cite{Latour1998}, GINAR(1) for short, is defined as $\{X_t, t \in \mathbb{Z}\}$ with
\begin{equation}
X_t = \alpha \bullet X_{t - 1} + \import_t \label{eq:ginar1}    
\end{equation}
and $\alpha > 0$. The \textit{imports} $\{\import_t\}$ are independent and identically distributed (i.i.d.) count random variables with mean $\nu > 0$ and variance $\sigma^2_\nu < \infty$, while $\bullet$ is the \textit{generalized thinning operator}. With $N \in \mathbb{N}_0$ and $\alpha > 0$ it is defined as $\alpha \bullet N = 0$ if $N = 0$ and 
\begin{equation}\label{eq:def_thinning}
    \alpha \bullet N = \sum_{i = 1}^N Z_i
\end{equation}
otherwise. Independently of $N$, the $Z_1, \dots, Z_N$ are i.i.d. draws from a count-valued \textit{offspring distribution} with mean $\alpha > 0$ and variance $\sigma^2_\alpha$. All thinnings in \eqref{eq:ginar1} are performed independently of each other, the imports $\{\import_t\}$ and the past of the process $\{X_t\}$, an assumption we will make throughout the paper unless relaxed explicity.

Model \eqref{eq:ginar1} can be read as an adaptation of the classic Gaussian AR(1) process with multiplication replaced by generalized thinning. As GINAR(1) models are first-order conditionally linear autoregressive (CLAR) models \cite{Grunwald2000}, they preserve many stochastic properties of their continuous counterpart. Two particularly influential instances of the class are the Poisson INAR(1) \cite{McKenzie1985, Al-Osh1987} and INARCH(1) \cite{Ferland2006, Fokianos2009} models (the naming of the latter being somewhat controversial, see Remark 4.1.2 in \cite{Weiss2018}). While in both $\import_t \stackrel{\text{i.i.d.}}{\sim} \text{Pois}(\nu)$ is assumed, the offspring distributions differ. In the INAR(1), given by
\begin{equation}
X_t = \alpha \circ X_{t - 1} + \epsilon_t,    
\end{equation}
binomial thinning $\circ$ \cite{Steutel1979} is used, which results from $Z_i \stackrel{\text{i.i.d.}}{\sim} \text{Bern}(\alpha)$. For the INARCH(1), Poisson thinning $\star$  with $Z_i \stackrel{\text{i.i.d.}}{\sim} \text{Pois}(\alpha)$ is assumed instead.

\begin{remark}
    GINAR(1) models can be thought of as simple epidemic processes \cite{Cardinal1999}. Each of $X_t$ infectives present in a population at time $t$ causes on average $\alpha$ new infectives (``offspring'') at $t + 1$ and then recovers. Infections from outside sources are imported at rate $\nu$. As will be seen in Section \ref{subsec:interpretation_epidemic}, similar interpretations also apply to extended models.
\end{remark}

\noindent While the GINAR(1) representation of the Poisson INARCH(1) is well-known (e.g., \cite[p.56]{Weiss2018}), the model is usually defined in terms of a conditional generalized linear regression model (GLM). It then becomes $\{X_t, t \in \mathbb{N}\}$ with
\begin{align}
X_t \mid X_{t - 1}, \dots, X_0, \lambda_0 & \sim \text{Pois}(\lambda_t) \label{eq:inarch1}\\
\lambda_t & = \nu + \alpha X_{t - 1}\label{eq:lambda_inarch1}
\end{align}
and fixed starting values $\lambda_0, X_0$. This formulation is attractive as it can be extended to the Poisson INGARCH(1, 1) model \cite{Ferland2006, Fokianos2009}, where \eqref{eq:lambda_inarch1} becomes
\begin{equation}
\lambda_t = \nu + \alpha X_{t - 1} + \beta \lambda_{t - 1}    \label{eq:ingarch11} 
\end{equation}
with $0 \leq \beta$. The \textit{feedback term} $\beta \lambda_{t - 1}$ here leads to an ARMA(1, 1) autocorrelation function.

\subsection{Compound Poisson distributions}

To handle overdispersion in a flexible way, Poisson (G)INAR and INGARCH models are commonly extended using compound Poisson (CP) distributions. A random variable $Y$ is said to follow a CP distribution \cite[Chapter 3]{Feller1968} if it can be written as a randomly stopped sum $
Y = \sum_{i = 1}^N Z_i$, where $N$ is Poisson distributed and $Z_1, \dots, Z_N \stackrel{\textnormal{i.i.d.}}{\sim} G(\clustersize)$ independently of $N$. We assume throughout that the \textit{cluster distribution} $G(\clustersize)$ has a single parameter $\clustersize$ and support $\{1, \dots, \zeta\}$, where $\zeta$ is the \textit{order} of the CP distribution and $\zeta = \infty$ is allowed. For simplicity we identify $\theta$ with the mean of $G$ and in analogy to \eqref{eq:def_thinning} use the shorthand
\begin{equation}
\clustersize * N = \sum_{i = 1}^N Z_i, \ \ \ Z_1, \dots, Z_N \stackrel{\textnormal{i.i.d.}}{\sim} G(\clustersize). \label{eq:def_compund_thinning}
\end{equation}
We assume $0 < \clustersize < \infty$ and denote the variance of $G(\theta)$ by $\sigma^2_\clustersize < \infty$. Setting $N \sim \text{Pois}(\mu/\clustersize)$, we obtain a CP distribution with mean $\mathbb{E}(Y) = \mu$ and variance $\text{Var}(Y) = \mu \times(\sigma^2_\clustersize/\clustersize + \clustersize)$.
\begin{remark}\label{remark:clustering}
We use the term ``\textit{cluster distribution}'' as in disease modelling and ecology, CP distributions are often applied to phenomena that occur in clusters. The number $N$ of clusters is then Poissonian, while the number of units per cluster follows $G(\clustersize)$.
\end{remark}
Two popular CP distributions are the Hermite and negative binomial.  Hermite imports are attractive in INAR models as the resulting marginal distributions are often also Hermite \cite{Fernandez-Fontelo2017, Weiss2015b}. A random variable $Y$ is Hermite distributed if it can be written as ${Y = A_1 + 2A_2}$ where independently $A_1 \sim \text{Pois}(\lambda_1), A_2 \sim \text{Pois}(\lambda_2)$. In slight variation of \cite{Gupta1974} we parameterize the distribution by its mean $\mu = \lambda_1 + 2\lambda_2$ and a dispersion parameter $\disp = 2\lambda_2/(\lambda_1 + 2\lambda_2) \in [0, 1]$, implying $\text{Var}(Y) = (1 + \disp)\mu$. The probability mass function is then
$$
\text{Pr}(Y = y) = \exp\left[\mu\left(-1 + \frac{\disp}{2}\right)\right] \mu^y(1 - \disp)^y \sum_{j = 0}^{[y/2]} \frac{\disp^j}{2^j\mu^j(1 - \disp)^{2j}(y - 2j)!j!}, y = 0, 1, 2, \dots
$$
where $[y/2]$ is the integer part of $y/2$. The Hermite is a CP distribution of order 2, i.e., the cluster distribution has support $\{1, 2\}$ (see Supplement \ref{suppl:param_hermite}).

For the negative binomial distribution we likewise use a parameterization via its mean $\mu$ and a dispersion parameter $\disp > 0$, given by
$$
\text{Pr}(Y = y) = \frac{\Gamma(1/\disp + y)}{y!\Gamma(1/\disp)} \left(\frac{1/\disp}{1/\disp + \mu}\right)^{1/\disp} \left(\frac{\mu}{1/\disp + \mu}\right)^y,
$$
while $\text{Var}(Y) = (1 + \disp\mu)\mu$. The negative binomial distribution is equivalent to a CP distribution with a logarithmic cluster distribution \cite{Weiss2018}, see Supplement \ref{suppl:negbin_compound}.

\subsection{CP-(G)INAR($p$) and INGARCH($p, q$) models}

In the (generalized) INAR framework,  the extension of model \eqref{eq:ginar1} to the CP case \cite{Schweer2014} is straightforward as only the import distribution is replaced by a CP with mean $\nu$ and variance $\sigma^2_\nu$. Hermite and negative binomial innovations have been considered e.g., by \cite{Fernandez-Fontelo2017} and \cite{Pedeli2015}. A higher-order GINAR($p$) model is obtained by setting \cite{Dion1995}
\begin{equation}
    X_t = \sum_{i = 1}^p \alpha_i \bullet X_{t - i} \ + \ \import_t. \label{eq:ginar_p}
\end{equation}
Slightly generalizing \cite{Latour1998}, we allow for dependent offspring $(\alpha_1 \bullet X_t, \dots, \alpha_p \bullet X_t)$, see Definition \ref{def:ingarma}. This is because we will extend the INAR($p$) by Alzaid and Al-Osh \cite{Alzaid1990}, where
\begin{equation}
X_t = \sum_{i = 1}^p \alpha_i \circ X_{t - i} + \import_t\label{eq:inar_p}
\end{equation}
with $\alpha_1, \dots, \alpha_p \geq 0, 0 < \sum_{i = 1}^p \alpha_i < 1$ is combined with multinomial thinnings,
\begin{equation}
(\alpha_1 \circ X_t, \dots, \alpha_p \circ X_t) \sim \text{Mult}(X_t, \alpha_1, \dots, \alpha_p).\label{eq:multinomial_thinning}
\end{equation}
We note that an equally well-known INAR($p$) model with independent thinning operations has been proposed by Du and Li \cite{Du1991}, but it is less fruitful within our framework.

The GINAR($p$) model can be generalized further to the $\text{GINARMA}_{\text{DGL}}(p, q)$, given by
\begin{align}\label{eq:dion_inarma}
        X_t = \sum_{i = 1}^p \alpha_i \bullet X_{t - i} + \sum_{j = 1}^q \delta_j \bullet \import_{t - j} \ \ + \ \ \import_t.
\end{align}
We here add the initials of its authors -- Dion, Gauthier and Latour -- to the notation to distinguish this INARMA model from our own suggestion presented later on. 

To extend \eqref{eq:ingarch11} to a CP-INGARCH($p, q$) model we adopt notation from Wei\ss\ et al \cite[Sec. 2]{Weiss2017}. The model is then defined as $\{X_t, t \in \mathbb{N}\}$ with
\begin{align}
N_t \ & \mid \ X_{t - 1}, \dots, X_{1 - p}, \lambda_0, \dots, \lambda_{1 - q} \sim \text{Pois}(\lambda_t/\clustersize) \label{eq:N_CP_original}\\
X_t & = \sum_{i = 1}^{N_t} Z_{t, i} \ \ \text{where} \ \  Z_{t, 1}, \dots, Z_{t, N_t} \stackrel{\textnormal{i.i.d.}}{\sim} G(\clustersize)
\label{eq:X_CP_original}\\
\lambda_t & = \nu \ + \ \sum_{i = 1}^p \alpha_i X_{t - i} \ + \ \sum_{j = 1}^q \beta_j \lambda_{t - j}.\label{eq:lambda_CP_original}
\end{align}
Here, $\lambda_{1 - q}, \dots, \lambda_0 \geq 0$ and $X_{1 - p}, \dots, X_0 \in \mathbb{N}_0$ are again fixed and we assume $\alpha_1, \dots, \alpha_p,$ $\beta_1, \dots, \beta_q \geq 0, \sum_{i = 1}^p \alpha_i > 0$. Given the past, $X_t$ then follows a CP distribution with mean $\lambda_t$ and variance $\lambda_t \times (\sigma^2_\theta/\theta + \theta)$. This is more restrictive than in \cite{Goncalves2015} where $\clustersize_t$ is a function of $\lambda_t$. It nonetheless contains e.g., the negative binomial (\cite{Weiss2018, Xu2012}), generalized Poisson \cite{Xu2012, Zhu2012} and Neyman Type A \cite{Goncalves2015a} INGARCH models. For the Hermite and negative binomial cases we provide details in Supplements \ref{suppl:hermite_ingarch} and \ref{suppl:negbin_ingarch}. We note that despite its name, the INGARCH($p, q$) has an ARMA($p, q$) correlation structure \cite[Remark 4.1.3]{Weiss2018}; nonetheless it behaves quite differently than the $\text{GINARMA}_{\text{DGL}}(p, q)$, see Section \ref{subsubsec:general_imports}.

\section{A new GINARMA($p, q$) model}
\label{sec:ginarma}

\subsection{Model definition}

We now propose an alternative GINARMA extension of model \eqref{eq:ginar_p}. Rather than directly replacing the multiplications in the Gaussian ARMA($p, q$) model
\begin{align}
        X_t \ = \ \sum_{i = 1}^p \alpha_i X_{t - i} \ \ + \ \ \sum_{j = 1}^q \delta_j \import_{t - j} \ \ + \ \ \nu \ \ + \ \ \import_t, \ \ \ \import_t \stackrel{\textnormal{i.i.d.}}{\sim} \text{N}(0, \sigma^2_\import)\label{eq:classic_arma}
\end{align}
by thinnings as in \eqref{eq:dion_inarma}, we use the following reformulation. Setting
$$
\kappa_i = \frac{\alpha_i + \delta_i}{1 + \sum_{j = 1}^q \delta_j}, \ \ \beta_j = -\delta_j, \ \ \tau = \frac{\nu}{1 + \sum_{j = 1}^q \delta_j},
$$
with $\delta_j = 0$ for $j > q$, an equivalent of \eqref{eq:classic_arma} is (see Supplement \ref{sec:arma})
\begin{align}
X_t & = \Bigl(1 - \sum_{j = 1}^q \beta_j \Bigr) \times E_t \ \ + \ \ \nu \ \ + \ \ \import_t. \label{eq:X_arma_reformulation} 
\end{align}
Here, $\{E_t\}$ is an auxiliary process defined as
\begin{align}
E_t & = \sum_{j = 1}^q \beta_j\times E_{t - j} \ \ + \ \ \sum_{i = 1}^p \kappa_i \times X_{t - i}.\label{eq:E_arma_reformulation}
\end{align}
Based on this structure, we define our model as follows. Note that it contains an additional compounding step, which is useful to accommodate CP-INGARCH models in the class.

\begin{definition}\label{def:ingarma}
The GINARMA($p, q$) model is a stochastic process $\{X_t, t \in \mathbb{N}\}$ with
\begin{align}
X_t & = \clustersize * \left[\Bigl(1 - \sum_{j = 1}^q \beta_j \Bigr) \circ E_t + \import_t\right]\label{eq:X_general}\\
E_t & = \sum_{j = 1}^q \beta_j \circ E_{t - j} \ \ + \ \ \sum_{i = 1}^p \kappa_i \bullet X_{t - i}\label{eq:E_general}
\end{align}
and $\kappa_i, \beta_j \geq 0; 0 < \sum_{i = 1}^q \kappa_i; 0 \leq \sum_{j = 1}^q \beta_j < 1$. Specifically, the following is assumed.
\begin{itemize}
    \item[(i)] The sequence $\{\import_t\}$ consists of i.i.d. realizations from an integer-valued import distribution with mean $0 < \tau < \infty$ and variance $\sigma^2_\tau < \infty$.
    \item[(ii)] The offspring $(\kappa_1 \bullet X_t, \dots, \kappa_p \bullet X_t)$ result from a generalized thinning operation $\bullet$ and can be dependent. Specifically, we assume $(\kappa_1 \bullet X_t, \dots, \kappa_p \bullet X_t) = \mathbf{Z}_{t, 1} + \dots + \mathbf{Z}_{t, X_t}$, where, independently of $X_t$, the vectors $\mathbf{Z}_{t, j}$ are i.i.d. on $\mathbb{N}_0^p$ with $\mathbb{E}(\mathbf{Z}_{t, j}) = (\kappa_1, \dots, \kappa_p)$ and finite variances $\sigma^2_{\kappa_1}, \dots, \sigma^2_{\kappa_p}$. To avoid dealing with e.g., purely even-valued offspring distributions we assume $\textnormal{Prob}(\kappa_i \bullet 1 = 1) > 0$ if $\kappa_i > 0$.
    \item[(iii)] The thinnings of $E_t$ are coupled via
    \end{itemize}
\begin{equation}
    \left[\beta_1 \circ E_t, \dots, \beta_q \circ E_t, \Bigl(1 - \sum_{j = 1}^q \beta_j \Bigr) \circ E_t\right] \biggm| E_t \sim \textnormal{Mult}\left(E_t; \beta_1, \dots, \beta_q, 1 - \sum_{j = 1}^q \beta_j\right).\label{eq:mult_general}
\end{equation}
\begin{itemize}
    \item[(iv)] As in equation \eqref{eq:def_compund_thinning}, $\theta \:*$ denotes a compounding step with a clustering distribution $G$. Its mean and variance are $0 < \theta < \infty$ and $\sigma^2_\theta < \infty$, respectively.
    \item[(v)] Apart from the dependencies between different thinnings of the same $X_t$ or $E_t$ introduced in (ii) and (iii), all thinnings and compoundings are performed independently of each other, the past of the process and the import sequence $\{\varepsilon_t\}$.
    \item[(vi)] Unless stated otherwise, the values $X_{1 - p}, \dots, X_0,$ $E_{1 - q}, \dots, E_0$ are fixed. In some instances we will initialize them with the respective stationary distributions instead.
\end{itemize}
\end{definition}

\noindent Assumption (iii) may seem arbitrary at first sight, but is central to obtaining appealing stochastic properties. Notably, it ensures that the model reduces to the $\text{GINAR}(p)$ if the compounding step $\theta \: *$ is omitted and $\beta_1 = \dots = \beta_q = 0$. As we will discuss in Section \ref{sec:alternative_formulation}, many INGARCH models can be obtained via Poisson offspring (setting $\bullet$ to Poisson thinning $\star$), while Bernoulli offspring (setting $\bullet$ to binomial thinning $\circ$) yield a new extension of the INAR class; see Section \ref{sec:extension_inar}.

\subsection{Interpretation as a stochastic epidemic process}
\label{subsec:interpretation_epidemic}

Formulation \eqref{eq:X_general}--\eqref{eq:mult_general} can be interpreted as a discrete-time model for the spread of an infectious disease, which provides a useful language and intuition for the following. We first illustrate this for $p = q = 1$ and omitting the compounding step, i.e., we consider
\begin{equation}
X_t = (1 - \beta) \circ E_t + \import_t; \ \ \ \ E_t = \beta \circ E_{t - 1} + \kappa \bullet X_{t - 1}. \label{eq:compounding_omitted}
\end{equation}
A graphical display of the following interpretation is provided in the top panel of Figure \ref{fig:ingarch_flowchart_poisson}.

\begin{enumerate}
\item $X_t$ is the number of infectious individuals at time $t$. These stay infectious for one time period and independently cause new infections with mean $\kappa$ and variance $\sigma^2_\kappa$.
\item Individuals newly infected at time $t$ do not necessarily become infectious already at $t + 1$. Instead, they enter into an ``\textit{exposed pool}'' $\juv_{t + 1}$.
\item At each time $t$, each of the $E_t$ exposed individuals can either remain in the exposed pool (with probability $\beta$) or advance to infectiousness (with probability $1 - \beta$).
\item An exposed individual from $E_t$ advancing to infectiousness becomes part of $X_t$.
\item At each time $t$, $\import_t$ individuals get infectious due to external sources.
\end{enumerate}
The reproductive number, i.e., mean number of new infections caused by one infected is $\Reff = \kappa$. The latent period, defined as the number of time points an infected spends in the exposed pool, is geometrically distributed with mean $1/(1 - \beta)$. The same holds for the generation time, i.e., time between the start of infectiousness of one individual and that of a second individual infected by the first.

When allowing $p, q > 1$, $X_t$ can be seen as the number of individuals becoming \textit{newly} infectious at time $t$. These are contagious over $p$ time steps, with $(\kappa_1, \dots, \kappa_p)$ the infectivity profile. In the exposed pool, individuals can ``move forward'' up to $q$ time periods at once, leading to more complex latent period distributions; see Figure \ref{fig:ingarch_flowchart_poisson}, middle panel.

When adding the compounding step $\theta \: *$ from equation \eqref{eq:X_general} to \eqref{eq:compounding_omitted}, $E_t$ and $\varepsilon_t$ can be thought of as clusters of exposed individuals, each containing a $G(\clustersize)$-distributed number of members (see Remark \ref{remark:clustering}). All members of a cluster turn infectious simultaneously. The effective reproductive number then becomes $R_e = \kappa\clustersize$; see Figure \ref{fig:ingarch_flowchart_poisson}, bottom panel. 

The above mechanisms resemble classic epidemic models like the SEIR (susceptible-exposed-infectious-removed), with the difference that immunity due to infection is ignored. See \cite{Bauer2018} for a similar argument on the INARCH(1) and \cite{Bjoernstad2002} for a related model accounting for immunity. This simplification is appropriate e.g., for vaccine-preventable diseases with high, but not complete vaccination levels in a population \cite{DeSerres2000}. In this situation, only minor outbreaks seeded by imported cases occur, which do not meaningfully reduce the number of remaining susceptibles. We will return to such a setting in our case study in Section \ref{sec:real_data}.

\subsection{Stochastic properties in the case $p = q = 1$}
\label{eq:stochastic_properties_general}

Various stochastic properties can be obtained by noting that if $p = q = 1$, the process $\{E_t\}$ is a Galton-Watson branching process with immigration. 

\begin{lemma}
\label{lemma:galton_watson}
The process $\{E_t\}$ from \eqref{eq:E_general} with $p = q = 1$ can be expressed as
\begin{equation}
\juv_t = \sum_{k = 1}^{\juv_{t - 1}} B_{t - 1, k} \ \ + \ \ \import^*_t.
\label{eq:galton_watson}
\end{equation}
Here we set $\import^*_t = \kappa\ \bullet (\clustersize * \import_{t - 1})$ and independently for each $k = 1, \dots, E_{t - 1}$
\begin{align}
B_{t - 1, k} & = \begin{cases}
1 & \text{with probability } \beta\\
\kappa \bullet (\clustersize * 1) & \text{with probability } 1 - \beta.
\label{eq:Z_t_i}
\end{cases}
\end{align}
\end{lemma}

\begin{corollary}
If $\kappa\theta < 1$, the Markov chain $\{E_t\}$ is moreover aperiodic and irreducible.
\label{corollary:aperiodic_irreducible}
\end{corollary}

\begin{proposition}
\label{proposition:finite_moments}
For $p = q = 1$, the processes $\{E_t\}$ and $\{X_t\}$ from equations \eqref{eq:X_general}--\eqref{eq:E_general} have unique limiting-stationary distributions if $\kappa\clustersize < 1$. The limiting-stationary moments are finite up to order $r$ if the same is true for $\{\import_t\}$, $\kappa \bullet 1, \clustersize * 1$. 
\end{proposition}

\begin{lemma}
\label{lemma:general_moments}
Given they exist, the limiting-stationary means and variances of $\{E_t\}$ and $\{X_t\}$ in a GINARMA(1, 1) process are
$$
\mu_E = \frac{\kappa\tau\clustersize}{1 - \beta - (1 - \beta)\kappa\clustersize}, \ \ \ \mu_X = \frac{\tau\clustersize}{1 - \kappa\clustersize},
$$
\begin{align*}
\sigma^2_E & = \frac{\sigma^2_\kappa \clustersize \tau + (\sigma^2_\clustersize\tau + \sigma^2_\tau \clustersize^2) \times \kappa^2 + \mu_E \times (1 - \beta) \times \{\beta (1 - \kappa\clustersize)^2 + \sigma^2_\kappa\clustersize + \sigma^2_\clustersize\kappa^2\}}{1 - \{\beta + (1 - \beta)\kappa\clustersize\}^2},\\
\sigma^2_X & = (1 - \beta)\mu_E \sigma^2_\clustersize + \clustersize^2(1 - \beta)\{\beta\mu_E + (1 - \beta)\sigma^2_E\} +
  \tau\sigma^2_\clustersize + \sigma^2_\tau\clustersize^2.
\end{align*}
The autocovariance functions of $\{E_t\}$ and $\{X_t\}$ are of AR(1) and ARMA(1, 1) type, as
\begin{align*}
\gamma_E(d) & = \{\beta + (1 - \beta)\kappa\clustersize\}^d \times \sigma^2_E\\
\gamma_X(d) & = \{\beta + (1 - \beta)\kappa\clustersize\}^{d - 1} \times (1 - \beta)\clustersize \times \{\clustersize\beta(1 - \beta)(\sigma^2_E - \mu_E) + \kappa\sigma^2_X\}.
\end{align*}
\end{lemma}

\noindent Combining arguments from Pakes \cite{Pakes1971} and Meitz and Saikkonen \cite{Meitz2008}, it can be shown that $\{\juv_t\}$ is geometrically ergodic under mild conditions, which translates to $\{X_t\}$.

\begin{proposition}
\label{proposition:geometric_ergodicity}
For $p = q = 1$, the joint process $\{(X_t, E_t)\}$ from \eqref{eq:X_general}--\eqref{eq:E_general} is geometrically ergodic if $\kappa\clustersize < 1$ and $\sigma^2_\tau, \sigma^2_\kappa, \sigma^2_\clustersize < \infty$. If the initial value $E_t$ is generated from its stationary distribution, the process is moreover $\beta$-mixing with geometrically decaying coefficients.
\end{proposition}

\subsection{A new thinning-based display of CP-INGARCH models}
\label{sec:alternative_formulation}

\begin{figure}
\center
(a) GINARMA(1, 1) model without compounding step:

\medskip

\includegraphics[scale=0.85]{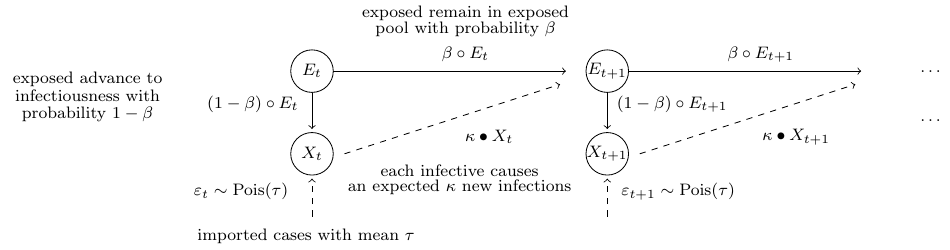}

\medskip

(b) GINARMA($2$, $2$) model without compounding step:

\medskip

\includegraphics[scale = 0.85]{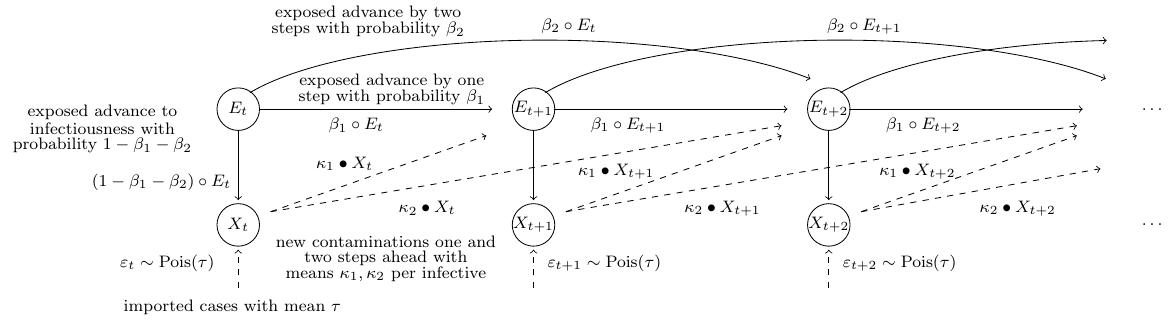}

\smallskip

(c) GINARMA(1, 1) model with compounding step:

\smallskip

\includegraphics[scale = 0.85]{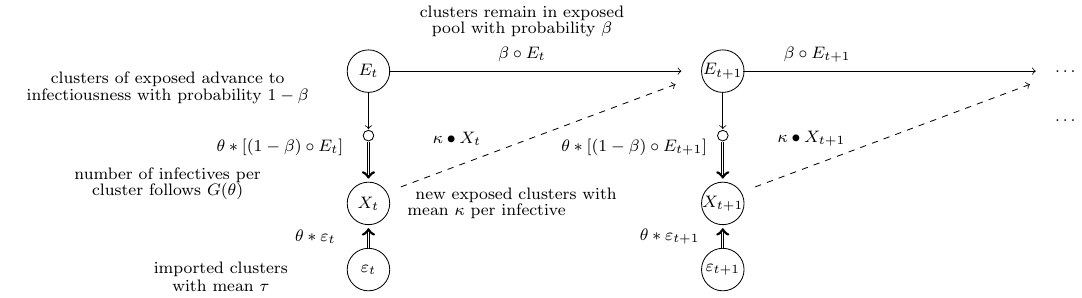}

\smallskip

\caption{Interpretation of GINARMA models as stochastic epidemic processes. (a) GINARMA(1, 1) without a compounding step, see equation \eqref{eq:compounding_omitted}. (b) GINARMA(2, 2) without a compounding step. (c) Full GINARMA(1, 1) model including compounding step. Solid lines represent multinomial thinning, dashed lines generalized thinning or immigration, and double lines thinning with a cluster distribution $G(\theta)$. In the bottom panel, the small circles represent the intermediate step $(1 - \beta) \circ E_t$, which is then subject to compounding.}
\label{fig:ingarch_flowchart_poisson}
\end{figure}

If $\sum_{j = 1}^q \beta_j < 1$, various CP-INGARCH models can be obtained by using Poisson offspring in our GINARMA class; see Supplement \ref{appendix:proofs_equivalence} for the derivations. With $\tau = \nu/(1 - \beta)$ and $\kappa = \alpha/(1 - \beta)$, the Poisson INGARCH(1, 1) model \eqref{eq:ingarch11} can be represented as
\begin{align}
X_t & = (1 - \beta) \circ E_t + \import_t, \ \ \ \  \label{eq:X_t_thinning_Poisson}
E_t = \beta \circ E_{t - 1} \ + \ \kappa \star X_{t - 1}, \ \ \ \ \import_t \stackrel{\text{i.i.d.}}{\sim} \text{Pois}(\tau). 
\end{align}
The two thinnings of $E_t$ are coupled as in \eqref{eq:mult_general} while $E_0 \sim \text{Pois}(\eta), \eta = (\lambda_0 - \tau)/(1 - \beta)$. We do not assume $\alpha + \beta < 1$, i.e., do not require stationarity.  A technical condition (that also applies in the following) is  $\lambda_{0} \geq \tau$ so that $\eta \geq 0$. This, however, is natural as $\lambda_t \geq \tau$ holds for all $t \geq 1$ if $\lambda_{0} \geq \tau$. The model structure corresponds to the top panel in Figure \ref{fig:ingarch_flowchart_poisson}.

\begin{remark}
    In epidemic modelling, Poisson offspring are widely used \cite{Farrington1999, Kucharski2014}. They arise e.g., when approximating a Reed–Frost chain binomial model for a large population \cite{Bauer2018}.
\end{remark}
 
\noindent Formulation \eqref{eq:X_t_thinning_Poisson} can be extended to the Poisson INGARCH($p$, $q$) case by setting
\begin{align}
X_t & = \left(1 - \sum_{j = 1}^q \beta_j \right) \circ E_t + \import_t,  \label{eq:Xt_thinning_pq}\ \ \ \ 
E_t = \sum_{j = 1}^q \beta_j \circ E_{t - j} \ \ + \ \ \sum_{i = 1}^p \kappa_i \star X_{t - i}
\end{align}
with $\tau = \nu / (1 - \sum_{j = 1}^q\beta_j), \kappa_i = \alpha_i / (1 - \sum_{j = 1}^q\beta_j)$. For initialization we need to set ${E_{m} \stackrel{\text{ind.}}{\sim} \text{Pois}(\eta_m)}$ with $\eta_m = (\lambda_m - \tau )/(1 - \sum_{j = 1}^q\beta_j), m = 1 - q, \dots, 0$.
This corresponds to the structure displayed in the middle panel of Figure \ref{fig:ingarch_flowchart_poisson}.

A CP-INGARCH(1, 1) model as in \eqref{eq:N_CP_original}--\eqref{eq:lambda_CP_original} is obtained by extending equation \eqref{eq:X_t_thinning_Poisson} to
\begin{align}
X_t & = \clustersize * \left[(1 - \beta) \circ E_t + \import_t\right],\label{eq:X_t_CP}
\end{align}
where $*$ denotes thinning with the clustering distribution $G$ (see equation \eqref{eq:def_compund_thinning}). Here, we need to set $\tau = (\nu/\clustersize)/(1 - \beta), \kappa = (\alpha/\clustersize)/(1 - \beta)$ and $\eta = (\lambda_0/\theta - \tau)/(1 - \beta)$.
This extension corresponds to the bottom panel of Figure \ref{fig:ingarch_flowchart_poisson}.

Proposition \eqref{proposition:geometric_ergodicity} thus implies geometric ergodicity of CP-INGARCH(1, 1) models, a topic that has received much attention (e.g., \cite{Davis2021, Fokianos2009, Goncalves2015} and references therein). Typically, the employed arguments are more sophisticated than what we use, the difficulty being that the state space of $\lambda_{t + 1} \mid \lambda_{t}$ depends on $\lambda_t$. We circumvent this via a fully discrete display. We note that a construction of the Poisson INGARCH(1, 1) using a ``cascade of thinning operations'' was already introduced in \cite{Ferland2006}, but is more complex than our representation.

\section{Extending the INAR class}
\label{sec:extension_inar}

\subsection{Defining a new INARMA($p$, $q$) model} \label{subsec:inarma}

Despite the parallels between the Poisson INAR(1) and INARCH(1) models seen in Section \ref{subsec:poisson_ginar_inarch}, it is not obvious how the INGARCH($p, q$) recursion \eqref{eq:lambda_CP_original} could be transposed to the INAR case. In the thinning-based representation \eqref{eq:Xt_thinning_pq}, however, we can simply swap all Poisson thinnings for multinomial thinnings. This leads to a new extension of the INAR($p$) model \eqref{eq:inar_p}--\eqref{eq:multinomial_thinning}, which as we shall see has attractive stochastic properties. Omitting the compounding step $\theta \: *$ from Definition \eqref{def:ingarma}, we define our INARMA($p, q$) process as
\begin{align}
X_t & = \Bigl(1 - \sum_{j = 1}^q \beta_j \Bigr) \circ E_t + \import_t\label{eq:X_inarma}\\
E_t & = \sum_{j = 1}^q \beta_j \circ E_{t - j} \ \ + \ \ \sum_{i = 1}^p \kappa_i \circ X_{t - i}\label{eq:E_inarma}.
\end{align}
In addition to the constraints from Definition \ref{def:ingarma}, we assume that $\sum_{i = 1}^p \kappa_i < 1$ and set
\begin{align}
(\kappa_1 \circ X_t, \dots, \kappa_p \circ X_t) \ \mid \ X_t & \sim \text{Mult}\left(X_t; \kappa_1, \dots, \kappa_p\right),\label{eq:C_t_inarma}
\end{align}
thus paralleling equation \eqref{eq:multinomial_thinning}. The imports $\{\import_t\}$, thinnings of $E_t$ and the initialization with $X_{1 - p}, \dots, X_0, E_{1 - q}, \dots, E_0$ are handled as in Definition \ref{def:ingarma}.

\begin{remark}
In terms of the interpretation from Section \ref{subsec:interpretation_epidemic}, an infected can cause at most one new infection in the INAR/INARMA model. This case is occasionally studied in theory, but corresponds to an unusual practical setting. As argued by \cite{Farrington1999}, the disease would need to be of very modest infectivity (low $\kappa_i$), or infectives would need to be isolated systematically after a first event of onward transmission.
\end{remark}

\subsection{Properties of the INARMA(1, 1) model}

\subsubsection{General import distributions}
\label{subsubsec:general_imports}

Many statements from Section \ref{eq:stochastic_properties_general} simplify considerably for the INARMA(1, 1) model.

\begin{lemma}
\label{lemma:moments_inarma_11}
As we assumed $\kappa < 1$, the limiting stationary mean, variance and autocorrelation function of an INARMA(1, 1) process $\{X_t\}$ are given by
\begin{align}
\mu_X & = \frac{\tau}{1 - \kappa}\label{eq:mu_X}\\
\sigma^2_X & = \frac{\kappa(1 + \beta)}{1 + \xi} \times \frac{\tau}{1 - \kappa}\label{eq:sigma2_X}  \ \ \ + \ \ \  \left(1 - \frac{\kappa(1 + \beta)}{1 + \xi}\right) \times \frac{\sigma^2_\tau}{1 - \kappa} \\ 
\gamma_X(d) & = (1 - \beta)\kappa\xi^{d - 1} \times\left(1 + \frac{\kappa\beta(\sigma^2_\tau - \tau)}{(1 + \beta) \{(1 - \kappa)\sigma^2_\tau + \kappa\tau\} + (1 - \beta)\kappa\sigma^2_\tau}\right) \times \sigma^2_X.\label{eq:rho_X}
\end{align}
Here, we use the shorthand
\begin{equation}
\xi = \gamma_X(2)/\gamma_X(1) = \beta + (1 - \beta)\kappa.
\label{eq:define_xi}
\end{equation}
\end{lemma}

\noindent If relative to a Poisson, the import distribution is overdispersed ($\sigma^2_\tau > \tau$) or underdispersed ($\sigma^2_\tau < \tau$), respectively, the same thus holds for the marginal distribution of $\{X_t\}$. Moreover, for overdispersed (underdispersed) imports, the autocorrelations will be stronger (weaker) than for identical $\tau, \beta, \kappa$ and equidispersed imports. It is easily shown that $\gamma_X(2) \geq \gamma_X(1)^2$ always holds in \eqref{eq:rho_X}, with equality for $\beta = 0$. Like the INGARCH(1, 1), our INARMA(1, 1) thus has a ``longer memory'' than an INAR(1), and there is no instance with an MA(1) structure. This differs from the $\text{INARMA}_\text{DGL}(1, 1)$, i.e., model \eqref{eq:dion_inarma} with binomial thinning, which implies $\gamma_X(2) \geq \gamma_X(1)^2$ and contains the INMA(1) as a special case.

\begin{remark}
    In the INARMA(1, 1), $\{E_t\}$ is an INAR(1) process, while $\{X_t\}$ is an INAR($\infty$) with geometrically decaying autoregressive parameters; see Supplementary Remark \ref{remark:inarma11}.
    \label{remark:embedded_inar}
\end{remark}

\subsubsection{Poisson imports}

The Poisson INARMA(1, 1) model has been discussed in \cite{Bracher2019}, but for completeness some results are repeated and extended.
\begin{lemma}
\label{proposition:bivariate_poisson}
If $\{X_t\}$ is a Poisson INARMA(1, 1) process with $E_0 \sim \text{Pois}[\kappa\tau/(1 - \xi)]$, then the process is strictly stationary with Poisson marginals.  Expressions \eqref{eq:mu_X}--\eqref{eq:rho_X} simplify to 
\begin{equation}
\mu_X = \sigma^2_X = \frac{\tau}{1 - \kappa}; \ \ \rho_X(d) = (1 - \beta)\kappa\xi^{d - 1}.\label{eq:moments_poisson_11}
\end{equation}
Moreover, for $t \in \mathbb{N}$ and $d > 0$ it holds that
\begin{equation}
(X_t, X_{t + d}) \sim \textnormal{BPois}\{\rho_X(d)\mu_X, [1 - \rho_X(d)]\mu_X, [1 - \rho_X(d)]\mu_X\},\label{eq:bivariate_poisson}
\end{equation}
with BPois the bivariate Poisson distribution as defined in \cite{Johnson1997}.
\end{lemma}
Further particularities of the Poisson INARMA(1, 1) process include that it is time-reversible and closed to binomial thinning. The latter means that if $\{X_t\}$ is a Poisson INARMA(1, 1), then so is $\{\tilde{X}_t\}$ with
$\tilde{X}_t = \pi \stackrel{\textnormal{ind}}{\circ} X_t$; see \cite{Bracher2019} and Supplementary Remark \ref{remark:underreporting}.

\begin{remark}
    In terms of the thinning-based formulation \eqref{eq:X_t_thinning_Poisson}, the limiting-stationary second-order properties of the Poisson INGARCH(1, 1) model are $\mu_X = \tau/(1 - \kappa)$ and
    \begin{align*}
    \sigma^2_X = \left(1 + \frac{(1 - \beta)^2 \kappa^2}{1 - \xi^2}\right) \times \mu_X,\ 
    \rho_X(d) = \left(1 + \frac{\kappa\beta(1 - \beta)}{1 - \xi^2 + \kappa^2(1 - \beta)^2}\right) \times (1 - \beta)\kappa\xi^{d - 1}.
    \end{align*}
    Poisson rather than binomial offspring thus lead to higher dispersion and stronger autocorrelations than in the Poisson INARMA(1, 1), but the ACFs are proportional.
\end{remark}

\subsubsection{Compound Poisson imports}

As the class of CP distributions is closed to binomial thinning and summation, the marginal distributions of INAR(1) models with CP imports are CP \cite{Schweer2014}. Remark \ref{remark:embedded_inar} thus implies that in the CP-INARMA(1, 1) process, $\{E_t\}$ has CP marginals, and it is easy to show that $X_t$ inherits this property. The order of the CP distribution, too, will be inherited \cite{Weiss2015b}, meaning that Hermite imports lead to Hermite marginals; see Supplementary Remark \ref{remark:hermite}.

A structural difference between CP-INARMA and CP-INGARCH models is that the latter feature a compounding step acting on the Poisson imports and offspring, see equation \eqref{eq:X_t_CP}. In CP-INARMA models, on the other hand, no compounding step is used, but the import distribution becomes a CP. Overdispersion thus enters purely via the imports.

\subsection{Properties of the Poisson INARMA($p$, $q$) model}

For $p, q > 1$, the model becomes considerably more complex and relevant properties can only be established for the Poisson case.

\begin{corollary}
\label{corollary:moments_inarma_pq}
The limiting-stationary mean of a Poisson INARMA($p, q$) process $\{X_t\}$ is
$$
\mu_X = \sigma^2_X = \frac{\tau}{1 - \sum_{i = 1}^p \kappa_i}
$$
while the autocorrelation function can be computed recursively via
\begin{align*}
\rho(0) = 1 \ \ \text{and} \ \ \rho_X(d) & = \left(1 - \sum_{j = 1}^q \beta_j\right) \times \left(\sum_{i = 1}^d \rho_X(d - i) \times s_i\right),\\
\text{where} \ \ s_i & = \sum_{k = 1}^{\min\{i, p\}} \kappa_k \pi_{i - k} \ \ \textnormal{ for } i = 1, 2, \dots \\
\pi_0 = 1, \ \ \pi_k & = \sum_{l = 1}^{\min\{k, q\}} \beta_l \pi_{k - l} \ \ \textnormal{ for } k = 1, 2, \dots
\end{align*}
Moreover, the bivariate Poisson property \eqref{eq:bivariate_poisson} from Lemma \ref{proposition:bivariate_poisson} still holds.
\end{corollary}

\begin{corollary}
\label{corollary:embedded_inar_p_main}
In the Poisson INARMA($p, q$) model, the process $\{E_t\}$ is a Poisson INAR$(\max\{p, q\})$ process as defined in equations \eqref{eq:inar_p}--\eqref{eq:multinomial_thinning}; see Lemma \ref{lemma:embedded_inar} in the Supplement for details.
\end{corollary}

\begin{remark}
Following arguments from \cite{Alzaid1990}, it can be shown that the ACF of the INARMA($p, q$) model with general import distribution bears resemblance with, but is not identical to that of a Gaussian ARMA($\max\{p, q\}, \max\{p, q\}$) model; see Supplementary Remark \ref{remark:details_acf_arma}.
\label{remark:arma}
\end{remark}

\subsection{Inference in the case $p = q = 1$}

\subsubsection{An algorithm for likelihood evaluation}

Evaluating the likelihood function for an INARMA(1, 1) model is considerably more difficult than for the INGARCH(1, 1) model, where equation \eqref{eq:inarch1} facilitates computations.  We suggest an adaptation of the forward algorithm \cite{Zucchini2009} to this end, resembling an existing procedure \cite{Weiss2019} for the $\text{INARMA}_{\text{DGL}}(1, 1)$ model. To facilitate notation in the following we introduce the shorthand $A_t = (1 - \beta) \circ E_t$ which enables us to write
\begin{align*}
X_t & = \underbrace{A_t}_{(1 - \beta) \circ E_t} \ + \ \import_t, \ \ \ \ \ 
E_t = \underbrace{E_{t - 1} - A_{t - 1}}_{\beta \circ E_{t - 1}} \ +\ \kappa \circ X_{t - 1}.
\end{align*}
Now denote the sequence of observed values by $\{x_1, \dots, x_T\}$, with $T$ the length of the time series. As a first step, a sufficiently large support $\mathcal{E} = \{0, 1, \dots, M\}$ needs to be chosen for $\{X_t\}$ and $\{E_t\}$. In practice we set $M$ to the maximum of $1.2 \times \max(x_1, \dots, x_t)$ and the 0.999 quantiles of the stationary distributions of $E_t$ and $X_t$ under the respective parameters. As $E_t \geq A_t$, $\mathcal{E}$ implies a support $\mathcal{E}^* = \{(0, 0), (1, 0), (1, 1), (2, 0), \dots, (M, M - 1), (M, M)\}$ for the tuple $(E_t, A_t)$. Moreover we introduce the following shorthands: $\Pr(Y = y \ \mid \ X_{ < t})$ is the probability that $Y = y$ provided that $X_{t - 1} = x_{t - 1}, X_{t - 2} = x_{t - 2}, \dots, X_1 = x_1$; for $t= 1$ this corresponds to the marginal distribution of $Y$ under some suitable initialization. $\Pr(Y = y \ \mid \ X_{ \leq t})$ is defined analogously, but $X_t = x_t$ is also included in the condition.

\begin{algorithm}\label{algorithm}
The algorithm is initialized by setting $E_1 \sim \text{Pois}(\eta)$, with $\eta$ treated like an extra parameter. Then the following steps are iterated for $t = 1, \dots, T$.
\begin{enumerate}
\item For each tuple $(e_t, a_t) \in \mathcal{E}^*$ compute
$$
\Pr(E_t = e_t, A_t = a_t \ \mid \ X_{< t}) = \Pr(A_t = a_t \ \mid \ E_t = e_t)\times \Pr(E_t = e_t \ \mid \ X_{< t}).
$$
\item Compute and store
\begin{align*}
\Pr(X_t = x_t \ \mid \ X_{< t}) & = \sum_{(e_t, a_t) \in \mathcal{E}^*} \Pr(X_t = x_t \ \mid \ E_t = e_t, A_t = a_t) \times \Pr(E_t = e_t, A_t = a_t \ \mid \ X_{< t})\\
& = \sum_{(e_t, a_t) \in \mathcal{E}^*} \Pr(\import_t = x_t - a_t) \times \Pr(E_t = e_t, A_t = a_t \ \mid \ X_{< t}).
\end{align*}
\item For $(e_t, a_t) \in \mathcal{E}^*$ compute
\begin{align*}
\Pr(E_t = e_t, A_t = a_t \ \mid \ X_{\leq t}) & = \frac{\Pr(E_t = e_t, A_t = a_t, X_t = x_t \ \mid \ X_{< t})}{\Pr(X_t = x_t \ \mid \ X_{< t})} \\
& = \frac{\Pr(E_t = e_t, A_t = a_t \ \mid \ X_{< t}) \times \Pr(\import_t = x_t - a_t)}{\Pr(X_t = x_t \ \mid \ X_{< t})}.
\end{align*}
\item For $l_t \in \mathcal{E}$ compute
$$
\Pr(E_t - A_t = l_t \ \mid \ X_{\leq t}) = \sum_{(e_t, a_t) \in \mathcal{E}^*: e_t - a_t = l_t} \Pr(E_t = e_t, A_t = a_t \ \mid \ X_{\leq t}).
$$
\item For $e_{t + 1} \in \mathcal{E}$ compute
\begin{align*}
\Pr(E_{t + 1} = e_{t + 1} \ \mid \ X_{\leq t}) & = \sum_{l_t \in \mathcal{E}} \Pr(E_{t + 1} = e_{t + 1} \ \mid \ E_t - A_t = l_t, X_{\leq t}) \times \Pr(E_t - A_t = l_t \ \mid \ X_{\leq t})\\
& = \sum_{l_t \in \mathcal{E}} \Pr(\kappa \circ x_t = e_{t + 1} - l_t) \times \Pr(E_t - A_t = l_t \ \mid \ X_{\leq t}).
\end{align*}
\end{enumerate}
\end{algorithm}
The values stored in Step 2 of each iteration serve to evaluate the (conditional) likelihood of the observed time series as
$$
\Pr(X_1 = x_1, \dots X_t = x_T) = \Pr(X_1 = x_1) \times \prod_{t = 2}^T \Pr(X_t = x_t \ \mid \ X_{< t}).
$$
Maximization of the log-likelihood is done using the Nelder-Mead method as implemented in the \texttt{R} function \texttt{optim}. We use moment estimators (see next section) to initialize the optimization. All parameters are handled on suitable transformed scales allowing for unconstrained optimization. For parameters constrained to the unit interval we use logit transformations, for parameters which can take any positive value we use the natural logarithm. Standard errors are estimated via the inverse observed Fisher information (obtained by numerical differentiation) with subsequent application of the delta method. Fitted values and Pearson residuals can be obtained using the probabilities $\Pr(X_t = x_t \ \mid \ X_{< t}), x_t \in \mathcal{E}$.

As the likelihood function is not available in closed form, establishing consistency or asymptotic normality of the estimators is not straightforward (a typical proof strategy relying on threefold continuous differentiability of the log-likelihood function \cite{Fokianos2009}). Indeed, results on the asymptotics of maximum likelihood estimators seem to be lacking even for INAR models with general import distributions. Saddlepoint approximations as sugested for INAR($p$) models \cite{Pedeli2015} may be a useful alternative for fast and principled inference. 

\subsubsection{Moment-based estimation}

The suggested likelihood evaluation method can get slow even for moderately high count values. As a computationally fast alternative we consider moment-based estimators, see Supplement \ref{suppl_moment_based_estimation}. In the case of Poisson innovations, consistency and asymptotic normality can be established. For general innovation distributions, the estimators do not have a closed form, but can be evaluated by solving a cubic equation numerically. We note that least squares estimation as often employed for INAR(1) models is hampered by the difficulty of computing $\mathbb{E}(X_t \mid X_{t - 1},  \dots, X_1)$, which essentially requires application of Algorithm \ref{algorithm}.

\subsubsection{Simulation study}

To assess the behaviour of our estimators we specify three simulation scenarios:

\begin{enumerate}
\item Scenario 1: $\tau = 1, \beta = 0.5, \kappa = 0.5$. In the Poisson case this implies $\mu_X = \sigma^2_X = 2, \rho_X(d) = 0.25 \times 0.75^{d - 1}$. In the negative binomial / Hermite cases we set $\disp = 0.5$, resulting in $\mu_X = 2; \sigma^2_X = 2.57; \rho_X(d) = 0.26\times 0.75^{d - 1}$.
\item Scenario 2: $\tau = 1, \beta = 0.2, \kappa = 0.6$. In the Poisson case this implies $\mu_X = \sigma^2_X = 2.5; \rho_X(d) = 0.48 \times 0.68^{d - 1}$. In the negative binomial / Hermite cases we set $\disp = 0.7$, resulting in $\mu_X = 2.5; \sigma^2_X = 3.5; \rho_X(d) = 0.5\times 0.68^{d - 1}$.
\item Scenario 3: $\tau = 1, \beta = 0.1, \kappa = 0.8$. In the Poisson case this implies $\mu_X = \sigma^2_X = 5; \rho_X(d) = 0.72 \times 0.82^{d - 1}$. In the negative binomial / Hermite cases we set $\disp = 0.9$, resulting in $\mu_X = 5; \sigma^2_X = 7.32; \rho_X(d) = 0.74\times 0.82^{d - 1}$.\end{enumerate}

\noindent As we chose $\tau = 1$ in all settings we get the same second-order properties for the same values of $\disp$ in the Hermite and the negative binomial versions. Note however that this is not generally the case. We simulated 1000 time series for each scenario and different lengths of time series $T \in \{250, 500, 1000\}$. The results for maximum-likelihood and moment-based estimation can be found in Table \ref{tab:sim_ml} and Supplementary Table \ref{tab:sim_moments}, respectively.

\begin{table}[h!]
\spacingset{1.3} 
\scriptsize

\caption{\footnotesize Simulation results for maximum likelihood estimators in the Poisson, Hermite and negative binomial settings, scenarios 1--3, $T \in \{250, 500, 1000\}$ and 1000 runs. In 0.3\% of runs estimates and/or standard errors $\widehat{\text{se}}$ could not be evaluated due to numerical problems.}
\label{tab:sim_ml}
\center
\begin{tabular}{p{0.3cm} @{\hskip 0.7cm} p{0.45cm} p{0.45cm} p{0.45cm} p{0.6cm} @{\hskip 0.6cm} p{0.45cm} p{0.45cm} p{0.45cm} p{0.6cm} @{\hskip 0.6cm} p{0.45cm} p{0.45cm} p{0.45cm} p{0.6cm} @{\hskip 0.6cm} p{0.45cm} @{\hskip 0.45cm} p{0.45cm} p{0.6cm} p{0.65cm}}
\hline 
\multicolumn{17}{c}{Poisson}\\
\hline 
$T$ & \multicolumn{4}{c}{$\tau$} & \multicolumn{4}{c}{$\disp$} & \multicolumn{4}{c}{$\beta$} & \multicolumn{4}{c}{$\kappa$}\\
\hline
& true & mean & se & mean & true & mean & se & mean & true & mean & se & mean & true & mean & se & mean \\
& & & & of $\widehat{\text{se}}$ & & & & of $\widehat{\text{se}}$ & & & & of $\widehat{\text{se}}$ & & & & of $\widehat{\text{se}}$\\
\hline
\input{table/sim_pois_sc1.tex}

\input{table/sim_pois_sc2.tex}

\input{table/sim_pois_sc3.tex}\\
\hline
\multicolumn{17}{c}{Hermite}\\
\hline 
$T$ & \multicolumn{4}{c}{$\tau$} & \multicolumn{4}{c}{$\disp$} & \multicolumn{4}{c}{$\beta$} & \multicolumn{4}{c}{$\kappa$}\\
\hline
& true & mean & se & mean & true & mean & se & mean & true & mean & se & mean & true & mean & se & mean \\
& & & & of $\widehat{\text{se}}$ & & & & of $\widehat{\text{se}}$ & & & & of $\widehat{\text{se}}$ & & & & of $\widehat{\text{se}}$\\
\hline
\input{table/sim_herm_sc1.tex}

\input{table/sim_herm_sc2.tex}

\input{table/sim_herm_sc3.tex}\\
\hline 
\multicolumn{17}{c}{Negative binomial}\\
\hline 
$T$ & \multicolumn{4}{c}{$\tau$} & \multicolumn{4}{c}{$\disp$} & \multicolumn{4}{c}{$\beta$} & \multicolumn{4}{c}{$\kappa$}\\
\hline
& true & mean & se & mean & true & mean & se & mean & true & mean & se & mean & true & mean & se & mean \\
& & & & of $\widehat{\text{se}}$ & & & & of $\widehat{\text{se}}$ & & & & of $\widehat{\text{se}}$ & & & & of $\widehat{\text{se}}$\\
\hline
\input{table/sim_negbin_sc1.tex}

\input{table/sim_negbin_sc2.tex}

\input{table/sim_negbin_sc3.tex}

\end{tabular}
\end{table}

Overall, both fitting procedures yield approximately unbiased estimates for $\tau, \beta$ and $\kappa$, with some small-sample biases. The dispersion parameters $\disp$ are subject to some biases and their estimation can become instable if $\hat{\tau}$ is small. The maximum-likelihood estimators have considerably smaller standard errors than their moment-based counterparts, reflecting the well-known inefficiency of moment estimators in models with MA components. In the maximum likelihood scheme, the estimated standard errors are mostly in good agreement with the observed standard errors, but in some instances underestimate the true variability.

\section{Application: childhood diseases in Bavaria}
\label{sec:real_data}

\begin{figure}[h!]
\center
\includegraphics[scale=0.71]{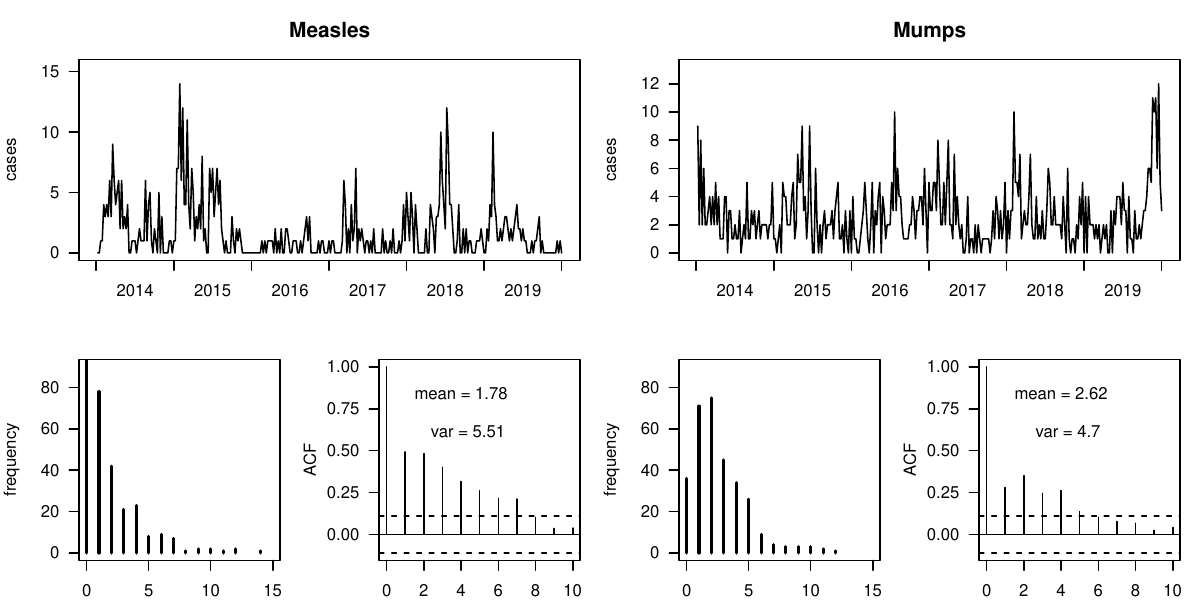}
\caption{Top: weekly counts of reported measles and mumps cases in Bavaria, 2014--2019. Bottom: marginal distributions and autocorrelation functions.}
\label{fig:data}
\end{figure}

We now apply various instances of the introduced model class to two time series of infectious disease counts. We consider weekly numbers of reported measles and mumps cases in the German state of Bavaria, 2014--2019. Measles and mumps are vaccine-preventable childhood diseases and have become rare in Western Europe. While both diseases exhibit seasonal patterns in the absence of vaccination, during the considered period they only occurred sporadically. They thus match the setting described in Section \ref{subsec:interpretation_epidemic} well, and are indeed commonly modelled using subcritical branching processes \cite{Chiew2014, DeSerres2000}. The data, available from Robert Koch Institute (\texttt{\url{https://survstat.rki.de}}), are displayed in Figure \ref{fig:data}. Both series exhibit slowly decaying autocorrelation functions and some degree of overdispersion.

\begin{figure}[h!]
\center
\includegraphics[scale=0.76]{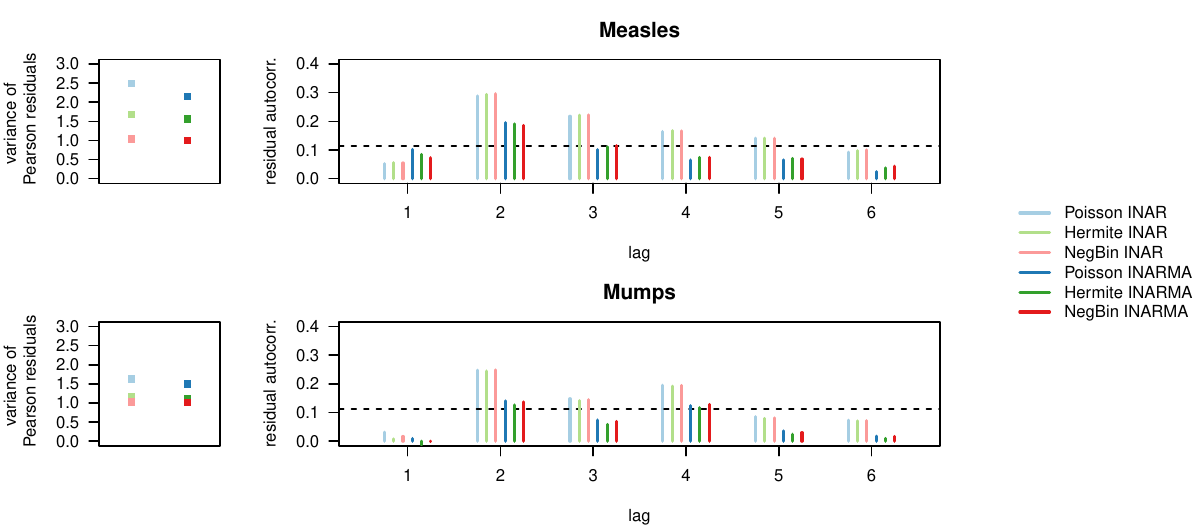}
\vspace{-5mm}

\caption{Analysis of Pearson residuals of INAR and INARMA models. Left: Variances, which should be close to 1. Right: Autocorrelations, which should be small. The dashed line shows $2/\sqrt{T}$, i.e., the 97.5\% quantile for the empirical ACF of white noise.
}
\label{fig:residuals}
\end{figure}

Table \ref{tab:fits_real_data_epi} summarizes the fits of the INARCH(1), INGARCH(1, 1), INAR(1) and INARMA(1, 1) models. Each of them was applied in the Poisson, Hermite and negative binomial version. To make the results more easily comparable across the different models, we present them in terms of the epidemiological interpretation from Section \ref{subsec:interpretation_epidemic}; results for the original parameterizations are shown in Supplementary Table \ref{tab:fits_real_data}. The mean generation times obtained from the INGARCH and INARMA models are in good agreement with commonly used estimates from the literature (slightly below 2 weeks for measles; 18 days for mumps, \cite{Bjoernstad2002, Vink2014}). The estimated reproductive numbers are highest for the Hermite and negative binomial INGARCH models, which feature the most flexible offspring mechanisms. Here, they are around 0.7 and 0.6 for measles and mumps, respectively. While no comparable estimates for Germany exist, these values seem plausible in light of estimates from Australia, a country with somewhat higher vaccination coverage ($\Reff = 0.47$ to 0.65 for measles depending on the exact method for 2009--2011, \cite{Chiew2014}). 

In terms of the Akaike information criterion (AIC), both a more flexible autocorrelation structure (i.e., INGARCH or INARMA) and accounting for overdispersion considerably improves model fits. The INGARCH approach, where overdispersion enters both via the import and offspring distributions, leads to better results than the INARMA, where the offspring distribution is always Bernoulli. Figure \ref{fig:residuals} shows an analysis of the Pearson residuals of the INAR and INARMA models. The Pearson residuals are too dispersed for the Poisson version (variance exceeding 1); the negative binomial version can remedy this, while the Hermite model only partly does so. The INAR models show pronounced residual autocorrelation at lags 2 and 3, which is largely remedied by the INARMA versions. For the INARCH and INGARCH models, the picture is similar, see Supplementary Figure \ref{fig:residuals_ingarch}. Graphical representations of the fits are provided in Supplementary Figures \ref{fig:fit_measles} and \ref{fig:fit_mumps}.

Despite the good agreement with literature estimates, we emphasize that aggregate-level analyses like the above should not be overinterpreted; see discussion in the next Section.


\begin{table}[h!]
\footnotesize
\caption{\footnotesize Model fits for measles and mumps. For better comparability, estimates are presented in terms of the epidemiological interpretation of the models. IC $ = \tau\clustersize$ is the mean weekly number of imported cases; $\Reff = \kappa\clustersize$ is the effective reproductive number; GT $ = 1/(1 - \beta)$ is the mean generation time in weeks; CS $ = \clustersize$ is the mean cluster size. The definitions of IC, $\Reff$, GT, CS refer to the general formulation \eqref{eq:X_general}--\eqref{eq:E_general} with $p = q = 1$. Entries with asterisks$^*$ are not estimated but implied by the model definition.}
\label{tab:fits_real_data_epi}
\center
\begin{tabular}{p{3.5cm} @{\hskip 0.8cm} p{0.45cm} p{0.45cm} p{0.45cm} p{0.45cm} p{0.45cm} @{\hskip 2cm} p{0.45cm} p{0.45cm} p{0.45cm} p{0.45cm} p{1cm}}
\hline
 & & \multicolumn{3}{c}{Measles} & & \multicolumn{5}{c}{Mumps}\\
\hline 
Model & IC & $\Reff$ & GT & CS & AIC & IC & $\Reff$ & GT & CS & AIC\\
\hline
\input{table/tab_inarch_epi.tex} \\
\input{table/tab_ingarch_epi.tex} \\
\input{table/tab_inar_epi.tex} \\
\input{table/tab_inarma_epi.tex}

\end{tabular}
\end{table}

\section{Discussion}
\label{sec:discussion}

In this paper we introduced an overarching class of count time series models, which includes many popular CP-INAR and INGARCH models. Each of them is characterized by an immigration, an offspring and a clustering distribution. We gave particular attention to a new INARMA($p, q$) model which mirrors the INGARCH($p, q$) formulation. Numerous other instances could be examined, for example models based on other thinning operators \cite{Joe1996, Scotto2015}. Other potential avenues are the inclusion of covariates and multivariate extensions.

We note that our class only comprises linear CP-INGARCH models with a time-constant clustering distribution. For instance, the negative binomial INGARCH model by \cite{Zhu2011}, which features a parameter $\theta_t$ that depends on $\lambda_t$, is not contained; nor are log-linear models \cite{Fokianos2009} or other variations where the linearity assumption in \eqref{eq:ingarch11} is relaxed.

As noted before, the introduced INARMA(1, 1) class only allows for ACFs which from lag 2 onwards decay \textit{more slowly} than in the corresponding INAR(1) model. For other INARMA models suggested in the literature \cite{Dion1995}, the converse is true. It would be desirable to construct a model able to accommodate both patterns.

Concerning the real-data application, several caveats are needed. Firstly, our aggregate analysis glosses over population heterogeneities, ignoring e.g., that non-vaccination may be clustered in certain groups. Given the sparse data, we pragmatically assumed constant $R_e$ within and across seasons. We consider this acceptable for the pre-COVID-19 period, but it would certainly not hold for the years since. Routine surveillance data are moreover subject to many biases, including reporting delays and underreporting, which can distort estimates of $R_e$ \cite{Bracher2021}. These aspects can moreover vary over time e.g., due to changes in healthcare seeking or testing practices. In practice, branching process models are usually not fitted to surveillance counts alone, but also data e.g., on the type of infection (imported/domestic), and estimates based on different data types are compared to assess robustness \cite{Chiew2014}. This will yield more reliable estimates than we provide in our illustrative example.

\textbf{Reproducibility:}
An R package implementing the presented estimation method is available at \url{https://github.com/jbracher/rinarma}. Data and code to reproduce all results are available at \url{https://github.com/jbracher/ginarma}.

\textbf{Acknowledgements:}
We would like to thank Mirko Armillotta, Konstantinos Fokianos, Melanie Schienle and Christian Wei{\ss} for discussions on earlier versions of the paper. Both authors were supported by the German Research Foundation (DFG), project 512483310.

{
\footnotesize
\bibliographystyle{plain}

}

\newpage

\begin{center}
{\LARGE Supplementary material for \textit{A class of count time series models uniting compound Poisson INAR and INGARCH models}}

\bigskip

\end{center}

\appendix

\setcounter{page}{1}
\renewcommand{\thepage}{S\arabic{page}}
\renewcommand{\thesection}{\Alph{section}}
\renewcommand{\thetable}{S\arabic{table}}
\renewcommand{\theequation}{S\arabic{equation}}
\renewcommand{\thelemma}{S\arabic{lemma}}
\renewcommand{\theremark}{S\arabic{remark}}
\renewcommand{\theproposition}{S\arabic{proposition}}
\renewcommand{\thecorollary}{S\arabic{corollary}}
\renewcommand\thefigure{S\arabic{figure}}

\setcounter{figure}{0}
\setcounter{table}{0}
\setcounter{figure}{0}
\setcounter{equation}{0}
\setcounter{remark}{0}
\setcounter{lemma}{0}
\setcounter{proposition}{0}
\setcounter{corollary}{0}

\section{Details on Section \ref{sec:preliminaries}: Compound Poisson INAR and INGARCH models}

In the following we present some additional description of the two considered compound Poisson distributions and the associated INGARCH models. As going back between classical, regression-type and compound-Poisson parameterizations can be tedious we hope that this more detailed account will be a useful addition.

\subsection{The Hermite distribution and associated models}
\label{suppl:hermite}

\subsubsection{Parameterization and representation as a compound Poisson}
\label{suppl:param_hermite}

As mentioned in the main text, a random variable $Y$ is Hermite-distributed if it can be written as
$$
Y = A_1 + 2A_2,
$$
where independently $A_1 \sim \text{Pois}(\lambda_1), A_2 \sim \text{Pois}(\lambda_2).$ Expressed through $\lambda_1, \lambda_2$, the probability mass function is \cite{Kemp1965}
\begin{equation}
\text{Pr}(Y = y) = \exp(-\lambda_1 - \lambda_2) \times \sum_{j = 0}^{[y/2]} \frac{\lambda_1^{y - 2j}\lambda_2^{2j}}{(y - 2j)!j!}, y = 0, 1, 2, \dots \label{eq:pmf_hermite1}    
\end{equation}
where $[y/2]$ is the integer part of $y/2$.

In the main text, we use an alternative parameterization where we characterize the distribution via its mean $\mu = \lambda_1 + 2\lambda_2$ and a dispersion parameter $\disp = 2\lambda_2/(\lambda_1 + 2\lambda_2)$. This implies the probability mass function
\begin{equation}
\text{Pr}(Y = y) = \exp\left[\mu\left(-1 + \frac{\disp}{2}\right)\right] \mu^y(1 - \disp)^y \sum_{j = 0}^{[y/2]} \frac{\disp^j}{2^j\mu^j(1 - \disp)^{2j}(y - 2j)!j!}, y = 0, 1, 2, \dots \label{eq:pmf_hermite2}    
\end{equation}
and $\text{Var}(Y) = (1 + \disp)\mu$. This version is preferred because relevant statements on INAR and INARMA models take a particularly simple form in terms of this parameterization. We note that e.g., Gupta and Jain \cite{Gupta1974} use $d = 1 + \disp$ as the dispersion parameter, which corresponds to the index of dispersion.

In the classical formulation \eqref{eq:pmf_hermite1}, the $\text{Herm}(\lambda_1, \lambda_2)$ distribution can be displayed as a CP distribution by setting
\begin{align}
    N & \sim \text{Pois}(\lambda_1 + \lambda_2)\\
    Y & = \sum_{i = 1}^N Z_i\\
    Z_i & = \begin{cases}
        1 & \text{ with probability } \frac{\lambda_1}{\lambda_1 + \lambda_2}\\
        2 & \text{ with probability } \frac{\lambda_2}{\lambda_1 + \lambda_2}.\\
    \end{cases}
\end{align}
This implies
\begin{equation}
\clustersize = \mathbb{E}(Z_i) = \frac{\lambda_1 + 2\lambda_2}{\lambda_1 + \lambda_2}. \label{eq:theta_herm}
\end{equation}

Returning to our alternative parameterization \eqref{eq:pmf_hermite2}, after some simple algebra we obtain that $\text{Herm}[\text{mean} = \mu, \text{disp} = \disp]$ can equivalently be written as
\begin{align}
    N & \sim \text{Pois}(\mu/\clustersize)\\
    Y & = \sum_{i = 1}^N Z_i\\
    Z_i & = \begin{cases}
        1 & \text{ with probability } 2 - \clustersize\\
        2 & \text{ with probability }\clustersize - 1\\
    \end{cases}
\end{align}
where
$$
\clustersize = \frac{2}{2 - \disp}.
$$

\subsubsection{Hermite INGARCH model}
\label{suppl:hermite_ingarch}

To obtain a Hermite INGARCH model, we set
$$
Z_i = \begin{cases}
        1 & \text{ with probability } 2 - \clustersize\\
        2 & \text{ with probability }\clustersize - 1\\
    \end{cases}
$$
in expression \eqref{eq:N_CP_original}--\eqref{eq:lambda_CP_original}, with $\clustersize$ as specified in equation \eqref{eq:theta_herm}. This implies
$$
X_t \ \mid \ X_{t - 1}, \dots, X_{1 - p}, \lambda_0, \dots, \lambda_{1 - q} \sim \text{Herm}[\text{mean} = \lambda_t, \text{disp} = \disp].
$$
The conditional mean-variance relationship is linear, with
$$
\text{Var}(X_t \ \mid \ X_{t - 1}, \dots, X_{1 - p}, \lambda_0, \dots, \lambda_{1 - q}) = \lambda_t \times \left(1 + \disp\right).
$$
We estimate the parameters $\nu, \alpha, \beta, \disp$ as well as the initial value $\lambda_0$ via maximum likelihood and obtain the estimates for other parameterizations via transformation. The Hermite INARCH(1) model results by constraining $\beta = 0$.

\subsection{The negative binomial distribution and associated models}
\label{suppl:negbin}

\subsubsection{Parameterization and representation as a compound Poisson}
\label{suppl:negbin_compound}

The classical parametrization of the negative binomial distribution is via a size parameter $r$ and a success probability $\pi$. The probability mass function is then given by
$$
\text{Pr}(Y = y) = \binom{y + r - 1}{y} \pi^r (1 - \pi)^y, \ \ \ y = 0, 1, 2, \dots.
$$
For our purposes, a different parameterization often used in GLMs is more suitable. Its two parameters are the mean $\mu$ and a dispersion parameter $\disp$, with
\begin{align*}
    \mu = r(1 - \pi)/\pi,  \ \ \ \ \ \disp = 1/r.
\end{align*}
As indicated in the main manuscript this implies
\begin{equation}
\text{Pr}(Y = y) = \frac{\Gamma(1/\disp + y)}{y!\Gamma(1/\disp)} \left(\frac{1/\disp}{1/\disp + \mu}\right)^{1/\disp} \left(\frac{\mu}{1/\disp + \mu}\right)^y,\label{eq:def_nb_suppl}
\end{equation}
while $\text{Var}(Y) = \mu + \disp\mu^2$.

As described e.g., by \cite{Weiss2018}, the negative binomial distribution can be represented as a compound Poisson distribution as follows. If
\begin{align*}
    N & \sim \text{Pois}[-r\log(\pi)]\\
    Y & = \sum_{i = 1}^N Z_i\\
    Z_i & \stackrel{\textnormal{i.i.d.}}{\sim} \text{Log}(\pi)
\end{align*}
then $Y \sim \text{NegBin}(\text{size} = r, \text{success probability} = \pi)$. Here, $\text{Log}(\pi)$ denotes the logarithmic distribution with probability mass function
$$
\text{Pr}(Z = z) = \frac{(1 - \pi)^z}{-z\log(\pi)}.
$$
Returning to the GLM-like parameterization used in the main manuscript, some back and forth between the different parameterizations yields
\begin{align}
    N & \sim \text{Pois}(\mu/\clustersize)\\
    Y & = \sum_{i = 1}^N Z_i\\
    Z_i & \stackrel{\textnormal{i.i.d.}}{\sim} \text{Log}\left(\frac{1}{1 + \disp^*}\right)
\end{align}
with
$$
\theta = \frac{\disp^*}{\log(1 + \disp^*)}
$$
as an equivalent of
$$
Y \sim \text{NegBin}[\text{mean} = \mu, \text{disp} = \disp^*/\mu].
$$
Note that rather than fixing the dispersion parameter to a value $\disp$ as in \eqref{eq:def_nb_suppl} we introduced a parameter $\disp^*$ which denotes the product of the mean $\mu$ and the dispersion parameter, such that $\text{disp} = \disp^* / \mu$. This will be helpful in the next section.

\subsubsection{Negative binomial INGARCH model}
\label{suppl:negbin_ingarch}

To obtain a negative binomial INGARCH($p, q$) model along the lines of \cite{Xu2012}, we set
$$
\theta = \frac{\disp^*}{\log(1 + \disp^*)}
$$
and
$$
Z_i \stackrel{\textnormal{i.i.d.}}{\sim} \text{Log}\left(\frac{1}{1 + \disp^*}\right)
$$
in expression \eqref{eq:N_CP_original}--\eqref{eq:lambda_CP_original}. This implies
$$
X_t \ \mid \ X_{t - 1}, \dots, X_{1 - p}, \lambda_0, \dots, \lambda_{1 - q} \sim \text{NegBin}[\text{mean} = \lambda_t, \text{disp} = \disp^*/\lambda_t].
$$
The dispersion parameter $\disp/\lambda_t$ of the conditional negative binomial distribution thus depends on $\lambda_t$ (while the mean cluster size $\clustersize$ does not). This implies that
\begin{equation}
\text{Var}(X_t \ \mid \ X_{t - 1}, \dots, X_{1 - p}, \lambda_0, \dots, \lambda_{1 - q}) = \lambda_t \times \left(1 + \disp^*\right),\label{eq:mean_variance_nb}
\end{equation}
i.e., there is again a linear conditional mean variance relationship.

In practice we estimate the parameters $\nu, \alpha, \beta, \disp^*$ as well as the initial value $\lambda_0$ via maximum likelihood and obtain the estimates for other parameterizations via transformation. The negative binomial INARCH(1) model results by constraining $\beta = 0$.

\section{Proofs for Section \ref{sec:ginarma}: A new GINARMA model formulation and properties}
\label{suppl:proofs_general}

\subsection{Alternative display of the Gaussian ARMA($p$, $q$) model} \label{sec:arma}

We show how to display the classical ARMA($p$, $q$) process
\begin{align}
        X_t = \nu \ \ + \ \ \sum_{i = 1}^p \alpha_i \times X_{t - i} + \sum_{j = 1}^q \delta_j \times \import_{t - j} \ \ + \ \ \import_t, \ \ \ \import_t \stackrel{\textnormal{i.i.d.}}{\sim} \textnormal{N}(0, \sigma_{\import}^2)
\end{align}
in as in \eqref{eq:X_arma_reformulation}--\eqref{eq:E_arma_reformulation}. Let us assume without loss of generality that $p = q$. To begin with, we introduce an auxiliary process ${E_t}$
\begin{equation} \label{eq:E_arma_intermediate}
    E_t := \frac{1}{1 + \sum_{j = 1}^p \delta_j} \times \left( \nu + \sum_{i = 1}^p \alpha_i X_{t - i} + \sum_{j = 1}^q \delta_j \varepsilon_{t - j} - \frac{\nu}{1 + \sum_{j = 1}^p \delta_j} \right).
\end{equation}
Consequently, we can express $X_t$ using $E_t$ as
\begin{align}
    X_t & = \nu + \sum_{i = 1}^p \alpha_i X_{t - i} + \sum_{j = 1}^q \delta_j \import_{t - j} + \import_t \nonumber \\
    & = \left(1 + \sum_{j = 1}^p \delta_j\right) \times E_t + \frac{\nu}{1 + \sum_{j = 1}^p \delta_j} + \import_t.\label{eq:X_arma}
\end{align}
Moreover, we can solve \eqref{eq:X_arma} for $\import_t$ and obtain
\begin{equation}\label{eq:epsilon_arma}
    \import_t = X_t - \left(1 + \sum_{j = 1}^p \delta_j\right) \times E_t - \frac{\nu}{1 + \sum_{j = 1}^p \delta_j}.
\end{equation}
Plugging \eqref{eq:epsilon_arma} into \eqref{eq:E_arma_intermediate} we get
\begin{align*}
    E_t & = \frac{1}{1 + \sum_{j = 1}^p \delta_j} \left( \nu + \sum_{i = 1}^p \alpha_i X_{t - i} + \sum_{j = 1}^p \delta_j \import_{t - j} - \frac{\nu}{1 + \sum_{j = 1}^p \delta_j} \right) \\
    & = \frac{1}{1 + \sum_{j = 1}^p \delta_j} \Biggl( \nu + \sum_{i = 1}^p \alpha_i X_{t - i} + \sum_{j = 1}^p \delta_j \underbrace{\left\{X_{t - j} - \left(1 + \sum_{j = 1}^p \delta_j\right) \times E_{t - j} - \frac{\nu}{1 + \sum_{j = 1}^p \delta_j} \right\}}_{\varepsilon_{t - j}} \\
 & \ \ \ \ \ \ \ \ \ \ \ \ \ \ \ \ \ \ \ \ \ \ - \frac{\nu}{1 + \sum_{j = 1}^p \delta_j} \Biggr)\\
 & = \frac{ \sum_{i = 1}^p(\alpha_i + \delta_i) X_{t - i}}{1 + \sum_{j = 1}^p \delta_j} - \sum_{j = 1}^p \delta_j E_{t - j} + \frac{1}{1 + \sum_{j = 1}^p \delta_j}\underbrace{\left( \nu - \frac{\sum_{j = 1}^p \delta_j\nu}{1 + \sum_{j = 1}^p \delta_j} - \frac{\nu}{1 + \sum_{j = 1}^p \delta_j}\right)}_{= 0} \\
 & = \sum_{i = 1}^p \frac{\alpha_i + \delta_i}{1 + \sum_{j = 1}^p \delta_j} \times X_{t - i} - \sum_{j = 1}^p \delta_j E_{t - j}
\end{align*}
Defining $\kappa_i := \frac{\alpha_i + \delta_i}{1 + \sum_{j = 1}^p \delta_j}$, $\beta_j := -\delta_j$ and $\tau := \frac{\nu}{1 - \sum_{j = 1}^p \beta_j}$,
we can display the ARMA($p$, $q$) as

\begin{align*}
X_t & = \left(1 - \sum_{j = 1}^q \beta_j\right)\times E_t + \underbrace{\tau + \import_t}_{:= \import_t'},\\
E_t & = \sum_{j = 1}^q\beta_j E_{t - j} \ \ + \ \ \sum_{i = 1}^p \kappa_i X_{t - i}.
\end{align*}
Note that if we now replace the multiplication operations by corresponding thinning operations from Section \ref{sec:ginarma} and Section \ref{subsec:inarma} and put appropriate constraints on the coefficients $\kappa_i$, $\beta_j$ and $\tau$, we obtain the GINARMA($p$, $q$) model without compounding and the INARMA($p$, $q$) model respectively.

\subsection{Properties of the GINARMA model for $p = q = 1$}
\label{subsec:ginarma_properties}

In the following we will be concerned with the special case $p = q = 1$ of model \eqref{eq:X_general}--\eqref{eq:E_general}. We introduce some additional notation and write this model as
\begin{align}
X_t & = \clustersize * (\import_t + A_t)\label{eq:X_general_11}\\
E_t & = L_{t - 1} \ \ + \ \ C_{t - 1}\label{eq:E_general_11}
\end{align}
with
\begin{align}
    L_t & = \beta \circ E_t\\
    A_t & = (1 - \beta) \circ E_t \stackrel{\text{\eqref{eq:mult_general}}}{=} E_t - L_t \label{eq:A_general_11}\\
    C_t & = \kappa \bullet X_t \label{eq:C_general_11}.
\end{align}
As before, the imports $\{\import_t\}$ are independent samples from an integer-valued distribution.

\subsubsection{Embedded Galton-Watson branching process, lemma \ref{lemma:galton_watson}}
\label{subsubsec:galton_watson_derivation}

We can express $E_t$ as
\begin{align}
    E_t & = L_{t - 1} + \kappa \bullet X_{t - 1}\\
    & \stackrel{\eqref{eq:X_general_11}}{=} L_{t - 1} + \kappa \bullet \{ \clustersize * (\import_{t - 1} + A_{t - 1})\} \nonumber\\
    & = L_{t - 1} + \underbrace{\kappa \bullet ( \clustersize * \import_{t - 1})}_{\text{denote this by } \import^*_{t - 1}} + \kappa \bullet (\clustersize * A_{t - 1}) \nonumber \\
    & \stackrel{\eqref{eq:A_general_11}}{=} \import^*_{t - 1} + L_{t - 1} + \kappa \bullet \{\clustersize * (E_{t - 1} - L_{t - 1})\}.
    \label{eq:galton_watson_intermediate}
\end{align}
Now remember that $L_{t - 1}$ arises from binomial thinning of $E_{t - 1}$ and can be written as
\begin{align*}
L_{t - 1} & = \sum_{k = 1}^{E_t} R_{t - 1, k}, \ \ R_{t - 1, k}  \stackrel{\textnormal{i.i.d.}}{\sim} \text{Bernoulli}(\beta).
\end{align*}
Plugging this back into \eqref{eq:galton_watson_intermediate} we obtain
\begin{align*}
E_t & = \import^*_{t - 1} + \sum_{k = 1}^{E_{t - 1}}R_{t - 1, k} + \kappa \bullet \left\{\clustersize * \left(E_{t - 1} - \sum_{k = 1}^{E_{t - 1}}  R_{t - 1, k}\right) \right\}\\
& = \import^*_{t - 1} + \sum_{k = 1}^{E_{t - 1}}R_{t - 1, k} + \sum_{k = 1}^{E_{t - 1}} \kappa \bullet \{\clustersize * (1 - R_{t - 1, k})\}\\
&  = \import^*_{t - 1} + \sum_{k = 1}^{E_{t - 1}} \Big[R_{t - 1, k}  + \kappa \bullet \{\clustersize * (1 - R_{t - 1, k})\}\Big]\\
& = \import^*_{t - 1} + \sum_{k = 1}^{E_{t - 1}} B_{t - 1, k}
\end{align*}
with
$$
B_{t - 1, k} = \begin{cases}
1 & \text{with probability } \beta \ (\text{i.e., if } R_{t - 1, k} = 1)\\ 
\kappa \bullet (\clustersize * 1) & \text{with probability } 1 - \beta \ (\text{i.e., if } R_{t - 1, k} = 0). 
\end{cases}
$$
This concludes the proof.

\subsubsection{Irreducibility and aperiodicity of $\{E_t\}$, corollary \ref{corollary:aperiodic_irreducible}}

It is immediately visible from equations \eqref{eq:galton_watson}--\eqref{eq:Z_t_i} that $\{E_t\}$ is a time-homogeneous Markov chain. Recall that a time-homogeneous Markov chain $\{Y_t\}$ is called irreducible if all its states communicate, i.e., for each pair of possible states $i$ and $j$, there is an $n \geq 0$ such that
$$
\text{Prob}(Y_{t + n} = j \ \mid \ Y_{t} = i) > 0.
$$
We thus need to prove that for each pair $i, j \in \mathbb{N}_0$ there is an $n \geq 0$ such that 
$$
\text{Prob}(E_{t + n} = j \ \mid \ E_{t} = i) > 0
$$
for model \eqref{eq:galton_watson}--\eqref{eq:Z_t_i}. We prove this by distinction of cases.

\noindent \textbf{Case 1:} If $j \leq i$ we have
$$
\text{Prob}(E_{t + 1} = j \ \mid \ E_{t} = i) > 0
$$
for the following reasons.
\begin{itemize}
    \item There is a positive probability
\begin{equation}
\text{Prob}\left(\sum_{k = 1}^{i} B_{t, k} = j\right) > 0.
\label{eq:prob_all_B}    
\end{equation}
This is the case because there are positive probabilities
$$
p_0 = \text{Pr}(B_{t, k} = 0)
$$
and
$$
p_1 = \text{Pr}(B_{t, k} = 1).
$$
For $p_0$ this follows from $\kappa\theta < 1, \beta < 1$ (there is a positive probability that an exposed cluster will turn infectious, but that none of its members will cause new infections). For $p_1$ it is obvious if $\beta > 0$ (as there is a positive probability that the exposed will just remain exposed). For $\beta = 0$, the technical assumption $\text{Prob}(\kappa_i \bullet 1 = 1) > 0$ made below Definition \ref{def:ingarma} together with $0 < \theta$ ensures that $p_1$ is positive (there is a positive probability that an exposed cluster will turn infectious, that it will have a non-negative number of members and that exactly one member will cause one new infection and the others none). As both $p_0$ and $p_1$ are positive, the probability from \eqref{eq:prob_all_B} is positive with
$$
\text{Prob}\left(\sum_{k = 1}^{i} B_{t, k} = j\right) \geq \binom{i}{j} \times 
p_1^j \times p_0^{i - j} > 0.
$$
\item There is a positive probability
$$
\text{Prob}(\import^*_{t + 1} = 0) > 0.
$$
This holds because irrespective of the value of $\import_{t}$, there is a positive probability that $\kappa \bullet (\clustersize * \import_t) = 0$. This follows from the assumption $\kappa\clustersize < 1$ and the fact that thinning operations cannot yield negative outcomes.
\end{itemize}

\noindent \textbf{Case 2:} If $j > i$ denote $d = j - i$. We now note that for each $i$ we have
\begin{equation}
\text{Prob}(E_{t + 1} \geq i + 1 \ \mid \ E_{t} = i) > 0.\label{eq:one_up}    
\end{equation}
This is because of the two following aspects. Firstly, as in \eqref{eq:prob_all_B} there is a positive probability that $\sum_{k = 1}^{i} B_{t, k} = i$. Secondly, there is a positive probability that $\import^*_t \geq 1$ as $\kappa, \clustersize, \tau > 0$ by assumption. Repeated application of Equation \eqref{eq:one_up} implies that
\begin{equation*}
\text{Prob}(E_{t + d} \geq \underbrace{i + d}_{j} \ \mid \ E_{t} = i) > 0.  
\end{equation*}
There is thus a positive probability to move from $i$ to some $s \geq j$ in $d$ steps. For any such $s \geq j$ there is a positive probability of moving on to $j$ in just one step, see the case $i \geq j$ discussed above. We can thus always transition from $i$ to $j$ via some $s \geq j$ in $d + 1$ steps. This implies irreducibility of $\{E_t\}$. We note that if $\import_t$ can take arbitrarily large values, the argument can be considerably simplified as there is a positive probability of moving from $i$ to some $s \geq j$ in just one step.

As in \cite{Schweer2014} we establish aperiodicity of the irreducible Markov chain by noting that
\begin{equation*}
\text{Prob}(E_{t + 1} = i \ \mid \ E_{t} = i) > 0 \quad \text{for } i = 0, 1, 2, \dotsc.
\end{equation*}
This follows directly from equation \eqref{eq:prob_all_B} with $j = i$. This concludes the proof.

\subsubsection{Limiting-stationary distributions and moments, proposition \ref{proposition:finite_moments}}

To prove the statement we require the following property of randomly stopped sums \cite[Theorem 5.2]{Gut2009}.

\begin{lemma}
    Consider a randomly stopped sum $Y = \sum_{i = 1}^N Z_i$ of i.i.d. random variables $Z_i$ which are independent of $N$. The random variable $Y$ has finite $r$-th moments if $N$ and the $Z_i$ have finite $r$-th moments.\label{lemma:Gut}
\end{lemma}
Provided that it represents an irreducible and aperiodic Markov chain, a sub-critical Galton-Watson branching process with finite import mean has a proper limiting-stationary distribution \cite[Theorem on p.214; note that this condition is sufficient, but could be weakened somewhat]{Heathcote1966}. If the inititial, offspring and immigration distributions have finite $r$-th moments, this is also the case for the limiting-stationary distribution \cite[Sec. 4]{Lange1981}. The process $\{E_t\}$ thus has finite limiting-stationary moments up to order $r$ if the following conditions hold.
\begin{enumerate}
    \item The process $\{E_t\}$ is indeed sub-critical. This is ensured if $\kappa\clustersize <1$ and thus $\mathbb{E}(B_{t, k}) < 1$.
    \item The immigration term $\import^*_t = \kappa \bullet (\clustersize * \import_{t - 1})$ from equation \eqref{eq:galton_watson} has finite moments up to order $r$. Applying Lemma \ref{lemma:Gut} twice, a sufficient condition for this is that $\import_t$, $\clustersize * 1$ and $\kappa \bullet 1$ have finite moments up to order $r$.
    \item $B_{t, k}$ has finite moments up to order $r$. Again we can invoke Lemma \ref{lemma:Gut} twice to show that this is the case if $\clustersize * 1$ and $\kappa \bullet 1$ have finite moments up to order $r$.
    \item The initial value $E_0$ has finite moments up to order $r$. This is of course the case if we initialize the process with fixed values.
\end{enumerate}

\noindent Again under the assumption that $\import_t$, $\clustersize * 1$ and $\kappa \bullet 1$ have finite moments up to order $r$, it is straightforward to show that $A_t = (1 - \beta) \circ E_t$, $\clustersize * A_t$ and ultimately $X_t = \clustersize * (A_t + \import_t)$ likewise have finite limiting-stationary moments up to order $r$. This only requires repeated application of Lemma \ref{lemma:Gut} and concludes the proof.

The above result also implies the existence of higher-order moments of CP-INGARCH(1, 1) models with time-constant $\clustersize$. This had already been proven in \cite{Silva2016}, but the proof is quite involved. Our novel representation allows for a more condensed argument.

\subsubsection{Limiting-stationary means, variances and covariances, lemma~\ref{lemma:general_moments}}

\paragraph{Means}

As demonstrated in Lemma \ref{lemma:galton_watson}, $\{E_t\}$ has a representation as a Galton-Watson branching process with immigration, see equation \eqref{eq:galton_watson}. This also makes it a conditionally linear autoregressive (CLAR) model of order 1 as studied by \cite{Grunwald2000}. Specifically, as
$$
\mathbb{E}(B_{t, k}) = \beta + (1 - \beta)\kappa\clustersize,
$$
compare equation \eqref{eq:Z_t_i}, and
$$
\mathbb{E}(\import^*_t) = \tau\clustersize\kappa
$$
we have
\begin{equation}
\mathbb{E}(E_{t + 1} \ \mid \ E_t) = \tau\clustersize\kappa + E_t \times \{\beta + (1 - \beta)\kappa\clustersize\}.\label{eq:cond_mean_relationship}    
\end{equation}
This implies (\cite{Grunwald2000}, Proposition 1) that if $\kappa \clustersize < 1$ the limiting-stationary mean of $\{E_t\}$ is
$$
\mu_E = \frac{\tau\kappa\clustersize}{1 - \beta - (1 - \beta)\kappa\clustersize}.
$$
For the observable process $\{X_t\}$ we can then compute
\begin{equation}
\mu_X = (1 - \beta)\clustersize\mu_E + \clustersize \tau = \frac{\tau\clustersize}{1 - \kappa\clustersize},\label{eq:CLAR}
\end{equation}
where the simplification in the last step results after some simple algebra.

\paragraph{Variances}

In what follows we will repeatedly use two well-known relationships:
\begin{itemize}
\item If $Y = \sum_{i = 1}^N Z_i$ is a randomly stopped sum of identically and i.i.d. random variables with $Z_1, \dots, Z_N$ independent of $N$, then
\begin{equation}
\text{Var}(Y) = \text{Var}(Z_i)\mathbb{E}(N) +  \text{Var}(N)\mathbb{E}(Z_i)^2.\label{eq:variance_rss}
\end{equation}
\item Consider a generalized thinning operation $\alpha \bullet X = \sum_{i = 1}^X Z_i$ where, independently of $X$, the $Z_i$ are i.i.d. with expectation $\alpha$ and variance $\sigma^2_\alpha < \infty$. If this thinning is performed independently of $Y$ then
\begin{equation}
\text{Cov}(\alpha \bullet X, Y) = \alpha \text{Cov}(X, Y).\label{eq:multiplication_independence}
\end{equation}
\end{itemize}

\noindent To obtain the limiting-stationary variance of $\{E_t\}$, we need to study the conditional variance structure $\text{Var}(E_{t + 1} \ \mid \ E_{t})$ of the process. We first recall from Supplement \ref{subsubsec:galton_watson_derivation} that we can re-write equation \eqref{eq:Z_t_i} as
\begin{align}
R_{t - 1, k} & \stackrel{\textnormal{i.i.d.}}{\sim} \text{Bernoulli}(\beta)\\
B_{t - 1, k} & = \begin{cases}
1 & \text{if } R_{t - 1, k} = 1\\ 
\kappa \bullet (\clustersize * 1) & \text{if } R_{t - 1, k} = 0. 
\end{cases}
\end{align} 
Then we consider
$$
\text{Var}(B_{t, k}) = \text{Var}\{\mathbb{E}(B_{t, k} \ \mid R_{t, k})\} \ + \ \mathbb{E}\{\text{Var}(B_{t, k} \ \mid R_{t, k})\}.
$$
We treat the two summands separately, in a first step
\begin{align*}
\text{Var}\{\mathbb{E}(B_{t, k} \ \mid R_{t, k})\} & = \mathbb{E}\{\mathbb{E}(B_{t, k} \ \mid R_{t, k})^2\} - \mathbb{E}\{\mathbb{E}(B_{t, k} \ \mid R_{t, k})\}^2\\
& = \beta + (1 - \beta)\kappa^2\clustersize^2 - \{\beta + (1 - \beta)\kappa\clustersize\}^2\\
& [...]\\
& = \beta(1 - \beta)(1 - \kappa\clustersize)^2.
\end{align*}
In a second step and using relationship \eqref{eq:variance_rss} we obtain
\begin{align*}
\mathbb{E}\{\text{Var}(B_{t, k} \ \mid R_{t, k})\} & = \beta \times 0 + (1 - \beta)\times \text{Var}(\kappa \bullet \clustersize * 1)\\
& = (1 - \beta) \times (\sigma^2_\kappa\clustersize + \sigma^2_\clustersize\kappa^2).
\end{align*}
Bringing the two summands back together we then obtain
\begin{align*}
\text{Var}(B_{t, k}) & = \beta(1 - \beta)(1 - \kappa\clustersize)^2 + (1 - \beta) \times \{\sigma^2_\kappa\clustersize + \sigma^2_\clustersize\kappa^2\}\\
& = (1 - \beta) \times \{\beta (1 - \kappa\clustersize)^2 + \sigma^2_\kappa\clustersize + \sigma^2_\clustersize\kappa^2\}.
\end{align*}
Turning to the import distribution and using relationship \eqref{eq:variance_rss} twice, we moreover show that
\begin{align*}
\text{Var}(\import^*_t) & = \text{Var}(\kappa \bullet \clustersize * \import_{t - 1})\\
& = \sigma^2_\kappa \mathbb{E}(\clustersize * \import_{t - 1}) + \text{Var}(\clustersize * \import_{t - 1}) \times \kappa^2\\
& = \sigma^2_\kappa\clustersize\tau + (\sigma^2_\clustersize\tau + \sigma^2_\tau\clustersize^2)\times\kappa^2.
\end{align*}
We now have all the necessary pieces to write out the conditional variance structure as
\begin{align*}
\text{Var}(E_{t + 1} \ \mid \ E_{t}) & = \text{Var}(\import^*_t) + E_t \times \text{Var}(B_{t, k})\\
& = \sigma^2_\kappa\clustersize\tau + (\sigma^2_\clustersize \tau + \sigma^2_\tau)\kappa^2 + E_t \times (1 - \beta) \times \{\beta \times (1 - \kappa\clustersize)^2 + \sigma^2_\kappa\clustersize + \sigma^2_\clustersize\kappa^2\}.
\end{align*}
Using Proposition 2 from \cite{Grunwald2000} we then conclude that if $\kappa \clustersize < 1$
$$
\sigma^2_E = \frac{\sigma^2_\kappa\clustersize \tau + (\sigma^2_\clustersize\tau + \sigma^2_\tau\clustersize^2)\times\kappa^2 + \mu_E \times (1 - \beta) \times \{\beta (1 - \kappa\clustersize)^2 + \sigma^2_\kappa\clustersize + \sigma^2_\clustersize\kappa^2\}}{1 - \{\beta + (1 - \beta)\kappa\clustersize\}^2}.
$$
Now we turn to the variance of $\{X_t\}$, which can be obtained as
\begin{align*}
\text{Var}(X_t) & = \text{Var}\{\clustersize * [(1 - \beta) \circ E_t]\} + \text{Var}(\clustersize * \import_t).
\end{align*}
Considering again the two summands separately we obtain
\begin{align*}
\text{Var}\{\clustersize * (1 - \beta) \circ E_t\} & = \text{Var}\{(1 - \beta) \circ E_t\}\clustersize^2 + \sigma^2_\clustersize\mathbb{E}\{(1 - \beta) \circ E_t\}\\
& \stackrel{\eqref{eq:variance_rss}}{=} \{\beta(1 - \beta)\mu_E + \sigma^2_E(1 - \beta)^2\}\clustersize^2 + \sigma^2_\clustersize(1- \beta)\mu_E\\
& = (1 - \beta)\mu_E \sigma^2_\clustersize + \clustersize^2(1 - \beta)\{\beta\mu_E + (1 - \beta)\sigma^2_E\}
\end{align*}
and
\begin{align*}
\text{Var}(\clustersize * \import_t) & \stackrel{\eqref{eq:variance_rss}}{=} \sigma^2_\clustersize\tau + \sigma^2_\tau\clustersize^2,
\end{align*}
which in result gives us
$$
\text{Var}(X_t) = (1 - \beta)\mu_E \sigma^2_\clustersize + \clustersize^2(1 - \beta)\{\beta\mu_E + (1 - \beta)\sigma^2_E\} + \sigma^2_\clustersize\tau + \sigma^2_\tau\clustersize^2.
$$

\paragraph{Autocovariances}

From the CLAR(1) representation \eqref{eq:CLAR} of $\{E_t\}$ it follows that (\cite{Grunwald2000}, Proposition 4)
$$
\gamma_E(d) = \{\beta + (1 - \beta)\kappa\clustersize\}^d \times \sigma^2_E.
$$
For the autocovariance structure of $\{X_t\}$ we use the notation introduced at the beginning of Section \ref{subsec:ginarma_properties} and consider
\begin{align*}
\text{Cov}(X_t, E_{t + 1}) & = \text{Cov}(X_t, L_t + \kappa \bullet X_t)\\
& = \text{Cov}(X_t, L_t) + \text{Cov}(X_t, \kappa \bullet X_t)\\
& = \text{Cov}(\clustersize * A_t + \clustersize * \import_t, L_t) + \text{Cov}(X_t, \kappa \bullet X_t)\\
& = \text{Cov}(\clustersize * A_t, L_t) \ + \ \underbrace{\text{Cov}(\clustersize * \import_t, L_t)}_{= \ 0} \ + \ \text{Cov}(X_t, \kappa \bullet X_t).
\end{align*}
Considering the two non-zero summands separately, we get
\begin{align}
\text{Cov}(\clustersize * A_t, L_t) & \stackrel{\eqref{eq:multiplication_independence}}{=} \mathbb{E}\{\text{Cov}(\clustersize * A_t, L_t \ \mid \ E_t)\} + \text{Cov}\{\mathbb{E}(\clustersize * A_t \ \mid \ E_t), \mathbb{E}(L_t \ \mid \ E_t)\} \nonumber\\
& = \clustersize \mathbb{E}\{\underbrace{\text{Cov}(A_t, L_t \ \mid \ E_t)}_{\text{note that } A_t + L_t = E_t}\} \ + \ \text{Cov}\{(1 - \beta)\clustersize E_t, \beta E_t\} \nonumber\\
& = \clustersize \mathbb{E}\{\text{Cov}(E_t - L_t, L_t \ \mid \ E_t)\} \ + \ (1 - \beta)\clustersize\beta \text{Cov}\{E_t, E_t\} \nonumber\\
& = -\clustersize \mathbb{E}\{\underbrace{\text{Cov}(L_t, L_t \ \mid \ E_t)}_{L_t \mid E_t \sim \text{Bin}(E_t, \beta)}\} \ + \ (1 - \beta)\clustersize\beta \text{Cov}\{E_t, E_t\} \nonumber\\
& = - \clustersize \beta(1 - \beta)\mu_E \ + \ (1 - \beta)\clustersize\beta \sigma^2_E \nonumber\\
& = \beta(1 - \beta)\clustersize(\sigma^2_E - \mu_E)\label{eq:intermediate_step}
\end{align}
and
\begin{align*}
\text{Cov}(X_t, \kappa \bullet X_t) & = \underbrace{\mathbb{E}\{\text{Cov}(X_t, \kappa \bullet X_t \ \mid \ X_t)\}}_{= \ 0} \ + \ \text{Cov}\{\mathbb{E}(X_t \ \mid \ X_t), \mathbb{E}(\kappa \bullet X_t \ \mid \ X_t)\}\\
& = \text{Cov}(X_t, \kappa \bullet X_t) = \kappa\sigma^2_X.
\end{align*}
Putting these back together results in
$$
\text{Cov}(X_t, E_{t + 1}) = \clustersize\beta(1 - \beta)(\sigma^2_E - \mu_E) + \kappa\sigma^2_X.
$$
For $d = 2, 3, \dots$ we can now consider
\begin{align}
\text{Cov}(X_t, E_{t + d}) & = \underbrace{\mathbb{E}\{\text{Cov}(X_t, E_{t + d} \ \mid \ X_t, E_{t + d - 1})\}}_{= \ 0} \ + \ \text{Cov}\{\mathbb{E}(X_t \ \mid \ X_t, E_{t + d - 1}), \mathbb{E}(E_{t + d} \ \mid \ X_t, E_{t + d - 1})\}\nonumber\\
& \stackrel{\eqref{eq:cond_mean_relationship}}{=} \text{Cov}[X_t, \{\beta + (1 - \beta)\kappa\clustersize\} E_{t + d - 1}\}]\nonumber\\
& = \{\beta + (1 - \beta)\kappa\clustersize\}\times\text{Cov}(X_t, E_{t + d - 1})\nonumber\\
& = \{\beta + (1 - \beta)\kappa\clustersize\}^{d - 1}\times\text{Cov}(X_t, E_{t + 1}).\label{eq:intermediate_step2}
\end{align}
Finally, combining \eqref{eq:intermediate_step} and \eqref{eq:intermediate_step2}, we note that for $d = 1, 2, \dots$
\begin{align*}
\text{Cov}(X_t, X_{t + d}) & = \underbrace{\mathbb{E}\{\text{Cov}(X_t, X_{t + d} \ \mid \ X_t, E_{t + d})\}}_{= \ 0} \ + \ \text{Cov}\{\mathbb{E}(X_t \ \mid \ X_t, E_{t + d}), \mathbb{E}(X_{t + d} \ \mid \ X_t, E_{t + d})\}\\
& = \text{Cov}[X_t, (1 - \beta) \clustersize E_{t + d} + \tau\clustersize\}]\\
& = \text{Cov}[X_t, (1 - \beta) \clustersize E_{t + d}\}]\\
& = (1 - \beta)\clustersize\times \{\beta + (1 - \beta)\kappa\clustersize\}^{d - 1} \times \{\clustersize\beta(1 - \beta)(\sigma^2_E - \mu_E) + \kappa\sigma^2_X\}.
\end{align*}

\subsubsection{Geometric ergodicity, proposition \ref{proposition:geometric_ergodicity}}

We use again the representation \eqref{eq:galton_watson}--\eqref{eq:Z_t_i} of $\{E_t\}$ as a Galton-Watson branching process with immigration. Theory on branching processes with immigration, specifically Theorem 1 from Pakes \cite{Pakes1971} tells us that $\{\juv_t\}$ is geometrically ergodic if (i) $\mathbb{E}(B_{t, k}) < 1$, (ii) $\mathbb{E}[B_{t, k} \times \log(B_{t, k}) \ \mid \ B_{t, k} \geq 1] < \infty$, (iii) $\mathbb{E}(\import^*_t) < \infty$. These conditions are easily verified for $\{E_t\}$ provided that $\kappa\clustersize < 1$ and, as previously assumed, $\sigma^2_\tau, \sigma^2_\kappa, \sigma^2_\clustersize < \infty$. Note that in \cite[Assumption 1]{Pakes1971} there are additional technical conditions (iv) $0 < \text{Pr}(B_{t, k} = 0) < 1$, (v) $0 < \text{Pr}(B_{t, k} \leq 1) < 1$ and (vi) $0 < \text{Pr}(\varepsilon^*_{t} = 0) < 1$. Condition (iv) is implied by $\beta < 1$ and $0 < \kappa\clustersize < 1$. Condition (vi) is implied by $\tau > 0$ and $0 < \kappa\clustersize < 1$. Concerning condition (v), Schweer and Weiss \cite[footnote 3]{Schweer2014} have remarked that it is not actually required for the proof of Theorem 1 from Pakes \cite{Pakes1971}.

As in Fokianos et al \cite{Fokianos2009}, Proposition 1 from Meitz and Saikkonen \cite{Meitz2008} can then be used to show that geometric ergodicity of $\{E_t\}$ is inherited by the joint process $\{(\juv_t, L_t, A_t, X_t, C_t)\}$, see definitions in Supplementary Section \ref{subsec:ginarma_properties}. Even though it is in principle sufficient to initialize the process with $E_0$ and $X_0$ as in Section \ref{sec:ginarma}, we now assume that $\{(\juv_t, L_t, A_t, X_t, C_t)\}$ is initialized by a vector $(e_0, l_0, a_0, x_0, c_0)$ with all elements from $\mathbb{N}_0$ and $l_0 + a_0 = e_0, x_0 \geq a_0$. Geometric ergodicity of the joint process is then established by verifying two conditions (Assumption 1 in \cite{Meitz2008}):

\begin{enumerate}
\item Given $E_t$, $(L_t, A_t, X_t, C_t)$ is independent of all $\juv_u, L_u, A_u, X_u, C_u, u < t$. It is straightforward to see from equations \eqref{eq:X_general_11}--\eqref{eq:C_general_11} or the graphical representation in Figure \ref{fig:ingarch_flowchart_poisson} that this is the case.
\item There is an $n \geq 1$ such that for all $t > n$, the generation mechanism of $\juv_t \ \mid \ \juv_0 = e_0, L_0 = l_0, A_0 = a_0, X_0 = x_0, C_0 = c_0$ has the same structure as that of $\juv_t \ \mid \ \juv_n = \tilde{e}_n$, where $\tilde{e}_n$ is some function of $(e_0, l_0, a_0, x_0, c_0)$. As $(e_0, l_0, a_0, x_0, c_0)$ only impacts the further course of the process $\{\juv_t\}$ through $\juv_1 = c_0 + l_0$, this is the case for $n = 1$.
\end{enumerate}
This concludes the proof of geometric ergodicity. As noted by \cite[directly below their Proposition 1]{Meitz2008}, the joint process is moreover $\beta$-mixing with exponentially decreasing weights if it is initialized with its stationary distribution (in our case this means $E_0$ must be assigned the respective stationary distribution, and $L_0, A_0, X_0, C_0$ must be sampled based on $E_0$).

\subsection{Proofs for Subsection \ref{sec:alternative_formulation}: Thinning-based representation of CP-INGARCH models}
\label{appendix:proofs_equivalence}

In this section we provide the derivations of the alternative thinning-based representation of various INGARCH models. We will use the language of Subsection \ref{subsec:interpretation_epidemic} to facilitate the verbal description.

\subsubsection{Poisson INGARCH(1, 1)}
\label{subsubsec:derivation_poisson11}

We demonstrate that the process $\{X_t, t \in \mathbb{N}\}$ from \eqref{eq:X_t_thinning_Poisson} is equivalent to the Poisson INGARCH(1, 1) process \eqref{eq:ingarch11}. We start by writing out the thinning-based representation with some auxiliary processes analogously to Subsection \ref{subsec:ginarma_properties}. We specify
\begin{align}
    X_t & = A_t + \import_t \label{eq:X_thinned_suppl}\\
    E_t & = L_t + C_t,
\end{align}
where
\begin{align*}
    L_t & = \beta \circ E_t\\
    A_t & = (1 - \beta) \circ E_t \stackrel{\text{\eqref{eq:mult_general}}}{=} E_t - L_t \\
    C_t & = \kappa \star X_t\\
    \import_t & \stackrel{\textnormal{i.i.d.}}{\sim} \text{Pois}(\tau).
\end{align*}
Remember that for initialization we have
$$
E_0 \sim \text{Pois}(\eta).
$$
We start by decomposing $C_t, t \in \mathbb{N}$ and $E_0$ by when these individuals will become infectious, i.e., will transition from $E$ to $X$. We denote by $C_t^{(i)}$ the number of exposed persons caused by infectives from time $t$ and turning themselves infectious at $t + i$; and by $E^{(t)}_0$ the number of exposed individuals initially in the pool and turning infectious at time $t$. This implies
\begin{equation}
A_t = \sum_{i = 1}^{t} C_{t - i}^{(i)} \ \ + \ \ E_0^{(t)}
\label{eq:sums}
\end{equation}
for $ t = 1, 2, \dots$. A person infected by an infective from time $t$ (i.e., entering the exposed pool at time $t + 1$ via $C_t$) has a probability of 
\begin{equation}
\beta^{i - 1}(1 - \beta)\label{eq:geom_distr}
\end{equation}
to become infectious at time $t + i, i = 1, 2, \dots$, and thus be part of $C_t^{(i)}$ (it has to remain in the exposed pool $i - 1$ times and then turn infectious). The Poisson splitting property \cite{Kingman1993} implies that given $X_t$, the $C_t^{(i)}, i =1, 2, \dots$ are independently Poisson distributed,
$$
C_t^{(i)} \mid X_t \stackrel{\text{ind}}{\sim} \text{Pois}(\beta^{i - 1}[1 - \beta]\kappa X_t), i =1, 2, \dots
$$
We note that given $X_t$, $C_t^{(i)}$ does not have any impact on the further course of the process $\{X_t\}$ until time $t + i$. Also, given $X_t$, $C_t^{(i)}$ is independent of all preceding values $X_{t - 1}, X_{t - 2}, X_0$. We can thus extend the condition in the above and write
\begin{equation}
C_t^{(i)} \mid X_{t + i - 1}, \dots, X_0 \sim \text{Pois}(\beta^{i - 1}[1 - \beta]\kappa X_t). \label{eq:conditional_Cti}
\end{equation}
Now consider
\begin{align}
X_t & = \import_t \ \ + \ \ \underbrace{\sum_{i = 1}^{t} C_{t - i}^{(i)} \ \ + \ \ E_0^{(t)}}_{= A_t}, \label{eq:decomposition_Xt}
\end{align}
where we substituted $A_t$ in equation \eqref{eq:X_t_thinning_Poisson} using equation \eqref{eq:sums}. Because, given $X_{t - 1}, \dots, X_0$, the $C_{t - i}^{(i)}, i  =1, \dots, t$ only impact the further process from $t$ onwards, it is clear that they are all conditionally independent. The same holds for $E_0^{(t)}$, which is Poisson distributed with rate $\beta^{t}(1 - \beta)\eta$. Conditioned on $X_{t - 1}, \dots, X_0$, we thus have that $X_t$ is a sum of independent Poisson random variables. We can therefore write
$$
X_t \mid X_{t - 1}, \dots, X_0 \sim \text{Pois}(\lambda_t)
$$
where the conditional expectation is given by
\begin{align*}
& \lambda_t = \mathbb{E}(\import_t) \ \ + \ \ \sum_{i = 1}^t \mathbb{E}(C_{t - i}^{(i)} \ \mid \ X_{t- 1}, \dots, X_0) \ \ + \ \ \mathbb{E}(E_0^{(t)})\\
& \ \ \ = \ \ \tau \ \ \ \ \ + \ \ \sum_{i = 1}^t \beta^{i - 1}(1 - \beta)\kappa X_{t - i} \ \ \ \ \ \ \ \ \ \ + \ \ \beta^{t}(1 - \beta)\eta.
\end{align*}
We can then re-write $\lambda_t$ as
\begin{align*}
\lambda_t & = (1 - \beta)\tau + (1 - \beta)\kappa X_{t - 1} + \beta \left\{\tau +    \sum_{i = 2}^t \beta^{i - 1}(1 - \beta)\kappa X_{t - i}  \ \ + \ \ \beta^{t - 1}(1 - \beta)\eta\right\}\\
& = \underbrace{(1 - \beta)\tau}_{\nu} \ + \ \underbrace{(1 - \beta)\kappa}_{\alpha} X_{t - 1} \ + \ \beta \lambda_{t - 1}
\end{align*}
for $t \geq 2$. This is the form a Poisson INGARCH(1, 1) model with parameters $\nu = (1 - \beta)\tau, \alpha = (1 - \beta)\kappa$ and $\beta$. We conclude by considering the initialization of the process, where we have
\begin{align*}
\lambda_1 = \mathbb{E}(X_1 \ \mid \ X_0) & = \tau + (1 - \beta)\kappa X_0 + \beta(1 - \beta)\eta\\
& =  (1 - \beta)\tau + (1 - \beta)\kappa X_0 + \beta  \left\{\tau + (1 - \beta) \times \eta \right\},
\end{align*}
meaning that we have to set $\lambda_0 =\tau + (1 - \beta) \times \eta$ for initialization. The formulas $\tau = \nu/(1 - \beta), \kappa = \alpha/(1 - \beta)$ and $\eta = (\lambda_0 - \tau)/(1 - \beta)$ from the main manuscript result from solving the respective equations for the parameters of the thinning-based parameterization.

\subsubsection{Poisson INGARCH($p, q$)}
\label{subsubsec:poissonpq_proof}

We use an argument similar to the one from Subsection \ref{subsubsec:derivation_poisson11} to demonstrate that the Poisson INGARCH($p, q$) model 
\begin{align}
X_t \ \mid \ X_{t - 1}, \dots, X_{1 - p}, \lambda_0, \dots, \lambda_{1 - q} & \sim \text{Pois}(\lambda_t)\\
\lambda_t & = \nu + \sum_{i = 1}^p \alpha_i X_{t - i} + \sum_{j = 1}^q \beta_j \lambda_{t - j}
\end{align}
and the thinning-based formulation \eqref{eq:Xt_thinning_pq} are equivalent. We re-write the latter as
\begin{align*}
    X_t & = A_t + \import_t\\
    E_t & = \sum_{i = 1}^p C_{t - i, i} \ + \ \sum_{j = 1}^q L_{t - j, j}
\end{align*}
with
\begin{align*}
    C_{t, i} & = \kappa_i \star X_t\\
    (L_{t, 1}, \dots, L_{t, q}, A_t) \ \mid \ E_t  & \sim \text{Mult}\left[E_t, \beta_1, \dots, \beta_q, 1 - \sum_{j = 1}^q \beta_j \right] \\
    \import_t & \stackrel{\textnormal{i.i.d.}}{\sim} \text{Pois}(\tau).
\end{align*}
For the initialization we fix $X_{1 - p}, \dots, X_0$ and set $E_{m} \stackrel{\text{ind}}{\sim} \text{Pois}(\eta_m)$ with $\eta_m > 0$ for $m = 1 - q, \dots, 0$.

Again we denote by $C_t^{(i)}, i = 1, 2, \dots$ the number of persons infected by infectives from time $t$ and becoming themselves infectious at $t + i$. Extending on the notation from the Poisson INGARCH(1, 1) case, we denote by $E^{(i)}_m, m = 1 - q, \dots, 0, i = 1, 2, \dots$ the number of individuals entering the exposed pool via the initialization at time $m$ and turning infectious at time $m + i$. Generalizing equation \eqref{eq:decomposition_Xt} we then have
$$
X_t = \import_t \ \ + \ \ \underbrace{\sum_{i = 1}^{t + p - 1} C_{t - i}^{(i)} \ \ + \ \ \sum_{m = 1 - q}^0 E_m^{(t - m)}}_{= A_t}
$$
for $t = 1, 2, \dots$. Arguments identical to those from the previous section imply that given $X_{t - 1}, \dots, X_{1 - p}$ all summands in the above equation are independently Poisson distributed, so that $X_t$, too, is conditionally Poisson with a rate $\lambda_t$.

Paralleling equation \eqref{eq:conditional_Cti}, the conditional expectation of $C_t^{(i)}$ is given by
\begin{equation}
\mathbb{E}(C_t^{(i)} \ \mid \ X_{t + i - 1}, \dots, X_{1 - p}) = \left(1 - \sum_{l = 1}^q \beta_l \right) \times \left(\sum_{k = 1}^p\ \kappa_k \pi_{i - k}\right) X_t,\label{eq:ELt}
\end{equation}
where we denote by $\pi_j$ the probability that an individual entering the exposed pool at time $t$ is also in the pool at time $t + j$. The reasoning behind this relationship is that the $X_t$ infectives from time $t$ generate exposures entering at times $t + 1, \dots, t + p$ with rates $\kappa_1, \dots, \kappa_p$, respectively. The exposed individuals then have to also be present in the exposed pool exactly $i - 1, \dots, i - p$ time points later, respectively (which happens with probabilities $\pi_{i - 1}, \dots, \pi_{i - p}$), and then leave it (which happens with probability $1 - \sum_{l = 1}^q \beta_l$).

For the $\pi_j$, the recursion
\begin{align}
\pi_j = \sum_{l = 1}^q \beta_l \pi_{j - l}, \label{eq:recursion_pi}
\end{align}
with $\pi_0 = 1$ and $\pi_k = 0$ for $k < 0$ holds. This is because an individual which entered the exposed pool at time $t$ can arrive in $E_{t + j}, j \geq 1$ by a move from any of $E_{t + j - 1}, \dots E_{t + j - q}$ (even though some of these moves may not be possible if $j < q$; this will be reflected in $\pi_{j - l} = 0$). To do so, the individual needs to have arrived at the respective $E_{t + j - l}$ (which it does with probability $\pi_{j - l}$) and then make an $l$-step jump into $E_{t + j}$ (this happens with probability $\beta_l$).

We can now consider
\begin{align}
\lambda_t = \mathbb{E}(X_t \ \mid \ X_{t - 1}, \dots, X_{1 - p}) = & \ \tau 
\ + \ \sum_{i = 1}^{t + p - 1}\mathbb{E}(C^{(i)}_{t - i} \ \mid \ X_{t - 1}, \dots, X_{1 - p})
\ \label{eq:lambda_t_pq_recursion1}\\
& \ \ + \sum_{m = 1 - q}^0 \mathbb{E}(E_{m}^{(t - m)}  \ \mid \ X_{t - 1}, \dots, X_{1 - p}).\nonumber
\end{align}

Focusing on the second summand and plugging in equation \eqref{eq:ELt}, we obtain
\begin{align*}
\sum_{i = 1}^{t + p - 1}\mathbb{E}(C^{(i)}_{t - i} \ \mid \ X_{t - 1}, \dots, X_{1 - p}) = & \left(1 - \sum_{l = 1}^q \beta_l\right) \times \left( \sum_{i = 1}^{t + p - 1} \sum_{k = 1}^p \kappa_k \pi_{i - k} X_{t - i}\right)\\
= & \left(1 - \sum_{l = 1}^q \beta_l\right) \times \left(\sum_{k = 1}^p \sum_{i = k}^{t + p - 1} \kappa_k \pi_{i - k} X_{t - i}\right).
\end{align*}
Note that in the last step we can start the last sum from $i = k$ rather than $i = 1$ as $\pi_{i - k} = 0$ for $i < k$. We can then further decompose this sum into
\begin{align}
= & \left(1 - \sum_{l = 1}^q \beta_l\right) \times \left(\underbrace{\sum_{k = 1}^p \kappa_k \pi_0 X_{t - k}}_{\text{corresponds to } i = k; \text{ note: } \pi_0 = 1} \ + \ \sum_{k = 1}^p \sum_{i = k + 1}^{t + p - 1} \kappa_k \pi_{i - k} X_{t - i}\right)\nonumber\\
= & \left(1 - \sum_{l = 1}^q \beta_l\right) \times \left\{\sum_{k = 1}^p \kappa_k X_{t - k} \ + \ \sum_{k = 1}^p \sum_{i = k + 1}^{t + p - 1} \kappa_k \times \underbrace{\left(\sum_{j = 1}^q \beta_j \pi_{i - k - j}\right)}_{\text{using equation \eqref{eq:recursion_pi}}} \times X_{t - i}\right\}\nonumber\\
= & \left(1 - \sum_{l = 1}^q \beta_l\right) \times \left\{\sum_{k = 1}^p \kappa_k X_{t - k} \ + \ \sum_{j = 1}^q \beta_j \times \left(\sum_{k = 1}^p \sum_{i = k + 1}^{t + p - 1} \kappa_k \times \pi_{i - k - j} \times X_{t - i}\right)\right\}\nonumber\\
= & \left(1 - \sum_{l = 1}^q \beta_l\right) \times \left\{\sum_{k = 1}^p \kappa_k X_{t - k} \ + \ \sum_{j = 1}^q \beta_j \times \left(\sum_{k = 1}^p \sum_{i = k + 1 - j}^{(t - j) + p - 1} \kappa_k \times \pi_{i - k} \times X_{(t - j) - i}\right)\right\}\nonumber\\
= & \left(1 - \sum_{l = 1}^q \beta_l\right) \times \left\{\sum_{k = 1}^p \kappa_k X_{t - k} \ + \ \sum_{j = 1}^q \beta_j \times \left(\sum_{k = 1}^p \sum_{i = 1}^{(t - j) + p - 1} \kappa_k \times \pi_{i - k} \times X_{(t - j) - i}\right)\right\}\label{eq:substitution_ELt}
\end{align}
where in the last step we can let the last sum start at $i = 1$ rather than $i = k + 1 - j$ as $\pi_{i - k} = 0$ for $i = 1, \dots, k - j$.

For the third term from equation \eqref{eq:lambda_t_pq_recursion1} we pursue a similar recursive argument:
\begin{align}
\sum_{m = 1 - q}^0 \mathbb{E}(E_{m}^{(t - m)}  \ \mid \ X_{t - 1}, \dots, X_{1 - p}) & = \left(1 - \sum_{l = 1}^q \beta_l\right) \times \sum_{m = 1 - q}^0 \pi_{t - m}\eta_m\nonumber\\
& = \left(1 - \sum_{l = 1}^q \beta_l\right) \times \sum_{m = 1 - q}^0 \sum_{j = 1}^q \beta_j \pi_{t - j - m}\eta_m\nonumber\\
& =  \left(1 - \sum_{l = 1}^q \beta_l\right) \times\sum_{j = 1}^q \beta_j \times \left(\sum_{m = 1 - q}^0 \pi_{(t - j) - m}\eta_m\right).\label{eq:substitution_ESt}
\end{align}
Plugging the terms from \eqref{eq:substitution_ELt} and \eqref{eq:substitution_ESt} into \eqref{eq:lambda_t_pq_recursion1} we then get

\begin{align*}
\lambda_t = & \ \ \tau \ + \ \left(1 - \sum_{l = 1}^q \beta_l\right) \times \Bigg\{\sum_{k = 1}^p \kappa_k X_{t - k} \ + \ \sum_{j = 1}^q \beta_j \times \left(\sum_{k = 1}^p \sum_{i = 1}^{(t - j) + p - 1} \kappa_k \times \pi_{i - k} \times X_{(t - j) - i}\right) \\
& \ \ \ \ \ \ \ \ \ \ \ \ \ \ \ \ \ \ \ \ \ \  \ \ \ \ \ \ \ \ \ \ \ + \ \sum_{j = 1}^q \beta_j \times \left(\sum_{m = 1 - q}^0 \pi_{(t - j) - m}\eta_m\right)\Bigg\}.
\end{align*}
This can be re-ordered to
\begin{align*}
\lambda_t & = \left(1 - \sum_{l = 1}^q \beta_l\right) \times \tau \ +  \ \left(1 - \sum_{l = 1}^q \beta_l\right) \times \left(\sum_{k = 1}^p \kappa_k X_{t - k}\right)\\
& \ \ \ + \sum_{j = 1}^q \beta_j \times \Biggl\{\tau \ + \ \left(1 - \sum_{l = 1}^q \beta_l\right) \times \left(\sum_{k = 1}^p \sum_{i = 1}^{(t - j) + p - 1} \kappa_k \times \pi_{i - k} \times X_{(t - j) - i}\right) \  \\
&  \ \ \ \ \ \ \ \ \ \ \ \ \ \ \ \ \ \ \underbrace{\ \ \ \ \ \ \ \  + \left(1 - \sum_{l = 1}^q \beta_l\right) \times \left(\sum_{m = 1 - q}^0 \pi_{(t - j) - m}\eta_m \right) \Biggr\} \ \ \ \ \ \ \ \ \ \ \ \ \ \ \ \ \ \ }_{= \lambda_{t - j}}\\
& = \nu \ \ + \ \ \sum_{k = 1}^p \alpha_k X_{t - k} \ \ + \ \ \sum_{j = 1}^q \beta_j \lambda_{t - j},
\end{align*}
where
$$
\nu = \left(1 - \sum_{l = 1}^q \beta_l\right) \times \tau, \ \ \ \alpha_k = \left(1 - \sum_{l = 1}^q \beta_l\right) \times \kappa_k, \ \ \ k = 1, \dots, p.
$$
This is the form of a Poisson INGARCH($p, q$) model as defined in equations \eqref{eq:N_CP_original}--\eqref{eq:lambda_CP_original} (omitting the compounding step). Concerning the initialization, it can be shown that one needs to set $\lambda_m = \tau + (1 - \sum_{j = 1}^q\beta_j) \times \eta_m, m = 1 - q, \dots, 0$. This can be done using essentially the same argument as in Section \ref{subsubsec:derivation_poisson11}, but we omit the somewhat lengthy details. The equations
\begin{align*}
    \tau & = \nu / (1 - \sum_{j = 1}^q\beta_j)\\
    \kappa_i & = \alpha_i / (1 - \sum_{j = 1}^q\beta_j), i = 1, \dots, p\\
    \eta_m & = (\lambda_m - \tau )/(1 - \sum_{j = 1}^q\beta_j), m = 1 - q, \dots, 0
\end{align*}
provided in the manuscript result again from simply solving the relationships between the two parameterizations for the respective parameters.

\subsubsection{Compound Poisson INGARCH(1, 1)}
\label{subsubsec:compound_proof}

The thinning-based formulation of the compound Poisson INGARCH(1, 1) results from replacing \eqref{eq:X_thinned_suppl} by
$$
\clustersize * (A_t + \import_t).
$$
Setting $N_t = \import_t + A_t$, the same arguments as in Section \ref{subsubsec:derivation_poisson11} can be used to show that
$$
N_t \mid X_{t - 1}, \dots, X_0 \sim \text{Pois}(\lambda_t/\clustersize)
$$
where the conditional expectation is given by
\begin{align*}
& \lambda_t/\clustersize = \mathbb{E}(\import_t) \ \ + \ \ \sum_{i = 1}^t \mathbb{E}(C_{t - i}^{(i)}) \ \ + \ \ \mathbb{E}(E_0^{(t)})\\
& \ \ \ \ \ \ \ \ = \tau \ \ + \ \ \sum_{i = 1}^t \beta^{i - 1}(1 - \beta)\kappa X_{t - i} \ \ + \ \ \beta^{t}(1 - \beta)\eta\\
\Leftrightarrow \ \ & \lambda_t = \clustersize \tau \ \ + \ \ \clustersize \times \sum_{i = 1}^t \beta^{i - 1}(1 - \beta)\kappa X_{t - i} \ \ + \ \ \clustersize\beta^{t}(1 - \beta)\eta
\end{align*}
We can then re-write $\lambda_t$ as
\begin{align*}
\lambda_t & = \clustersize(1 - \beta)\tau + \clustersize(1 - \beta)\kappa X_{t - 1} + \beta \left[\clustersize\tau +   \clustersize \sum_{i = 2}^t \beta^{i - 1}(1 - \beta)\kappa X_{t - i}  \ \ + \ \ \clustersize\beta^{t - 1}(1 - \beta)\eta\right]\\
& = \underbrace{\clustersize(1 - \beta)\tau}_{\nu} \ + \ \underbrace{\clustersize(1 - \beta)\kappa}_{\alpha} X_{t - 1} \ + \ \beta \lambda_{t - 1}
\end{align*}
for $t \geq 2$. Combined with the relationship $X_t = \clustersize * (\import_t + A_t) = \clustersize * N_t$ this is the form a CP-INGARCH(1, 1) model as introduced in \eqref{eq:N_CP_original}--\eqref{eq:lambda_CP_original}. Concerning the initialization of the process, the same argument as in Subsection \ref{subsubsec:derivation_poisson11} implies that we have to set $\lambda_0 = \clustersize\times \left\{\tau + (1 - \beta) \times \eta\right\}$.

The equations
\begin{align*}
    \tau & = (\nu/\clustersize)/(1 - \beta)\\
    \kappa & = (\alpha/\clustersize)/(1 - \beta)\\
    \eta & = (\lambda_0/\theta - \tau)/(1 - \beta)
\end{align*}
provided in the manuscript result from some simple shifting around of terms.

\section{Details on Section \ref{sec:extension_inar}: Extending the INAR class}

Here we collect some additional statements on the INARMA model class. Proofs both for statements from the main manuscript and this appendix are listed in Subsection \ref{subsec:proofs_inarma}.

\subsection{Details on the INARMA(1, 1) model with general import distribution}

\begin{remark}
\label{remark:inarma11}
For an INARMA(1, 1) process, representation \eqref{eq:galton_watson} of $\{E_t\}$ reduces to an INAR(1) model as
\begin{equation}
E_t = \{\beta + (1 - \beta)\kappa\} \circ E_{t - 1} + \kappa \circ \import_{t - 1}.\label{eq:E_INAR}
\end{equation}
This follows directly from Lemma \ref{lemma:galton_watson}.

The process $\{X_t\}$ also has a purely autoregressive formulation
\begin{equation}
X_t = \import_t \ + \ \sum_{i  = 1}^t \alpha_i \circ X_{t - i} \ + \ (\alpha_t/\kappa) \circ E_0\label{eq:purely_autoregressive}
\end{equation}
where $\alpha_i = \kappa\{1 - \beta\}\beta^{i - 1}$. Here, in slight abuse of notation we set
\begin{align*}
(\alpha_1 \circ X_t, \alpha_2 \circ X_t, \dots) \ \mid \ X_t & \sim \textnormal{Mult}(X_{t}; \alpha_1, \alpha_2, \dots),\\
(\alpha_1/\kappa \circ E_0, \alpha_2/\kappa \circ E_0, \dots) \ \mid  \ E_0 & \sim \textnormal{Mult}(E_{0}; \alpha_1/\kappa, \alpha_2/\kappa, \dots).
\end{align*}
This amounts to an INAR($p$) model as defined by \cite{Alzaid1990}, but with an infinite number of lags, geometrically decaying autoregressive parameters and a somewhat peculiar initialization (note that the term $(\alpha_t/\kappa) \circ E_0$ becomes negligible for large $t$). The proof is straightforward and follows the outline from Section 4 in \cite{Bracher2019}.
For the INGARCH(1, 1) model, the corresponding expression has been mentioned by  \cite{Lu2021}.
\end{remark}

\subsection{Details on the Poisson INARMA(1, 1) model}
\begin{remark}
Assume that $\{X_t\}$ is a Poisson INARMA(1, 1) process with parameters $\tau, \kappa, \beta$. The binomially thinned process $\{\tilde{X}_t\}$ where
$$
\tilde{X}_t = \pi \circ X_t
$$
independently for each $t$ is then equivalent in distribution to another INARMA(1, 1) process $\{Y_t\}$ with parameters
\begin{align*}
    \tau_Y & = \frac{\tau\pi}{1 - (1 - \pi)\kappa}\\
    \kappa_Y & = \frac{\pi\kappa}{1 - (1 - \pi)\kappa}\\
    \beta_Y & = \beta + (1 - \pi)\kappa(1 - \beta).
\end{align*}
This has been shown in \cite{Bracher2019a} (using representation \eqref{eq:purely_autoregressive}). The above expressions for the parameterization used here follow from simple algebra.

In a similar manner it can be shown that each Poisson INARMA(1, 1) process $\{X_t\}$ with parameters $\tau, \kappa, \beta$ is equivalent to a binomially thinned INAR(1) model $\{\tilde{Y_t}\}$. The latter is given by
\begin{equation}
Y_t = \xi \circ Y_{t - 1} + \import^*_t \label{eq:equiv_INAR1_Y}
\end{equation}
and
\begin{equation}
\tilde{Y}_t = \{(1 - \beta)\kappa/\xi\} \circ Y_t \label{eq:equiv_INAR1_Ytilde}  
\end{equation}
independently for each $t$, with $\xi = \beta + (1 - \beta)\kappa$ as in equation \eqref{eq:define_xi}. The import distribution is given by
$$
\import^*_t \stackrel{\textnormal{i.i.d.}}{\sim} \textnormal{Pois}(\tau\xi/\kappa),
$$
see again \cite{Bracher2019a} for the derivation. This representation as an imperfectly observed INAR(1) process is an interesting parallel to the Gaussian ARMA(1, 1) process, which is equivalent to as a mismeasured AR(1) process \cite{Staudenmayer2005}. We note that binomially thinned INAR(1) processes have first been studied by \cite{Fernandez-Fontelo2016}.

Finally we note that the representation \eqref{eq:equiv_INAR1_Y}--\eqref{eq:equiv_INAR1_Ytilde} implies that the time reversibility of the Poisson INAR(1) translates to the Poisson INARMA(1, 1). This is in line with the time-reversibility of the INAR($p$) model by \cite{Alzaid1990} as established by \cite{Schweer2015}.
\label{remark:underreporting}
\end{remark}

\subsection{Details on Hermite INARMA(1, 1) processes}
\label{subsec:hermite_inarma}

\begin{remark}
\label{remark:hermite}
If the import distribution of an INARMA(1, 1) model is given by
$$
\import_t \stackrel{\textnormal{i.i.d.}}{\sim} \textnormal{Herm}(\tau, \disp),
$$
and the initial value $E_0$ is assigned the distribution
$$
E_0 \sim \textnormal{Herm}\left(\frac{\kappa\tau}{1 - \xi}, \frac{\kappa\disp}{1 + \xi}\right),
$$
then the process $\{X_t\}$ is strictly stationary with $E_t$ marginally following the same distribution as $E_0$ and
$$
X_t \sim \textnormal{Herm}\left(\frac{\tau}{1 - \kappa},
(1 - \kappa)\times\disp + \kappa\times\frac{(1 - \beta)\kappa\disp}{1 + \beta + (1 - \beta)\kappa}
  \right).
$$
\end{remark}

\subsection{Details on the Poisson INARMA($p, q$) model}

\begin{lemma}
\label{lemma:embedded_inar}
Under model \eqref{eq:X_inarma}--\eqref{eq:C_t_inarma} with $\import_t \stackrel{\textnormal{i.i.d.}}{\sim} \textnormal{Pois}\left(\tau\right)$, the process $\{E_t\}$ is a Poisson INAR$(\max\{p, q\})$ process as defined in equations \eqref{eq:inar_p}--\eqref{eq:multinomial_thinning}. Setting $p = q$ without loss of generality, it is given by
\begin{equation}
E_t = \import^*_t + \sum_{i = 1}^p \alpha_i^* \circ E_{t - i}
\label{eq:inarma_pq_poisson}
\end{equation}
where
\begin{align*}
(\alpha_1^* \circ E_t, \dots, \alpha_p^* \circ E_t) \ \mid \ E_t & \sim \textnormal{Mult}(E_t, \alpha^*_1, \dots, \alpha^*_p),\\
\alpha^*_i & = \beta_i + \left(1 - \sum_{j = 1}^q \beta_j\right) \times \kappa_i\\  
\import^*_t & \stackrel{\textnormal{i.i.d.}}{\sim} \textnormal{Pois}\left(\tau \times \sum_{i = 1}^p \kappa_i\right).
\end{align*}
\end{lemma}

\begin{remark}
\label{remark:details_acf_arma}
    The Gaussian ARMA($p$, $p - 1$) process is defined as
    $$
	Y_t = \sum_{i = 1}^p \alpha_i Y_{t - i} + \sum_{i = 1}^{p - 1} \delta_i \import_{t - i} + \import_t, \quad \import_t \stackrel{\textnormal{i.i.d.}}{\sim}N(0, \sigma^2_{\import})
	$$
    and its ACF can be written as
    \begin{align}\label{eq:arma_acf}
        \gamma_\textnormal{Gauss}(d) & = \sum_{i = 1}^p \alpha_i \gamma_\textnormal{Gauss}(d - i) + \sum_{i = d}^{p - 1} \textnormal{Cov}(\delta_i \import_{t + d - i}, Y_t) + \mathbb{1}_{\{d = 0\}}\sigma^2_{\import}.
    \end{align}
    
    \noindent Alzaid and Al-Osh in \cite{Alzaid1990} defined the INAR($p$) model as
    \begin{align*}
        X_t = \sum_{i = 1}^p \phi_i \circ X_{t - i} + \import_t,\\
        \left[\phi_1 \circ X_t, \dots, \phi_p \circ X_t\right] \ \mid \ X_t & \sim \textnormal{Mult}\left(X_t; \phi_1, \dots, \phi_p\right),
    \end{align*}
    where $\import_t$ are iid non-negative integer valued random variables with mean $\mu_{\varepsilon}$ and variance $\sigma^2_{\varepsilon}$. They derived a formula for its covariance structure, given by
    \begin{equation}\label{eq:inar_acf}
        \gamma_\textnormal{INAR}(d) = \sum_{i = 1}^p \phi_i \gamma_\textnormal{INAR}(d - i) + \sum_{i = d + 1}^p \mu_\textnormal{INAR}(d - i, \phi_i) + \mathbb{1}_{\{d = 0\}} \sigma^2_{\import},
    \end{equation}
    where the terms $\mu_\textnormal{INAR}(d - i, \phi_i)$ are linear functions of the mean value of the process $\mu_X$. For more details see \cite{Alzaid1990}. We can spot that \eqref{eq:inar_acf} resembles \eqref{eq:arma_acf} in the sense that both consist of a sum of $p$ AR-terms ($\phi_i \gamma_\textnormal{INAR}(d - i)$ for INAR and $a_j \gamma_\textnormal{Gauss}(d - i)$ for ARMA) and a sum of $p - d - 1$ MA-terms ($\mu_\textnormal{INAR}(d - i, \phi_i)$ for INAR and $\textnormal{Cov}(\delta_j \import_{t + d - i}, Y_t)$ for ARMA). Note however that, unlike sometimes claimed in the literature, the autocovariance functions of both models are not actually equivalent and only share a certain resemblance.

    Following the same arguments as in \cite{Alzaid1990}, we can write the ACF of an INARMA($p$, $q$) model \eqref{eq:X_inarma}--\eqref{eq:C_t_inarma} in a similar fashion as
    \begin{equation}\label{eq:inarma_acf}
        \gamma_X(d) = \sum_{i = 1}^{\max\{p, q\}} \alpha_i^*\gamma_X(d - i) + \sum_{i = d}^{\max\{p, q\}} \mu(d - i, i) + \mathbb{1}_{\{d = 0\}}\sigma^2_{\tau}.
    \end{equation}	
    This corresponds to ${\max\{p, q\}}$ AR-terms $\alpha_i^*\gamma_X(d - i)$ and ${\max\{p, q\}} - d$ MA-terms $\mu(d - i, i)$, resembling the autocovariance structure of the Gaussian ARMA($\max\{p, q\}$, $\max\{p, q\}$). The term $\mu(-l, i)$ from equation \eqref{eq:inarma_acf} is defined as
        $$
        \mu(-l, i) := \begin{cases}
            \sum\limits_{j = 1}^{l - 1} \alpha_j^*\mu(j - l, i) + \left(1 - \sum\limits_{k = 1}^{\max\{p, q\}} \hspace{-5pt}\beta_k\right)\mu_{i, l} \quad &\text{for } 1 \leq i \leq \max\{p, q\} \text{ and } l > 0\\
            \mu(0, i) = - \beta_d \sigma^2_{\tau} - \beta_d\left(1 - \sum\limits_{k = 1}^{\max\{p, q\}} \hspace{-5pt} \beta_k\right) \mu_E \quad &\text{for } 1 \leq i \leq \max\{p, q\} \text{ and } l = 0\\
			\sum\limits_{k = 1}^{\max\{p, q\}} \hspace{-5pt} \beta_k \left(1 - \sum\limits_{k = 1}^{\max\{p, q\}} \hspace{-5pt} \beta_k\right) \mu_E \quad &\text{for } i = 0 \text{ and } l = 0
        \end{cases}
    $$
    and we use the notation
    \begin{align*}
        \alpha_i^* &:= \beta + \left(1 - \sum_{k = 1}^{\max\{p, q\}}\hspace{-7pt} \beta_k \right) \kappa_i,\\
        \mu_{i, l} &:= \Bigl(1 - \sum_{k = 1}^{\max\{p, q\}}\hspace{-7pt} \beta_k\Bigr)\beta_i\mathbb{1}_{\{i = l\}}\mu_E  + \\
        &\phantom{:=} + \kappa_l\biggl\{\Bigl(1 - \sum_{k = 1}^{\max\{p, q\}} \hspace{-7pt}\beta_k\Bigr)(\mathbb{1}_{\{i = l\}} - \kappa_i)\mu_X - \beta_i\sigma^2_{\tau} - \Bigl(1 - \sum_{k = 1}^{\max\{p, q\}} \hspace{-7pt}\beta_k\Bigr)\beta_i\mu_E  \biggr\}.
    \end{align*}
    If $p > q$, we set $\beta_{q + 1} = \beta_{q + 2} = \dotsc = \beta_{p} = 0$ and vice versa for the parameters $\kappa_l$.
    
    Since in the INARMA($p$, $q$) model we need to work with the additional process $\{E_t\}$, the derivations become somewhat bulky. For this reason we omit them here.
\end{remark}

\subsection{Moment-based estimation}
\label{suppl_moment_based_estimation}

\subsubsection{Procedure}

A computationally cheaper alternative to maximum-likelihood estimation is moment-based estimation, as discussed for the Poisson INAR(1) and INARCH(1) processes in \cite{Weiss2016}. For the Poisson INARMA(1, 1) process, solving the system of equations \eqref{eq:moments_poisson_11}  for $\tau, \beta$, and $\kappa$ yields the following moment estimators:
\begin{equation}
\hat{\tau} = \frac{\hat{\mu}_X\{\hat{\rho}_X(1) - \hat{\rho}_X(2)\}}{\hat{\rho}_X(1)^2 + \hat{\rho}_X(1) - \hat{\rho}_X(2)}, \ \ \ \hat{\beta} = \frac{\hat{\rho}_X(2)}{\hat{\rho}_X(1)} - \hat{\rho}_X(1), \ \ \ \hat{\kappa} = \frac{\hat{\rho}_X(1)^2}{\hat{\rho}_X(1)^2 + \hat{\rho}_X(1) - \hat{\rho}_X(2)}.\label{eq:moment_estimators}
\end{equation}

\begin{proposition}
If $\{X_t, t \in \mathbb{N}\}$ is a Poisson INARMA(1, 1) process with $0 < \kappa < 1$, the moment estimators $\hat{\tau},\hat{\beta}, \hat{\kappa}$ from equation \eqref{eq:moment_estimators} are consistent and asymptotically normal.\label{proposition:moments}
\end{proposition}

\noindent For innovation distributions other than the Poisson, solving equations \eqref{eq:mu_X}--\eqref{eq:rho_X} for the model parameters is somewhat tedious and boils down to solving a cubic equation.

\begin{lemma}\label{lemma:kappa_cubic}
    If $\{X_t, t \in \mathbb{N}\}$ is an INARMA(1, 1) process with parameters $\tau, \sigma^2_\tau, \kappa, \beta$ and limiting-stationary second order properties as given in Lemma \ref{lemma:moments_inarma_11}, it holds that
    \begin{equation}
        a\kappa^3 + b\kappa^2 + c\kappa + d = 0,
    \end{equation}
    where
    \begin{align}
        a & = (1 - \xi)(\mu_X + \sigma^2_X) + 2\gamma_X(1)\nonumber\\
        b & = - (1 - \xi)\{(2 + \xi)\sigma^2_X + \xi \mu_X\} -  2\times(2 + \xi)\times\gamma_X(1)\label{eq:abcd}\\
        c & = (1 - \xi^2)\sigma^2_X  + 3\times(1 + \xi)\times\gamma_X(1)\nonumber\\
        d & = -(1 + \xi)\gamma_X(1).\nonumber
    \end{align}  
\end{lemma}

The moment estimator for $\kappa$ is obtained by replacing $\mu_X, \sigma^2_X, \gamma_X(1), \xi$ by their empirical counterparts $\hat{\mu}_X,$ $\hat{\sigma}^2_X,$ $\hat{\gamma}_X(1),$ $\hat{\gamma}_X(2)/\hat{\gamma}_X(1)$ in equation \eqref{eq:abcd} and finding a solution $\hat{\kappa} \in [0, 1]$ numerically. Given $\hat{\kappa}$, the other parameter estimates can be computed via the relationships
    \begin{align*}
        \tau = \mu_X\times(1 - \kappa);\ \ \
        \beta = \frac{\xi - \kappa}{1 - \kappa};\ \ \
        \sigma^2_\tau & = (1 - \kappa) \times \frac{\sigma^2_X - \frac{\kappa\times(1 + \beta)}{1 + \xi} \times \mu_X}{1 - \frac{\kappa\times(1 + \beta)}{1 + \xi}}.
    \end{align*}
The dispersion parameter of the innovation distribution (e.g., $\disp$ in the Hermite or negative binomial distributions) can be obtained from $\hat{\tau}$ and $\hat{\sigma}^2_{\tau}$ via the respective mean-variance relationship.

There is no simple closed form for the solution of the cubic equation \eqref{eq:abcd}, and in fact it can have multiple real-valued solutions. However, in all cases we examined, there was only one solution for $\kappa$ from the unit interval, and thus only one solution which will lead to positive estimates for all model parameters. While we cannot provide a formal proof for this, based on extensive numerical studies, we conjecture that under the following regularity conditions there is always exactly one solution for $\hat{\kappa} \in [0, 1]$:
\begin{itemize}
    \item[(i)] $\hat{\mu}_X > 0$.
    \item[(ii)] $\hat{\sigma}^2_X > 0$.
    \item[(iii)] $\hat{\gamma}_X(1) > 0$.
    \item[(iv)] $\hat{\gamma}_X(2) \geq \hat{\gamma}_X(1)^2$ (or equivalently $\hat{\xi} \geq \hat{\gamma}(1)$).
\end{itemize}
As the negative binomial and Hermite immigration distributions cannot reflect underdispersion, when using these parameters it is further required that
\begin{itemize}
    \item[(v)] $\hat{\sigma}^2_X \geq \hat{\mu}_X$.
\end{itemize}
Whenever these conditions are not fulfilled this is an indication that the specified INARMA model may not be an appropriate choice. For the purpose of our simulation studies (see next section for results) we set $\hat{\sigma}^2_X = \hat{\mu}_X$ and $\hat{\xi} = \hat{\gamma}_X(1)$ if conditions (iv) and (v) were violated for a simulated time series. As estimation moreover becomes numerically instable if $\hat{\gamma}(1)$ and $\hat{\xi}$ are too close to 1, in practice we threshold them at 0.95.

It is known that moment-based estimators are subject to small-samples biases even for simpler INAR model \cite{Weiss2016}. Deriving these biases analytically seems too technically involved for INARMA models, but the simulation results in the next section illustrate the general problem.

\subsubsection{Simulation results for moment-based estimators}

\begin{table}[h!]
\spacingset{1.3} 
\scriptsize
\caption{Simulation results for moment-based estimators in the Poisson, Hermite and negative binomial settings, scenarios 1--3 with $T \in \{250, 500, 1000\}$ and 1000 runs.}
\label{tab:sim_moments}
\center
\medskip
\begin{tabular}{p{0.3cm} @{\hskip 0.7cm} p{0.45cm} p{0.45cm} p{0.45cm} @{\hskip 0.7cm} p{0.45cm} p{0.45cm} p{0.45cm} @{\hskip 0.7cm} p{0.45cm} p{0.45cm} p{0.45cm} @{\hskip 0.7cm} p{0.45cm} p{0.45cm} p{0.45cm}}
\hline
\multicolumn{13}{c}{Poisson}\\
\hline
$T$ & \multicolumn{3}{c}{$\tau$} & \multicolumn{3}{c}{$\disp$} & \multicolumn{3}{c}{$\beta$} & \multicolumn{3}{c}{$\kappa$}\\
\hline 
& true & mean & se & true & mean & se & true & mean & se & true & mean & se \\
\hline 
\input{table/sim_pois_sc1_moments.tex}

\input{table/sim_pois_sc2_moments.tex}

\input{table/sim_pois_sc3_moments.tex}\\
\hline\\
\multicolumn{13}{c}{Hermite}\\
\hline 
$T$ & \multicolumn{3}{c}{$\tau$} & \multicolumn{3}{c}{$\disp$} & \multicolumn{3}{c}{$\beta$} & \multicolumn{3}{c}{$\kappa$}\\
\hline 
& true & mean & se & true & mean & se & true & mean & se & true & mean & se \\
\input{table/sim_herm_sc1_moments.tex}

\input{table/sim_herm_sc2_moments.tex}

\input{table/sim_herm_sc3_moments.tex}\\
\hline\\
\multicolumn{13}{c}{Negative binomial}\\
\hline
$T$ & \multicolumn{3}{c}{$\tau$} & \multicolumn{3}{c}{$\disp$} & \multicolumn{3}{c}{$\beta$} & \multicolumn{3}{c}{$\kappa$}\\
\hline 
& true & mean & se & true & mean & se & true & mean & se & true & mean & se \\
\hline 
\input{table/sim_negbin_sc1_moments.tex}

\input{table/sim_negbin_sc2_moments.tex}

\input{table/sim_negbin_sc3_moments.tex}

\end{tabular}
\end{table}

\subsection{Proofs}
\label{subsec:proofs_inarma}

\subsubsection{Limiting-stationary moments of the INARMA(1, 1) with generic import distributions, lemma \ref{lemma:moments_inarma_11}}\label{subsec:moments_inarma_11}

While these properties could in principle be obtained from the more general results in Lemma \ref{lemma:general_moments}, it seems more instructive and not much more difficult to derive them from scratch. We start by noting a few well-known properties of the binomial thinning operator.

\begin{lemma}
\label{lemma:properties_thinning}
For the binomial thinning operator $\circ$ and an arbitrary integer-valued random variable $A$, the following hold:
\begin{align}
\mathbb{E}(\alpha \circ A) & = \alpha \times \mathbb{E}(A)\\
\textnormal{Var}(\alpha \circ A) & = \alpha(1 - \alpha)\mathbb{E}(A) + \alpha^2 \textnormal{Var}(A)\\
\textnormal{Cov}(\alpha \circ A, A) & = \alpha \text{Var}(A) \label{eq:cov_thinned}
\end{align}
Moreover, it the thinning $\beta \circ B$ is performed independently of $A$ then
$$
\textnormal{Cov}(A, \beta \circ B) = \beta\textnormal{Cov}(A, B).
$$
\end{lemma}

We then turn to the process $\{E_t\}$, which as stated in Remark \ref{remark:inarma11} can be represented as an INAR(1) process with
$$
E_t = \underbrace{\xi}_{\beta + (1 - \beta)\kappa} \ \circ E_{t - 1} + \underbrace{\import^*_t}_{\kappa \circ \import_{t - 1}}.
$$
Using Lemma \ref{lemma:properties_thinning} we obtain
$$
\mathbb{E}(\import^*_t) = \kappa \tau, \ \ \ \text{Var}(\import^*_t) = \kappa(1 - \kappa)\tau + \kappa^2 \sigma^2_\tau.
$$
Well-known properties of the INAR(1) model \cite{Weiss2018} then imply that
\begin{align*}
\mu_E & = \frac{\kappa\tau}{1 - \xi}\\
\sigma^2_E & = \frac{\text{Var}(\import^*_t) + \xi \times \mathbb{E}(\import^*_t)}{1 - \xi^2} = \frac{\kappa^2\sigma^2_\tau + \kappa(1 - \kappa)\tau + \xi\kappa\tau}{1 - \xi^2}\\
& = \frac{\kappa^2\sigma^2_\tau + \kappa(1 - \kappa + \xi)\tau}{1 - \xi^2} = \frac{\kappa^2\sigma^2_\tau + \kappa(1 - \beta\kappa + \beta)\tau}{1 - \xi^2}\\
\rho_E(d) & = \xi^d.
\end{align*}
We can then turn back to the moments of $\{X_t\}$, where after some simple algebra we obtain
$$
\mu_X = (1 - \beta)\mu_E + \tau = ... = \frac{\tau}{1 - \kappa}.
$$
For the variance we use Lemma \ref{lemma:properties_thinning} to note that
\begin{align*}
\sigma^2_X & = (1 - \beta)^2 \sigma^2_E + \beta(1 - \beta)\mu_E + \sigma^2_\tau.
\end{align*}
After some simple re-ordering of terms this leads to
\begin{align}
\sigma^2_X & = \frac{(1 - \beta^2 - 2\beta\kappa + 2\beta^2\kappa)\times \sigma^2_\tau + (1 - \beta^2)\kappa\tau}{1 - \xi^2}\label{eq:sigma2_X_intermediate}
\end{align}
Splitting this term at the summation in the numerator, we get
\begin{align*}
\frac{(1 - \beta^2)\kappa\tau}{1 - \xi^2} = \frac{(1 - \beta)(1 + \beta)\kappa\tau}{(1 - \xi)(1 + \xi)} = \frac{(1 - \beta)(1 + \beta)\kappa\tau}{(1 - \beta)(1 - \kappa)(1 + \xi)} = \frac{(1 + \beta)\kappa}{1 + \xi} \times \frac{\tau}{1 - \kappa}
\end{align*}
and
\begin{align*}
\frac{(1 - \beta^2 - 2\beta\kappa + 2\beta^2\kappa)\times \sigma^2_\tau}{1 - \xi^2} & = \frac{(1 - \beta^2 - 2\beta\kappa + 2\beta^2\kappa)(1 - \kappa)}{1 - \xi^2} \times \frac{\sigma^2_\tau}{1 - \kappa}\\
& = \left(1 - \frac{1 - \xi^2 - (1 - \beta^2 - 2\beta\kappa + 2\beta^2\kappa)(1 - \kappa)}{(1 - \xi)(1 + \xi)}\right) \times \frac{\sigma^2_\tau}{1 - \kappa}\\
& = \ \ \ \ \dots \ \ \ \text{[multiplying out and re-ordering terms]}\\
& = \left(1 - \frac{\kappa - \kappa^2 - \kappa\beta^2 + \beta^2\kappa^2}{1 - \xi^2}\right) \times \frac{\sigma^2_\tau}{1 - \kappa}\\
& = \left(1 - \frac{(1 - \beta^2)\kappa(1 - \kappa)}{(1 - \beta)(1 - \kappa)(1 + \xi)}\right) \times \frac{\sigma^2_\tau}{1 - \kappa}\\
& = \left(1 - \frac{(1 + \beta)\kappa}{1 + \xi}\right) \times \frac{\sigma^2_\tau}{1 - \kappa}.
\end{align*}
Plugging these back into \eqref{eq:sigma2_X_intermediate} we then get
\begin{align}
\sigma^2_X & = \frac{(1 + \beta)\kappa}{1 + \xi} \times \frac{\tau}{1 - \kappa} \ + \ \left(1 - \frac{(1 + \beta)\kappa}{1 + \xi}\right) \times \frac{\sigma^2_\tau}{1 - \kappa}.
\end{align}
To obtain the autocovariance function we use again the auxiliary notation
$$
A_t = E_t - L_t = \beta \circ E_t
$$
as in \ref{subsec:ginarma_properties}, \eqref{eq:A_general_11}. We then consider
\begin{align}
\text{Cov}(X_t, E_{t + 1}) & = \text{Cov}(X_t, E_t - A_t + \kappa \circ X_t) \nonumber \\
& = \text{Cov}(X_t, E_t) - \text{Cov}(X_t, A_t) + \text{Cov}(X_t, \kappa\circ X_t)\nonumber\\
& = \text{Cov}(\{1 - \beta\} \circ E_t + \import_t, E_t) - \text{Cov}(A_t + \import_t, A_t) + \text{Cov}(X_t, \kappa\circ X_t) \nonumber\\
& \stackrel{\eqref{eq:cov_thinned}}{=} \text{Cov}(\{1 - \beta\} \circ E_t, E_t) - \text{Var}(A_t) + \text{Cov}(X_t, \kappa\circ X_t) \nonumber \\
& \stackrel{\eqref{eq:cov_thinned}}{=} (1 - \beta) \sigma^2_E - \beta(1 - \beta)\mu_E - (1 - \beta)^2\sigma^2_E + \kappa\sigma^2_X \nonumber \\
& = \beta(1 - \beta)(\sigma^2_E - \mu_E) + \kappa\sigma^2_X.\nonumber
\end{align}
We rewrite this to a form we will use later.
\begin{align}
\text{Cov}(X_t, E_{t + 1}) & = \beta(1 - \beta)(\sigma^2_E - \mu_E) + \kappa\sigma^2_X \nonumber \\
& = \kappa\sigma^2_X + \beta(1 - \beta)\left(\frac{\kappa^2\sigma^2_{\tau} + \kappa\tau(1 - \beta\kappa + \beta)}{1 - \xi^2} - \frac{\kappa\tau}{1 - \xi} \right)\nonumber \\
& = \kappa\sigma^2_X + \beta(1 - \beta)\left(\frac{\kappa^2\sigma^2_{\tau} + \kappa\tau(1 - \beta\kappa + \beta) - (1 + \xi)\kappa\tau}{1 - \xi^2} \right)\nonumber \\
& = \kappa\sigma^2_X + \beta(1 - \beta)\left(\frac{\kappa^2\sigma^2_{\tau} + \kappa\tau\overbrace{\{1 - \beta\kappa + \beta - 1 - \beta - (1 - \beta)\kappa\}}^{\text{this reduces to }-\kappa}}{1 - \xi^2} \right)\nonumber \\
& = \kappa\sigma^2_X + \beta(1 - \beta)\kappa\left(\frac{\kappa\sigma^2_{\tau} - \kappa\tau}{1 - \xi^2} \times \frac{\sigma^2_X}{\sigma^2_X}\right)\nonumber \\
& \stackrel{\eqref{eq:sigma2_X_intermediate}}{=} \kappa\sigma^2_X + \beta(1 - \beta)\kappa\left(\frac{\kappa\sigma^2_{\tau} - \kappa\tau}{(1 - \beta^2 - 2\beta\kappa + 2\beta^2\kappa)\times \sigma^2_\tau + (1 - \beta^2)\kappa\tau}\times\sigma^2_X \right)\nonumber\\
& = \kappa\sigma^2_X + \beta\kappa\sigma^2_X\left(\frac{\kappa\sigma^2_{\tau} - \kappa\tau}{(1 + \beta - 2\beta\kappa)\times \sigma^2_\tau + (1 + \beta)\kappa\tau}\right)\nonumber\\
& = \kappa\sigma^2_X \times \left(1 + \frac{\kappa\beta(\sigma^2_{\tau} - \tau)}{(1 + \beta - \beta\kappa - \kappa + \kappa - \beta\kappa)\times \sigma^2_\tau + (1 + \beta)\kappa\tau}\right) \nonumber \\
& = \kappa\sigma^2_X \times \left(1 + \frac{\kappa\beta(\sigma^2_{\tau} - \tau)}{(1 + \beta) \{ (1 - \kappa) \sigma^2_{\tau} + \kappa\tau \} + (1 - \beta)\kappa\sigma^2_{\tau}}\right).\label{eq:acov1_intermediate}
\end{align}

\noindent Secondly we prove by induction that
\begin{align}
\text{Cov}(X_t, E_{t + d}) & = \xi^{d - 1} \text{Cov}(X_t, E_{t + 1}), \ \text{for all } d \in \mathbb{N}.
\label{eq:acov_induction}
\end{align}
For $d = 1$ the equation holds trivially. Suppose that it holds for all $d$ up to some $k$. Then for $k + 1$ we have
\begin{align*}
\text{Cov}(X_t, E_{t + k + 1}) & = \text{Cov}(X_t, \beta \circ E_{t + k} + \kappa \circ X_{t + k})\\
& = \beta \times \text{Cov}(X_t, E_{t + k}) + \kappa \times \text{Cov}(X_t, (1 - \beta) \circ E_{t + k} + \import_{t + k}) \\
& = \beta \times\text{Cov}(X_t, E_{t + k}) + \kappa (1 - \beta)\times \text{Cov}(X_t, E_{t + k}) + (1 - \beta) \underbrace{\times \text{Cov}(X_t, \import_{t + d})}_{=0}\\
& = \underbrace{\{\beta + \kappa(1 - \beta)\}}_{\xi} \times \text{Cov}(X_t, E_{t + k})\\
& =\xi^k \text{Cov}(X_t, E_{t + 1}),
\end{align*}

\noindent which concludes the induction step.

Finally we calculate the autocovariance for the process $X_t$ as
\begin{align}
\text{Cov}(X_t, X_{t + d}) & = \text{Cov}(X_t, (1 - \beta) \circ E_{t + d} + \import_{t + d}) \nonumber\\
& = (1 - \beta) \times \text{Cov}(X_t, E_{t + d}) + \times \underbrace{\text{Cov}(X_t, \import_{t + d})}_{=0} \nonumber\\
& = (1 - \beta) \xi^{d - 1} \times \text{Cov}(X_t, E_{t + 1})
\label{eq:acovX_intermediate}
\end{align}
Plugging equations \eqref{eq:acov_induction} and \eqref{eq:acov1_intermediate} into \eqref{eq:acovX_intermediate} gives
\begin{equation}
\gamma_X(d) = (1 - \beta)\xi^{d - 1} \kappa \sigma^2_X\times\left(1 + \frac{\kappa\beta(\sigma^2_\tau - \tau)}{(1 + \beta) \{(1 - \kappa)\sigma^2_\tau + \kappa\tau\} + (1 - \beta)\kappa\sigma^2_\tau}\right).
\end{equation}

\subsubsection{Stationary moments of Poisson INARMA(1, 1), lemma \ref{proposition:bivariate_poisson}} 

The stationary mean and autocorrelation function follow directly from the more general results in lemma \ref{lemma:moments_inarma_11}. The statement on the bivariate Poisson distribution is proven for the more general case of the INARMA($p, q$) as an intermediate step in the proof of Corollary \ref{corollary:moments_inarma_pq}, see Section \ref{suppl:moments_inarma_p_q}.

We here only briefly justify that the provided Poisson distribution is the strictly stationary rather than just limiting-stationary distribution of $\{X_t\}$. To this end remember that $\{E_t\}$ is a Poisson INAR(1) process (Remark \ref{remark:inarma11})
\begin{equation}
E_t = \{\beta + (1 - \beta)\kappa\} \circ E_{t - 1} + \underbrace{\kappa \circ \import_{t - 1}}_{\varepsilon^*_t},
\end{equation}
where the imports $\varepsilon^*_t = \kappa \circ \import_{t - 1}$ independently follow a Poisson distribution with rate $\kappa\tau$. The stationary distribution of $\{E_t\}$ is thus Poisson with rate $\kappa\tau/\{1 - \beta - (1 - \beta)\kappa\}$. If we initialize $E_0$ with this distribution and set $X_0 = (1 - \beta) \circ E_0 + \varepsilon_0$, the process $\{E_t\}$ will be strictly stationary. As the joint process $\{E_t, X_t\}$, too, only depends on the initialization of $E_0$, it is likewise initialized in its stationary regime and thus strictly stationary.

\subsubsection{Limiting-stationary moments of Poisson INARMA($p, q$) model, corollary \ref{corollary:moments_inarma_pq}}
\label{suppl:moments_inarma_p_q}

We consider model \eqref{eq:X_inarma}--\eqref{eq:C_t_inarma} with Poisson imports, i.e., $\import_t \stackrel{\textnormal{i.i.d.}}{\sim} \text{Pois}(\tau)$. From lemma \ref{lemma:embedded_inar} we know that $E_t$ is a Poisson INAR[$\max(p, q)$] model, meaning that $\{E_t\}$ marginally follows a Poisson distribution. Due to the closedness to binomial thinning and summation of the Poisson distribution it is clear that $\{X_t\}$ also follows a Poisson distribution. To obtain the mean of this distribution we first note that the stationary mean of the INAR($p$) process $\{E_t\}$ is
$$
\mu_E = \frac{\mathbb{E}(\import^*_t)}{1 - \sum_{i = 1}^{\max(p, q)} \alpha_i^*}.
$$
where as in lemma \ref{lemma:embedded_inar} we set
$$
\alpha^*_i = \beta_i + \left(1 - \sum_{j = 1}^q \beta_j\right) \times \kappa_i
$$
with $\kappa_i = 0$ for $i > p$ and $\beta_j = 0$ for $j > q$. The limiting-stationary mean of $\{X_t\}$ can then be obtained as
$$
\mu_X = \left(1 - \sum_{i = 1}^q \beta\right) \times \mu_E + \tau.
$$
After some algebra this simplifies to
$$
\mu_X = \frac{\tau}{1 - \sum_{i = 1}^p \kappa_i}.
$$
As $\{X_t\}$ is marginally Poisson, we conclude that $\sigma^2_X = \mu_X$.

To derive the autocorrelation function we use again an extended notation, setting
\begin{align*}
A_t & = \left(1 - \sum_{i = 1}^q \beta_j\right) \circ E_t\\
L_t & = \sum_{j = 1}^q \beta_j \circ E_t =  E_t - A_t.
\end{align*}
The Poisson splitting property implies that $A_t$ and $L_t$ are independently Poisson with rates $(1 - \sum_{i = 1}^q \beta_j) \times \mu_E$ and $(\sum_{i = 1}^q \beta_j) \times \mu_E$, respectively. We note that
$$
X_t = A_t + \import_t
$$
is likewise Poisson, as is $X_{t + d}$. As $L_t$ and $A_t$ are independent, so are $L_t$ and $X_t$.

Using the language from Section \ref{subsec:interpretation_epidemic}, we now consider $X_{t + d}, d > 0$ and distinguish these infectives by where we can track their ``chain of infection''. Each infective person in $X_{t + d}$ must fall into one of three mutually exclusive categories:
\begin{itemize}
    \item[(1)] direct or indirect offspring of cases imported after time $t$ (i.e., $\import_{t + 1}, \dots, \import_{t + d}$).
    \item[(2)] direct or indirect indirect offspring of an infective from time $t$ (i.e., $X_t$).
    \item[(3)] persons infected prior to time $t$, but not yet infectious at $t$ (i.e., in $L_t$).
\end{itemize}
In the remainder of this proof, we will write ``offspring'' for both direct and indirect offspring, i.e., there may be intermediate steps in the chain of infection. Whenever we only refer to direct offspring, this will be mentioned explicitly.

Each of the individuals from $L_t$, $X_t, \import_{t + 1}, \dots, \import_{t + d}$ can have zero or one offspring in $X_{t + d}$. Denoting the probability of having one offspring by $\rho(X_t, X_{t + d})$, $\rho(\import_{t + 1}, X_{t + d})$ etc., we thus have
$$
X_{t + d} \ \ = \ \ \rho(L_t, X_{t + d}) \circ L_t \ \ + \ \ \rho(X_t, X_{t + d}) \circ X_t \ \ + \ \ \sum_{i = 0}^{d} \rho(\import_{t + d - i}, X_{t + d}) \circ \import_{t + d - i}.
$$
As $L_t, X_t, \import_{t + 1}, \dots, \import_{t + d}$ are all independently Poisson, this is again a sum of independent Poisson random variables. We can thus easily see that (in agreement with Lemma \ref{proposition:bivariate_poisson}) $X_t$ and $X_{t + d}$ jointly follow a bivariate Poisson distribution
$$
(X_t, X_{t + d}) \sim \text{BPois}\bigl(\rho(X_t, X_{t + d})\mu_X, \{1 - \rho(X_t, X_{t + d})\}\mu_X, \{1 - \rho(X_t, X_{t + d})\}\mu_X\bigr).
$$
Notably, this implies that the correlation $\rho_X(d) = \text{Corr}(X_t, X_{t + d})$ is just $\rho(X_t, X_{t + d})$, justifying this notation.

It remains to compute the probability $\rho(X_t, X_{t + d})$ that an infective from $X_t$ has an offspring in $X_{t + d}$. The reasoning for this is as follows. For an infective from $X_t$ to have an offspring in $X_{t + d}$ that offspring must have been part of $E_{t + d}$ before. Given there is such an offspring in $E_t$, it will progress to $X_{t + d}$ with probability $1 - \sum_{j = 1}^q \beta_j$. Thus, using analogous notation as above we have
$$
\rho(X_t, X_{t + d}) = \left(1 - \sum_{j = 1}^q \beta_j \right) \times \rho(X_t, E_{t + d}).
$$
We are thus looking for the probability $\rho(X_t, E_{t + d})$ that an individual from $X_t$ has an offspring in $E_{t + d}$. We can find a recursion for this, but require some auxiliary quantities. We will split up all offspring of $X_t$ in $E_{t + d}$ by when their direct infector was infectious. We denote by $s_i$ the probability that an infective from $X_{t + d - i}$ has a \textit{direct} offspring in $E_{t + d}$ (i.e., there are no intermediate steps in the chain of infections). This enables us to write the total probability $\rho(X_t, X_{t + d})$ as the sum
$$
\rho(X_t, X_{t + d}) = \left(1 - \sum_{j = 1}^q \beta_j \right) \times \left(\sum_{i = 1}^d \underbrace{\rho(X_{t}, X_{t + d - i})}_{\substack{\text{infective from } X_t \\ \text{ has offspring in } X_{t + d - i}}} \times \underbrace{s_i}_{\substack{\text{infective from } X_{t + d - i} \\ \text{has direct offspring in } E_{t + d}}} \right).
$$
It remains to derive an expression for $s_i$. We will use a second auxiliary quantity $\pi_d$, which is the probability that an exposed individual present in $E_t$ moves on to $E_{t + d}$ without ever leaving the exposed pool (but potentially making intermediate steps inside the exposed pool). Using this, we can express $s_i$ as
$$
s_i = \sum_{k = 1}^{\min\{i, p\}} \hspace{-7pt}\kappa_k \pi_{i - k}.
$$
Finally, we need an expression for $\pi_k$, which can be obtained without further detours. We simply split up the individuals by the last time before $t + d$ when they made their appearance in the exposed pool (i.e., from which out of $E_{t + d - 1}, \dots E_{t + d - q}$ they jumped to $E_{t + d}$). This leads to the recursion
$$
\pi_k = \sum_{l = 1}^{\min\{k, q\}} \hspace{-7pt}\beta_l \pi_{k - l}.
$$
For the initialization, we need to set $\pi_0 = 1$ as the probability that an individual from $E_{t + d}$ appears in $E_{t + d}$ without leaving the exposed pool in the meantime is trivially equal to $1$.

\subsubsection{Embedded Poisson INAR($p$) process in the Poisson INARMA($p, q$) model, proof of corrolary \ref{corollary:embedded_inar_p_main} / lemma \ref{lemma:embedded_inar}}
\label{subsubsec:embedded_inar}

To show Lemma \ref{lemma:embedded_inar}, we will for simplicity and without loss of generality assume $p = q$. We thus consider a Poisson INARMA($p, p$) model, given by
\begin{align}
X_t & = \Bigl(1 - \sum_{j = 1}^p \beta_j \Bigr) \circ E_t + \import_t\\
E_t & = \sum_{j = 1}^p \beta_j \circ E_{t - j} \ \ + \ \ \sum_{i = 1}^p \kappa_i \circ X_{t - i}\\
\left[\beta_1 \circ E_t, \dots, \beta_p \circ E_t, \Bigl(1 - \sum_{j = 1}^p \beta_j \Bigr) \circ E_t\right] \ \mid \ E_t & \sim \textnormal{Mult}\left(E_t; \beta_1, \dots, \beta_p, 1 - \sum_{j = 1}^p \beta_j\right)\\
(\kappa_1 \circ X_t, \dots, \kappa_p \circ X_t) \ \mid \ X_t & \sim \text{Mult}\left(X_t; \kappa_1, \dots, \kappa_p\right).\label{eq:offsping_mul}
\end{align}
We will again use the epidemiological interpretation of the process as in Subsection \ref{subsec:interpretation_epidemic} (with multinomial offspring distributions). We now consider how the exposed pool gets ``renewed'', i.e., how individuals from $E_t$ can contribute to $E_{t + 1}, \dots, E_{t + p}$. A ``contribution'' by an exposed individual can either be that same exposed individual re-appearing in $E_{t + 1}, \dots, E_{t + p}$, or an offspring of said individual.

Each individual from $E_t$ can contribute ``directly'' to at most one out of $E_{t + 1}, \dots, E_{t + p}$. By a ``direct contribution'' to $E_{t + i}$ we mean that the individual has not contributed to the exposed pool between $t$ and $t + i$. Let us denote the number of individuals directly contributing from $E_t$ to $E_{t + i}$ by $C^*_{t, i}$. An individual from $E_t$ can become part of $C^*_{t, i}$ either by remaining exposed and moving directly to $E_{t + i}$ with probability $\beta_i$; or by advancing to the infectious pool $X_t$ with probability $\left(1 - \sum_{j = 1}^p \beta \right)$ and then generating another exposed individual at time $t + i$ with probability $\kappa_i$. The total probability of an individual from $E_t$ contributing directly to $E_{t + i}$ is thus
$$
\alpha^*_i = \beta_i + \left(1 - \sum_{j = 1}^p \beta \right)\kappa_i.
$$
Since the individuals behave independently of each another and each individual from $E_t$ can be part of at most one out of $C^*_{t, 1}, \dots, C^*_{t, p}$, we obtain a conditional multinomial distribution
$$
(C^*_{t, 1}, \dots, C^*_{t, p}) \ \mid E_t \ \sim \textnormal{Mult}(E_t, \alpha^*_1, \dots, \alpha^*_p).
$$

We can now use this to formulate a recursion for $E_t$. To this end we have to sum all direct contributions to $E_t$, which are those originating from $E_{t - 1}$ $(C^*_{t - 1, 1})$, from $E_{t - 2}$ $(C^*_{t - 2, 2})$ and so on up to $E_{t - p}$ $(C^*_{t - p, p})$. In addition to these, there are new exposed individuals caused by imported infections from previous times $\import_{t - j}, j = 1, \dotsc, p$. Each imported infective from $\import_{t - j}$ has a probability of $\kappa_j$ to directly contribute an exposed to $E_t$. We can thus write 
$$
E_t = \sum_{i = 1}^p C^*_{t - i, i} + \underbrace{\sum_{j = 1}^p \kappa_j \circ \import_{t - j}}_{\import^*_t},
$$
where $\kappa_1 \circ \import_{t}, \dots, \kappa_p \circ \import_{t}$ are likewise coupled by a multinomial distribution. Since $\import_t \stackrel{\textnormal{i.i.d.}}{\sim}\textnormal{Pois}(\tau)$ and due to the Poisson splitting property, $\import^*_t$ is just a sum of $p$ independent Poisson variables with means $\tau\kappa_i$, $i = 1, \dotsc, p$, respectively. This implies
$$
\import^*_t \sim \textnormal{Pois}\left(\tau \times \sum_{i = 1}^p \kappa_i\right).
$$
Moreover, due to the Poisson splitting property, it can be shown that $\{\import^*_t\}$ is a sequence of independent Poisson random variables. Note that for other immigration distributions this independence would not hold and we would obtain an INAR($p$) model with dependent imports. For $p > q$ the results remain valid with $\kappa_{q + 1} = \kappa_{q + 2} = \dotsc = \kappa_{p} = 0$. Analogously, the $\beta_j$ can be completed with zeros if $p < q$.

\subsubsection{Marginal distribution of Hermite INARMA(1, 1), remark \ref{remark:hermite}}
\label{suppl:derivation_hermite}

We consider the model INARMA(1, 1) model 
\begin{equation} \label{eq:inarma11}
    \begin{aligned}
        X_t & = (1 - \beta) \circ E_t + \import_t \\
        E_t & = \beta \circ E_{t - 1} + \kappa \circ X_{t - 1} \\
        \\
        \bigl[\beta \circ E_t, (1 - \beta) & \circ E_t\bigr] \ \mid \ E_t  \sim     \text{Mult}\left(E_t; \beta, 1 - \beta\right)
    \end{aligned}
\end{equation}
with the innovation distribution
$$
\import_t \stackrel{\textnormal{i.i.d.}}{\sim} \textnormal{Herm}(\tau, \disp),
$$
where we parameterize the Hermite distribution via its mean $\tau$ and a dispersion parameter $\disp$; see Section \ref{sec:preliminaries}. Recall that the mean and the variance are given by $\mathbb{E}(\import_t) = \tau$ and $\text{Var}(\import_t) = (1 + \disp)\tau$. Several useful properties of the Hermite distribution are listed in \cite{Fernandez-Fontelo2017}, see also references therein. Note that we adapted them to our notation, which requires merely some shifting around of terms.

\begin{lemma}
    The sum $A = B + C$ of two random variables $B \sim \text{Herm}(\mu_B, \disp_B)$ and $C \sim \text{Herm}(\mu_C, \disp_C)$ follows again a Hermite distribution with parameters $\mu_A = \mu_B + \mu_C$ and $\disp_A = (\mu_B\disp_B + \mu_C\disp_C)/(\mu_B + \mu_C)$.\label{lemma:sum_hermite}
\end{lemma}

\begin{lemma}
    The binomially thinned version $D = q \circ E$ of $E \sim \text{Herm}(\mu_E, \disp_E)$ follows again a Hermite distribution with parameters $\lambda_D = q\mu_E, \disp_D = q\disp_E$.\label{lemma:thinned_hermite}
\end{lemma}

It is already known that the marginal distribution of a Hermite INAR(1) model is again Hermite \cite{Weiss2015b}. The following lemma provides the details in terms of our parameterization.
\begin{lemma}
\label{lemma:inar_hermite}
    The marginal distribution of a Hermite INAR(1) model
    $$
    X_t  = \alpha \circ X_{t - 1} + \import_t, \quad \import_t \stackrel{\textnormal{i.i.d.}}{\sim} \textnormal{Herm}(\tau, \disp)
    $$
    is Hermite. Specifically, it is given by
    \begin{equation} \label{eq:inar_hermite}
        X_t \sim \textnormal{Herm}\left(\frac{\tau}{1 - \alpha}, \frac{\disp}{1 + \alpha}\right),
    \end{equation}
    provided that we initialize $X_0$ with the same distribution.
\end{lemma}
We will show the statement by induction. Let us assume that \eqref{eq:inar_hermite} holds for $X_{t - 1}$. Then from Lemma \ref{lemma:thinned_hermite} we have
$$
 (\alpha \circ X_{t - 1}) \sim \text{Herm}\left(\frac{\alpha\tau}{1 - \alpha}, \frac{\alpha\disp}{1 + \alpha}\right),
$$
and $\import_t \sim \text{Herm}(\tau, \disp)$. According to Lemma \ref{lemma:sum_hermite}, $X_t$ is then likewise Hermite-distributed with mean
$$
\frac{\alpha\tau}{1 - \alpha} + \tau = \frac{\tau}{1 - \alpha},
$$
and dispersion parameter
$$
\left( \frac{\alpha\tau}{1 - \alpha} \times \frac{\alpha \disp}{1 + \alpha} + \disp\tau \right) \Big/ \left( \frac{\tau}{1 - \alpha} \right) = \frac{\alpha^2\tau \disp + \disp\tau - \disp\tau\alpha^2}{1 - \alpha^2} \times \frac{1 - \alpha}{\tau} = \frac{\disp}{1 + \alpha}.
$$
So we have shown that $X_t$ follows the Hermite distribution in question if $X_{t - 1}$ does. Since by definition we know that $X_0$ follows this distribution, we can conclude by induction that the statement is true for all $t \in \mathbb{N}_0$.

We now return to the case of the Hermite INARMA(1, 1) model. Lemmas \ref{lemma:sum_hermite}--\ref{lemma:inar_hermite} in combination with equation \eqref{eq:E_INAR} from Remark \ref{remark:inarma11} implies that $\{E_t\}$ is a Hermite INAR(1) model. More specifically, we have
\begin{align*}
E_t & = \xi \circ E_{t - 1} + \import^*_t\\
\import^*_t = \kappa \circ \import_{t - 1} & \sim \text{Herm}\left(\kappa\tau, \kappa\disp\right),
\end{align*}
where we used the shorthand $\xi$ as defined in equation \eqref{eq:define_xi}.
If we initiate the process as 
$$
E_0 \sim \text{Herm}\left(\frac{\kappa\tau}{1 - \xi}, \frac{\kappa\disp}{1 + \xi}\right),
$$
then $\{E_t\}$ is strictly stationary with the same marginal distribution. As $X_t = (1 - \beta) \circ E_t + \import_t$ where $(1 - \beta) \circ E_t \perp \import_t$, the observable process $\{X_t\}$ then likewise has Hermite marginals with mean
$$
\frac{(1 - \beta)\kappa\tau}{1 - \xi} + \tau = \frac{\kappa\tau}{1 - \kappa} + \tau = \frac{\tau}{1 - \kappa},
$$
and dispersion parameter
\begin{align*}
    \left(\frac{(1 - \beta)\kappa\tau}{1 - \xi} \times \frac{(1 - \beta)\kappa\disp}{1 + \xi} + \tau\disp\right) \Big/ \left( \frac{\tau}{1 - \kappa}\right) & = [...] = \frac{(1 - \beta)\kappa^2\disp + \disp(1 + \xi)(1 - \kappa)}{1 + \xi} \\
    & = \disp(1 - \kappa) + \kappa\frac{(1 - \beta) \kappa \disp}{1 + \xi}.
\end{align*}
This follows from Lemmas \ref{lemma:sum_hermite} and \ref{lemma:thinned_hermite} and concludes the proof.

\subsubsection{Consistency and normality of moment estimators for Poisson imports: Proposition \ref{proposition:moments}}

The proof closely follows the proof of Theorem 2.1 in \cite{Weiss2016}, which in turn builds upon the proofs of Theorem 4.1.1 and Lemma A.5.1 from \cite{Schweer2014}. Consider the vector-valued process
\begin{align}
\mathbf{Y}_t & := \big[X_t - \mu_X, X^2_t - \mu_X(0), X_t X_{t + 1} - \mu_X(1) , X_t X_{t + 2} - \mu_X(2)\big]^\top\\
& \text{ with } \mu_X(d):= \mathbb{E}(X_t X_{t + d}),\label{eq:def_Y}
\end{align}
which satisfies $\mathbb{E}(\mathbf{Y}_t) = \mathbf{0}$.
According to Proposition \ref{proposition:geometric_ergodicity}, the joint process $\{E_t, X_t\}$ is $\beta$-mixing with exponentially decreasing weights if $E_t$ is initialized with its stationary distribution (which in the Poisson INARMA(1, 1) is just a Poisson distribution; see Lemma \ref{proposition:bivariate_poisson}). It is thus also $\alpha$-mixing with exponentially decaying weights. We moreover note that $\{X_t\}$ is marginally Poisson such that all its moments are finite.

Since $\{\mathbf{Y}_t\}$ emerges from a measurable function of $\{X_t\}$ in \eqref{eq:def_Y}, it is also $\alpha$-mixing with exponentially decreasing weights, and it is straightforward to show that all its moments are finite. Thus, as in Weiss \cite[Section 2]{Weiss2016} and Schweer and Weiss \cite[Section 3.4 and proof of theorem 4.1.1]{Schweer2014}, Theorem 1.7 of Ibragimov \cite{Ibragimov1962} is applicable to the vector-valued process $\{\mathbf{Y}_t\}$ and implies
$$
\frac{1}{\sqrt{T}} \sum_{t = 1}^T \mathbf{Y}_t \stackrel{\mathcal{D}}{\longrightarrow} \text{N}(\mathbf{0}, \mathbf{\Sigma}).
$$
While the exact entries of $\mathbf{\Sigma}$ are finite and could in principle be determined along the lines of \cite{Weiss2016}, the computations get very involved and are thus not pursued here. It is sufficient for our purposes to conclude that
$$
\big[\hat{\mu}_X, \hat{\mu}_X(0), \hat{\mu}_X(1), \hat{\mu}_X(2)\big] \stackrel{\mathcal{D}}{\longrightarrow} \text{N}\left\{\big[\mu_X, \mu_X(0), \mu_X(1), \mu_X(2)\big], \frac{1}{T}\mathbf{\Sigma}\right\}
$$
where
$$
\hat{\mu}_X = \frac{1}{T} \sum_{t = 1}^T X_t,\ \ \ \hat{\mu}_X(0) = \frac{1}{T} \sum_{t = 1}^T X^2_t, \ \ \ \hat{\mu}_X(1) = \frac{1}{T} \sum_{t = 1}^T X_tX_{t + 1}, \ \ \ \hat{\mu}_X(2) = \frac{1}{T} \sum_{t = 1}^T X_tX_{t + 2}.
$$
We thus have a normally distributed estimator of a vector of relevant moments. Following the same arguments as in \cite{Weiss2016}, repeated application of the Delta method can be used to show that
$$
\big[\hat{\mu}_X, \hat{\rho}_X(1), \hat{\rho}_X(2)\big]^\top,
$$
where
$$
\hat{\rho}_X(d) = \frac{\hat{\mu}_X(d) + \hat{\mu}_X^2}{\hat{\mu}_X(0) + \hat{\mu}_X^2} \ \ \text{for} \ \ d = 1, 2
$$
and ultimately $(\hat{\tau}, \hat{\beta}, \hat{\kappa})$ are likewise asymptotically normally distributed. Note that to apply the delta method we need to assume that $\kappa > 0$ as otherwise the relevant derivatives do not exist.

Consistency of the estimators follows directly from the strict stationarity and geometric ergodicity of the process $\{E_t, X_t\}$ (which ensures that the relevant sample moments converge to their theoretical counterparts) and Slutsky's theorem.

\subsubsection{Moment estimators for general import distributions: Lemma \ref{lemma:kappa_cubic}}

In the following we will solve equations \eqref{eq:mu_X}--\eqref{eq:rho_X} for the model parameters $\tau, \sigma^2_\tau, \kappa$ and $\beta$. For the sake of better readability, we will use $\xi = \gamma_X(2)/\gamma_X(1)$ as defined in equation \eqref{eq:define_xi}. 

We start by solving \eqref{eq:mu_X} and \eqref{eq:define_xi} for $\tau$ and $\beta$, respectively, i.e.
\begin{align}
    \tau & = (1 - \kappa)\mu_X,\label{eq:tau}\\
    \beta & = \frac{\xi - \kappa}{1 - \kappa}.\label{eq:beta}
\end{align}
Now we plug \eqref{eq:tau} and \eqref{eq:beta} into \eqref{eq:sigma2_X},
\begin{align*}
    \sigma^2_X & \stackrel{\phantom{\eqref{eq:define_xi}}}{=} \frac{\kappa(1 + \beta)}{1 + \beta + (1 - \beta)\kappa} \times \frac{\tau}{1 - \kappa} + \left(1 - \frac{\kappa(1 + \beta)}{1 + \beta + (1 - \beta)\kappa}\right) \times \frac{\sigma^2_\tau}{1 - \kappa} \\
    & \stackrel{\eqref{eq:define_xi}}{=} \frac{\kappa(1 + \beta)}{1 + \xi} \times \mu_X + \left(1 - \frac{\kappa(1 + \beta)}{1 + \xi}\right) \times \frac{\sigma^2_\tau}{1 - \kappa} \\
    & \stackrel{\phantom{\eqref{eq:define_xi}}}{=} \frac{\kappa(1 + \frac{\xi - \kappa}{1 - \kappa})}{1 + \xi} \times \mu_X + \left(\frac{1 + \xi - \kappa(1 + \frac{\xi - \kappa}{1 - \kappa})}{1 + \xi}\right) \times \frac{\sigma^2_\tau}{1 - \kappa}\\
    & \stackrel{\phantom{\eqref{eq:define_xi}}}{=} \frac{\kappa(1 - 2\kappa + \xi)}{(1 + \xi)(1 - \kappa)} \times \mu_X + \frac{1 + \xi - \kappa\frac{1 - 2\kappa + \xi}{1 - \kappa}}{(1 + \xi)(1 - \kappa)} \times \sigma^2_{\tau} \\
    & \stackrel{\phantom{\eqref{eq:define_xi}}}{=} \frac{\kappa(1 - 2\kappa + \xi)}{(1 + \xi)(1 - \kappa)}\times \mu_X + \frac{1 + \xi - 2\kappa - 2\kappa\xi + 2\kappa^2}{(1 + \xi)(1 - \kappa)^2} \times \sigma^2_{\tau}
\end{align*}
\noindent and solve for $\sigma^2_\tau$,
\begin{align}
    \sigma^2_{\tau} = \frac{(1 + \xi)(1 - \kappa)^2}{1 + \xi - 2\kappa - 2\kappa\xi + 2\kappa^2} \times \sigma^2_X - \frac{\kappa(1 - \kappa)(1 - 2\kappa + \xi)}{1 + \xi - 2\kappa - 2\kappa\xi + 2\kappa^2} \times \mu_X.
    \label{eq:sigma2_tau}
\end{align}
Finally, we plug \eqref{eq:tau}, \eqref{eq:beta} and \eqref{eq:sigma2_tau} into the expression \eqref{eq:rho_X} for $\gamma_X(1)$,

\begin{align*}
    \gamma_X(1) & \stackrel{\phantom{\eqref{eq:beta}, \eqref{eq:tau}}}{=} (1 - \beta)\kappa \times\left(1 + \frac{\kappa\beta(\sigma^2_\tau - \tau)}{(1 + \beta) \{(1 - \kappa)\sigma^2_\tau + \kappa\tau\} + (1 - \beta)\kappa\sigma^2_\tau}\right) \times \sigma^2_X\\
    & \stackrel{\phantom{\eqref{eq:beta}, \eqref{eq:tau}}}{=} (1 - \beta)\kappa \times\left(1 + \kappa\beta\frac{\sigma^2_\tau - \tau}{(1 + \beta - 2\beta\kappa) \sigma^2_{\tau} + (1 + \beta)\kappa\tau}\right) \times \sigma^2_X\\
    & \stackrel{\eqref{eq:beta}, \eqref{eq:tau}}{=} \frac{1 - \xi}{1 - \kappa} \times \kappa\sigma^2_X \times \left( 1 + \kappa \frac{\xi - \kappa}{1 - \kappa} \times \frac{\sigma^2_{\tau} - (1 - \kappa)\mu_X}{\frac{1 + \xi - 2\kappa - 2\kappa\xi + 2\kappa^2}{1 - \kappa}\sigma^2_{\tau} + (1 + \xi - 2\kappa)\kappa\mu_X}\right)\\
    & \stackrel{\phantom{(10}\eqref{eq:sigma2_tau}\phantom{7,)}}{=} \frac{1 - \xi}{1 - \kappa}\times\kappa\sigma^2_X \times \left( 1 + \kappa \frac{\xi - \kappa}{1 - \kappa} \times \frac{\frac{(1 + \xi)(1 - \kappa)^2}{1 + \xi - 2\kappa - 2\kappa\xi + 2\kappa^2}\sigma^2_X - \frac{\kappa(1 - \kappa)(1 - 2\kappa + \xi)}{1 + \xi - 2\kappa - 2\kappa\xi + 2\kappa^2}\mu_X - (1 - \kappa)\mu_X}{(1 + \xi)(1 - \kappa)\sigma^2_X}\right)\\
    & \stackrel{\phantom{\eqref{eq:beta}, \eqref{eq:tau}}}{=} \frac{1 - \xi}{1 - \kappa}\times\kappa\sigma^2_X \times \left( 1 + \kappa (\xi - \kappa) \times \frac{(1 + \xi)(1 - \kappa)\sigma^2_X - (1 + \xi)(1 - \kappa)\mu_X}{(1 + \xi)(1 - \kappa)(1 + \xi - 2\kappa - 2\kappa\xi + 2\kappa^2)\sigma^2_X}\right)\\
    & \stackrel{\phantom{\eqref{eq:beta}, \eqref{eq:tau}}}{=} \frac{1 - \xi}{1 - \kappa}\times\kappa \times \frac{\{(1 + \xi - 2\kappa - 2\kappa\xi + 2\kappa^2) + \kappa(\xi - \kappa)\}\sigma^2_X - \kappa^2(\xi - \kappa)\mu_X}{1 + \xi - 2\kappa - 2\kappa\xi + 2\kappa^2}\\
    & \stackrel{\phantom{\eqref{eq:beta}, \eqref{eq:tau}}}{=} \frac{(1 - \xi) \{(1 + \xi)\kappa - (2 + \xi)\kappa^2 + \kappa^3\}\sigma^2_X - (1 - \xi)(\kappa^2\xi - \kappa^3) \mu_X}{1 + \xi - 3(1 + \xi)\kappa + 2(2 + \xi)\kappa^2 - 2\kappa^3}.
\end{align*}

\noindent To obtain the cubic equation \eqref{eq:abcd}, we multiply both sides of the above equation by $(1 + \xi - 3(1 + \xi)\kappa + 2(2 + \xi)\kappa^2 - 2\kappa^3)$, move all terms to the right hand side and order them by power of $\kappa$.

\newpage

\section{Supplementary materials on real-data example}


\begin{figure}[h]
\includegraphics[scale=0.8]{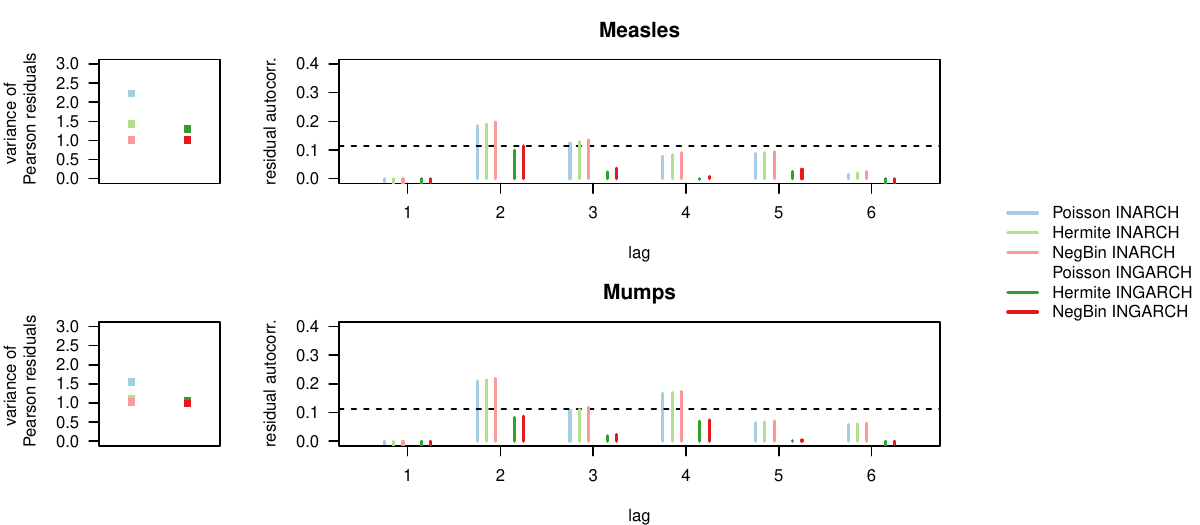}
\caption{Analysis of Pearson residuals of INARCH and INGARCH models with different immigration distributions. Left: Variance of Pearson residuals, which should ideally be close to 1. Right: Empirical autocorrelation function of Pearson residuals. Ideally, these should be close to 0, below the dashed line. The dashed line shows $2/\sqrt{T}$, which is the 97.5\% quantile for the empirical ACF of white noise.}
\label{fig:residuals_ingarch}
\end{figure}

\newpage

\begin{table}[h!]
\footnotesize
\caption{Summary of model fits to the measles and mumps data with parameters on the scale which is used for estimation.}
\label{tab:fits_real_data}
\center
\medskip
\begin{tabular}{p{3.5cm} @{\hskip 0.8cm} p{0.45cm} p{0.45cm} p{0.45cm} p{0.45cm} p{0.45cm} @{\hskip 2cm} p{0.45cm} p{0.45cm} p{0.45cm} p{0.45cm} p{1cm}}
\hline
 & & \multicolumn{3}{c}{Measles} & & \multicolumn{5}{c}{Mumps}\\
\hline 
Model & $\hat\nu$ & $\hat\alpha$ & $\hat{\beta}$ & $\hat\disp^*$ & AIC & $\hat\nu$ & $\hat\alpha$ & $\hat{\beta}$ & $\hat\disp^*$ & AIC \\
\hline\\
\input{table/tab_inarch.tex}\\
\input{table/tab_ingarch.tex}\\
\hline 
Model & $\hat\nu$ & $\hat\alpha$ & & $\hat\disp$ & AIC & $\hat\nu$ & $\hat\alpha$ &  & $\hat\disp$ & AIC \\
\hline\\
\input{table/tab_inar.tex}\\
\hline 
Model & $\hat\tau$ & $\hat\kappa$ & $\hat{\beta}$ & $\hat\disp$ & AIC & $\hat\tau$ & $\hat\kappa$ & $\hat{\beta}$ & $\hat\disp$ & AIC \\
\hline\\
\input{table/tab_inarma.tex}\\
\hline
\end{tabular}
\end{table}

\begin{figure}
    \centering
    \includegraphics[width=0.97\textwidth]{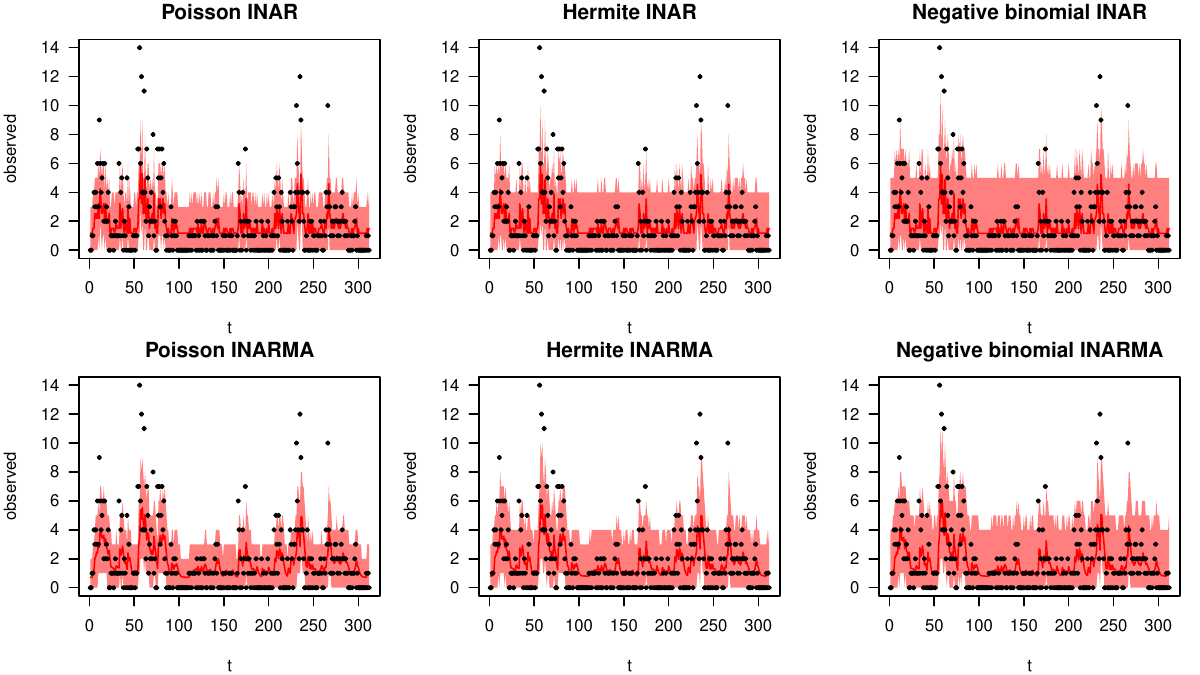}
    \caption{Fits of INAR and INARMA models to the measles data set. Red lines show fitted values, light red shaded areas 90\% in-samples prediction intervals.}
    \label{fig:fit_measles}
\end{figure}

\begin{figure}
    \centering
    \includegraphics[width=0.97\textwidth]{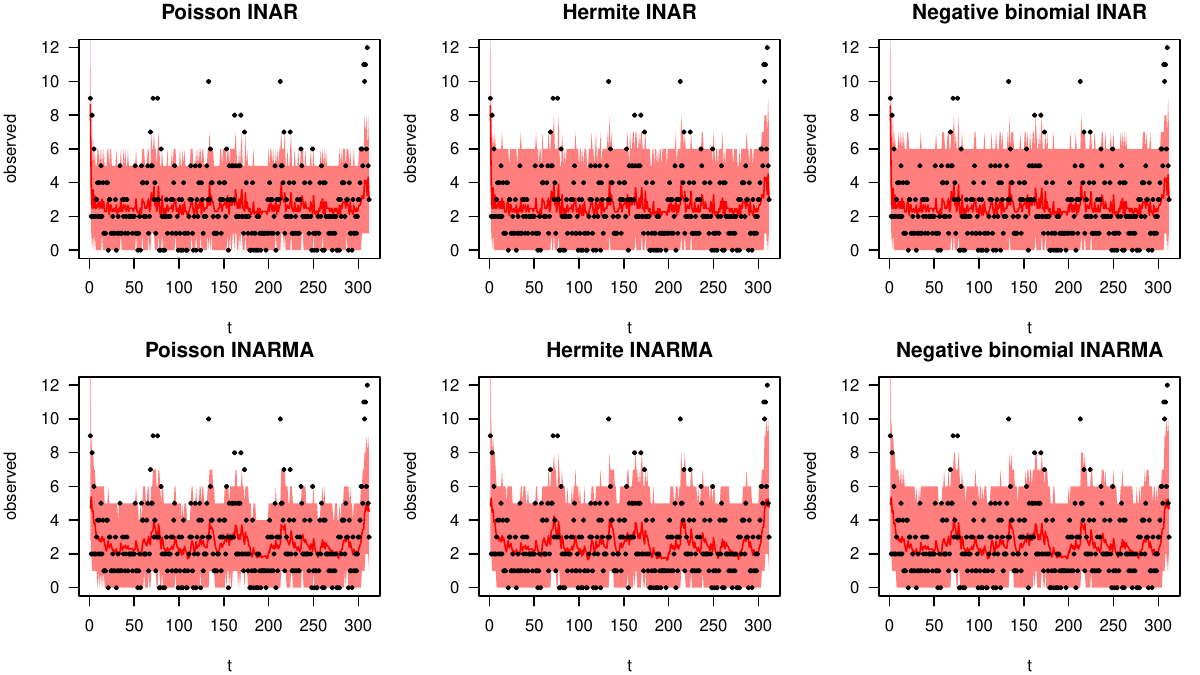}
    \caption{Fits of INAR and INARMA models to the mumps data set. Red lines show fitted values, light red shaded areas 90\% in-samples prediction intervals.}
    \label{fig:fit_mumps}
\end{figure}

\end{document}

%% file: table/sim_pois_sc1.tex
250 & 1.000 & 1.044 & 0.266 & 0.245 & - & - & - & - & 0.500 & 0.449 & 0.167 & 0.130 & 0.500 & 0.477 & 0.132 & 0.111 \\ 
  500 &  & 1.025 & 0.186 & 0.179 &  & - & - & - &  & 0.476 & 0.113 & 0.098 &  & 0.486 & 0.091 & 0.084 \\ 
  1000 &  & 1.015 & 0.131 & 0.128 &  & - & - & - &  & 0.486 & 0.073 & 0.070 &  & 0.493 & 0.064 & 0.061 \\

%% file: table/sim_pois_sc2.tex
 250 & 1.000 & 1.008 & 0.186 & 0.185 & - & - & - & - & 0.200 & 0.190 & 0.091 & 0.089 & 0.600 & 0.596 & 0.072 & 0.071 \\ 
  500 &  & 1.011 & 0.132 & 0.132 &  & - & - & - &  & 0.196 & 0.066 & 0.065 &  & 0.595 & 0.052 & 0.050 \\ 
  1000  &  & 1.005 & 0.092 & 0.092 &  & - & - & - &  & 0.198 & 0.048 & 0.046 &  & 0.598 & 0.035 & 0.035 \\

%% file: table/sim_pois_sc3.tex
 250& 1.000 & 1.009 & 0.182 & 0.184 & - & - & - & - & 0.100 & 0.097 & 0.039 & 0.040 & 0.800 & 0.797 & 0.036 & 0.035 \\ 
  500 &  & 1.007 & 0.123 & 0.128 &  & - & - & - &  & 0.099 & 0.028 & 0.028 &  & 0.798 & 0.024 & 0.025 \\ 
  1000 &  & 1.010 & 0.089 & 0.090 &  & - & - & - &  & 0.098 & 0.020 & 0.020 &  & 0.798 & 0.017 & 0.018 \\
  

%% file: table/sim_herm_sc1.tex
250 & 1.000 & 1.045 & 0.267 & 0.217 & 0.500 & 0.483 & 0.225 & 0.170 & 0.500 & 0.465 & 0.171 & 0.123 & 0.500 & 0.477 & 0.127 & 0.094 \\ 
  500 &  & 1.017 & 0.187 & 0.157 &  & 0.494 & 0.171 & 0.132 &  & 0.480 & 0.112 & 0.092 &  & 0.491 & 0.092 & 0.071 \\ 
  1000 &  & 1.022 & 0.129 & 0.122 &  & 0.488 & 0.113 & 0.105 &  & 0.491 & 0.075 & 0.069 &  & 0.489 & 0.062 & 0.057 \\

%% file: table/sim_herm_sc2.tex
250 & 1.000 & 1.032 & 0.189 & 0.178 & 0.700 & 0.656 & 0.195 & 0.160 & 0.200 & 0.191 & 0.092 & 0.086 & 0.600 & 0.586 & 0.069 & 0.062 \\ 
  500 &  & 1.019 & 0.125 & 0.123 &  & 0.679 & 0.128 & 0.117 &  & 0.193 & 0.063 & 0.063 &  & 0.591 & 0.046 & 0.044 \\ 
  1000 &  & 1.009 & 0.085 & 0.085 &  & 0.690 & 0.088 & 0.083 &  & 0.196 & 0.044 & 0.044 &  & 0.596 & 0.031 & 0.031 \\

%% file: table/sim_herm_sc3.tex
250 & 1.000 & 1.051 & 0.201 & 0.198 & 0.900 & 0.814 & 0.224 & 0.154 & 0.100 & 0.095 & 0.041 & 0.040 & 0.800 & 0.789 & 0.036 & 0.033 \\ 
  500 &  & 1.019 & 0.127 & 0.127 &  & 0.863 & 0.145 & 0.110 &  & 0.098 & 0.027 & 0.028 &  & 0.796 & 0.023 & 0.022 \\ 
  1000 &  & 1.008 & 0.087 & 0.086 &  & 0.884 & 0.099 & 0.081 &  & 0.101 & 0.020 & 0.020 &  & 0.798 & 0.015 & 0.016 \\ 
  

%% file: table/sim_negbin_sc1.tex
 250 & 1.000 & 1.017 & 0.255 & 0.208 & 0.500 & 0.599 & 0.668 & 0.549 & 0.500 & 0.474 & 0.168 & 0.122 & 0.500 & 0.491 & 0.121 & 0.093 \\ 
  500 &  & 1.013 & 0.173 & 0.159 &  & 0.532 & 0.269 & 0.136 &  & 0.485 & 0.109 & 0.093 &  & 0.494 & 0.084 & 0.073 \\ 
  1000 &  & 1.007 & 0.118 & 0.116 &  & 0.515 & 0.186 & 0.092 &  & 0.491 & 0.072 & 0.068 &  & 0.497 & 0.058 & 0.055 \\

%% file: table/sim_negbin_sc2.tex
250 & 1.000 & 1.019 & 0.177 & 0.179 & 0.700 & 0.713 & 0.379 & 0.309 & 0.200 & 0.191 & 0.091 & 0.088 & 0.600 & 0.590 & 0.067 & 0.064 \\ 
  500 &  & 1.012 & 0.122 & 0.125 &  & 0.708 & 0.258 & 0.192 &  & 0.196 & 0.063 & 0.064 &  & 0.595 & 0.047 & 0.045 \\ 
  1000 &  & 1.005 & 0.090 & 0.087 &  & 0.706 & 0.181 & 0.127 &  & 0.199 & 0.044 & 0.044 &  & 0.598 & 0.033 & 0.032 \\

%% file: table/sim_negbin_sc3.tex
250 & 1.000 & 1.024 & 0.190 & 0.192 & 0.900 & 0.911 & 0.469 & 0.501 & 0.100 & 0.099 & 0.043 & 0.041 & 0.800 & 0.795 & 0.035 & 0.034 \\ 
  500 &  & 1.019 & 0.128 & 0.131 &  & 0.895 & 0.312 & 0.299 &  & 0.097 & 0.029 & 0.029 &  & 0.796 & 0.024 & 0.024 \\ 
  1000 &  & 1.007 & 0.092 & 0.090 &  & 0.901 & 0.221 & 0.204 &  & 0.099 & 0.021 & 0.020 &  & 0.798 & 0.017 & 0.016 \\

%% file: table/tab_inarch_epi.tex
 Poisson INARCH & 0.83 & 0.54 & 1$^*$ & 1$^*$ & 1159.13 & 1.93 & 0.26 & 1$^*$ & 1$^*$ & 1274.26 \\ 
  Hermite INARCH & 0.85 & 0.52 & 1$^*$ & 1.37 & 1082.40 & 1.95 & 0.25 & 1$^*$ & 1.24 & 1249.33 \\ 
  NegBin INARCH & 0.88 & 0.51 & 1$^*$ & 1.51 & 1055.04 & 1.98 & 0.24 & 1$^*$ & 1.24 & 1244.75 \\ 
  

%% file: table/tab_ingarch_epi.tex
 Poisson INGARCH & 0.46 & 0.74 & 2.08 & 1$^*$ & 1096.91 & 1.04 & 0.60 & 2.98 & 1$^*$ & 1238.27 \\ 
  Hermite INGARCH & 0.51 & 0.72 & 2.01 & 1.31 & 1046.09 & 1.07 & 0.58 & 2.98 & 1.18 & 1224.43 \\ 
  NegBin INGARCH & 0.55 & 0.67 & 2.03 & 1.41 & 1028.23 & 1.11 & 0.57 & 2.96 & 1.19 & 1222.86 \\ 
  

%% file: table/tab_inar_epi.tex
 Poisson INAR & 1.17 & 0.34 & 1$^*$ & 1$^*$ & 1232.94 & 2.12 & 0.18 & 1$^*$ & 1$^*$ & 1283.22 \\ 
  Hermite INAR & 1.18 & 0.34 & 1$^*$ & 1$^*$ & 1122.68 & 2.07 & 0.20 & 1$^*$ & 1$^*$ & 1252.77 \\ 
  NegBin INAR & 1.17 & 0.34 & 1$^*$ & 1$^*$ & 1068.77 & 2.08 & 0.20 & 1$^*$ & 1$^*$ & 1245.57 \\ 
  

%% file: table/tab_inarma_epi.tex
 Poisson INARMA & 0.72 & 0.60 & 2.00 & 1$^*$ & 1166.26 & 1.38 & 0.47 & 2.50 & 1$^*$ & 1257.34 \\ 
  Hermite INARMA & 0.81 & 0.55 & 1.86 & 1$^*$ & 1094.07 & 1.30 & 0.50 & 2.51 & 1.00 & 1235.48 \\ 
  NegBin INARMA & 0.81 & 0.53 & 1.81 & 1$^*$ & 1046.65 & 1.41 & 0.46 & 2.46 & 1$^*$ & 1231.73 \\

%% file: table/sim_pois_sc1_moments.tex
250 & 1.000 & 0.988 & 0.400 & - & - & - & 0.500 & 0.457 & 0.236 & 0.500 & 0.506 & 0.198 \\ 
  500 &  & 0.960 & 0.343 &  & - & - &  & 0.482 & 0.182 &  & 0.519 & 0.171 \\ 
  1000 &  & 0.961 & 0.265 &  & - & - &  & 0.495 & 0.137 &  & 0.520 & 0.133 \\

%% file: table/sim_pois_sc2_moments.tex
250 & 1.000 & 1.029 & 0.246 & - & - & - & 0.200 & 0.198 & 0.119 & 0.600 & 0.588 & 0.096 \\ 
  500 &  & 1.020 & 0.182 &  & - & - &  & 0.199 & 0.083 &  & 0.592 & 0.070 \\ 
  1000 &  & 1.016 & 0.124 &  & - & - &  & 0.197 & 0.058 &  & 0.594 & 0.048 \\

%% file: table/sim_pois_sc3_moments.tex
250 & 1.000 & 1.104 & 0.298 & - & - & - & 0.100 & 0.102 & 0.053 & 0.800 & 0.778 & 0.060 \\ 
  500 &  & 1.051 & 0.200 &  & - & - &  & 0.101 & 0.035 &  & 0.790 & 0.039 \\ 
  1000 &  & 1.034 & 0.148 &  & - & - &  & 0.099 & 0.025 &  & 0.793 & 0.029 \\ 
  

%% file: table/sim_herm_sc1_moments.tex
250 & 1.000 & 0.963 & 0.379 & 0.500 & 0.571 & 0.357 & 0.500 & 0.483 & 0.243 & 0.500 & 0.518 & 0.185 \\ 
  500 &  & 0.945 & 0.317 &  & 0.573 & 0.293 &  & 0.494 & 0.186 &  & 0.527 & 0.157 \\ 
  1000 &  & 0.982 & 0.246 &  & 0.532 & 0.199 &  & 0.496 & 0.138 &  & 0.509 & 0.120 \\

%% file: table/sim_herm_sc2_moments.tex
250 & 1.000 & 1.037 & 0.242 & 0.700 & 0.676 & 0.326 & 0.200 & 0.198 & 0.122 & 0.600 & 0.585 & 0.092 \\ 
  500 &  & 1.025 & 0.163 &  & 0.684 & 0.227 &  & 0.197 & 0.084 &  & 0.589 & 0.063 \\ 
  1000 &  & 1.017 & 0.112 &  & 0.688 & 0.164 &  & 0.196 & 0.059 &  & 0.593 & 0.044 \\

%% file: table/sim_herm_sc3_moments.tex
250 & 1.000 & 1.104 & 0.284 & 0.900 & 0.819 & 0.526 & 0.100 & 0.098 & 0.053 & 0.800 & 0.778 & 0.054 \\ 
  500&  & 1.045 & 0.191 &  & 0.862 & 0.377 &  & 0.100 & 0.038 &  & 0.791 & 0.037 \\ 
  1000 &  & 1.028 & 0.137 &  & 0.878 & 0.276 &  & 0.100 & 0.027 &  & 0.795 & 0.026 \\
  

%% file: table/sim_negbin_sc1_moments.tex
250 & 1.000 & 0.962 & 0.382 & 0.500 & 0.901 & 1.034 & 0.500 & 0.476 & 0.239 & 0.500 & 0.520 & 0.184 \\ 
  500 &  & 0.953 & 0.311 &  & 0.801 & 0.775 &  & 0.492 & 0.186 &  & 0.524 & 0.154 \\ 
  1000 &  & 0.967 & 0.247 &  & 0.702 & 0.620 &  & 0.495 & 0.142 &  & 0.517 & 0.123 \\

%% file: table/sim_negbin_sc2_moments.tex
250 & 1.000 & 1.028 & 0.237 & 0.700 & 0.774 & 0.651 & 0.200 & 0.200 & 0.122 & 0.600 & 0.588 & 0.092 \\ 
  500 &  & 1.024 & 0.167 &  & 0.722 & 0.372 &  & 0.196 & 0.088 &  & 0.591 & 0.065 \\ 
  1000 &  & 1.008 & 0.117 &  & 0.718 & 0.252 &  & 0.201 & 0.060 &  & 0.597 & 0.045 \\

%% file: table/sim_negbin_sc3_moments.tex
250 & 1.000 & 1.103 & 0.286 & 0.900 & 0.876 & 0.830 & 0.100 & 0.099 & 0.056 & 0.800 & 0.780 & 0.055 \\ 
  500 &  & 1.056 & 0.195 &  & 0.892 & 0.551 &  & 0.097 & 0.040 &  & 0.789 & 0.038 \\ 
  1000 &  & 1.027 & 0.140 &  & 0.895 & 0.380 &  & 0.098 & 0.028 &  & 0.794 & 0.027 \\

%% file: table/tab_inarch.tex
 Poisson INARCH & 0.83 & 0.54 &  &  & 1159.13 & 1.93 & 0.26 &  &  & 1274.26 \\ 
  Hermite INARCH & 0.85 & 0.52 &  & 0.54 & 1082.40 & 1.95 & 0.25 &  & 0.39 & 1249.33 \\ 
  NegBin INARCH & 0.88 & 0.51 &  & 1.17 & 1055.04 & 1.98 & 0.24 &  & 0.52 & 1244.75 \\ 
  

%% file: table/tab_ingarch.tex
 Poisson INGARCH & 0.22 & 0.36 & 0.52 &  & 1096.91 & 0.35 & 0.20 & 0.66 &  & 1238.27 \\ 
  Hermite INGARCH & 0.25 & 0.36 & 0.50 & 0.47 & 1046.09 & 0.36 & 0.20 & 0.66 & 0.31 & 1224.43 \\ 
  NegBin INGARCH & 0.27 & 0.33 & 0.51 & 0.92 & 1028.23 & 0.38 & 0.19 & 0.66 & 0.40 & 1222.86 \\ 
  

%% file: table/tab_inar.tex
 INAR & 1.17 & 0.34 &  &  & 1232.94 & 2.12 & 0.18 &  &  & 1283.22 \\ 
  HINAR & 1.18 & 0.34 &  & 0.68 & 1122.68 & 2.07 & 0.20 &  & 0.50 & 1252.77 \\ 
  NBINAR & 1.17 & 0.34 &  & 1.81 & 1068.77 & 2.08 & 0.20 &  & 0.35 & 1245.57 \\ 
  

%% file: table/tab_inarma.tex
 Poisson INARMA & 0.72 & 0.60 & 0.50 &  & 1166.26 & 1.38 & 0.47 & 0.60 &  & 1257.34 \\ 
  Hermite INARMA & 0.81 & 0.55 & 0.46 & 0.74 & 1094.07 & 1.30 & 0.50 &  0.60 & 0.64 & 1235.48 \\ 
  NegBin INARMA & 0.81 & 0.53 & 0.45 & 3.17 & 1046.65 & 1.41 & 0.46 &  0.59 & 0.63 & 1231.73 \\ 
  